\documentclass[10pt]{article}

\usepackage[leqno]{amsmath}
\usepackage{amssymb}
\usepackage{latexsym}
\usepackage{amsthm}

\renewcommand{\theenumi}{\roman{enumi}}
\renewcommand{\labelenumi}{(\theenumi)}

\newcommand{\init}{\big\vert_{t = 0}}
 \newcommand{\half}{\frac{1}{2}}
\newcommand{\abs}[1]{\left\vert #1 \right\vert}
\newcommand{\bigabs}[1]{\bigl\vert #1 \bigr\vert}
\newcommand{\Bigabs}[1]{\Bigl\vert #1 \Bigr\vert}

\newcommand{\norm}[1]{\left\Vert #1 \right\Vert}
\newcommand{\bignorm}[1]{\bigl\Vert #1 \bigr\Vert}

\newcommand{\spacetimenorm}[3]{\norm{#1}_{#2,#3}}

\newcommand{\Spacetimenorm}[3]{\abs{#1}_{#2,#3}}

\newcommand{\Sobnorm}[2]{\norm{#1}_{H^{#2}}}
\newcommand{\Sobdotnorm}[2]{\norm{#1}_{\dot{H}^{#2}}}

\newcommand{\twonorm}[2]{\norm{#1}_{L^2#2}}
\newcommand{\bigtwonorm}[2]{\bignorm{#1}_{L^2#2}}
\newcommand{\inftynorm}[2]{\norm{#1}_{L^\infty#2}}

\newcommand{\mixednorm}[3]{\norm{#1}_{L_{t}^{#2}(L_{x}^{#3})}}
\newcommand{\mixed}[2]{L_{t}^{#1}(L_{x}^{#2})}
\newcommand{\Mixednorm}[3]{\norm{#1}_{\LL_{t}^{#2}(\LL_{x}^{#3})}}
\newcommand{\Mixed}[2]{\LL_{t}^{#1}(\LL_{x}^{#2})}
\newcommand{\bigMixednorm}[3]{\bignorm{#1}_{\LL_{t}^{#2}(\LL_{x}^{#3})}}
\newcommand{\bigmixednorm}[3]{\bignorm{#1}_{L_{t}^{#2}(L_{x}^{#3})}}

\newcommand{\C}{\mathbb{C}}

\newcommand{\scrH}{\mathcal{H}}
\newcommand{\LL}{\mathcal{L}}

\newcommand{\N}{\mathbb{N}}

\newcommand{\R}{\mathbb{R}}
\newcommand{\X}{\mathcal{X}}
\newcommand{\Y}{\mathcal{Y}}

\newcommand{\Schwartz}{\mathcal{S}}
\newcommand{\Fourier}{\mathcal{F}}

\newcommand{\innerprod}[2]{\left\langle \, #1 , #2 \, \right\rangle}
\newcommand{\biginnerprod}[2]{\bigl\langle \, #1 , #2 \, \bigr\rangle}

\newcommand{\hypwt}[2]{\bigabs{ \abs{#1} - \abs{#2} }}
\newcommand{\elwt}[1]{\left\langle #1 \right\rangle}

\DeclareMathOperator{\diag}{diag}

\DeclareMathOperator{\supp}{supp}

\swapnumbers
\newtheorem{theorem}{Theorem}[section]
\newtheorem{proposition}[theorem]{Proposition}
\newtheorem{lemma}[theorem]{Lemma}
\newtheorem{corollary}[theorem]{Corollary}
\newtheorem{definition}[theorem]{Definition}
\newtheorem{principle}[theorem]{Principle}
\newtheorem*{ClassicalLocalEx}{Classical Local Existence Theorem}
\newtheorem*{ThmA}{Theorem A}
\newtheorem*{ThmB}{Theorem B}
\newtheorem*{ThmC}{Theorem C}
\newtheorem*{ThmD}{Theorem D}
\newtheorem*{ThmE}{Theorem E}
\newtheorem*{ThmF}{Theorem F}
\newtheorem*{plaintheorem}{Theorem}
\newtheorem*{plainproposition}{Proposition} \newtheorem*{MainTheorem}{Main
Theorem}
\newtheorem*{GeneralWPConjecture}{General WP Conjecture}

\theoremstyle{definition}
\newtheorem{example}[theorem]{Example}

\theoremstyle{remark}
\newtheorem{remark}[theorem]{Remark}
\newtheorem*{remarkNoLabel}{Remark}

\title{Bilinear estimates and applications to nonlinear wave equations}
\author{Sergiu Klainerman and Sigmund Selberg} \date{}
\begin{document}

\numberwithin{equation}{section}

\setcounter{secnumdepth}{2}

\maketitle
\begin{abstract}
We undertake a systematic review of results proved in
\cite{Kl-Ma1,Kl-Ma3,Kl-Ma5,Kl-Se,Kl-Ta} concerning local well-posedness of
the Cauchy problem for certain systems of nonlinear wave equations, with
minimal regularity assumptions on the initial data. Moreover we give a
considerably simplified and unified treatment
of these results and provide also complete proofs for large data.
The paper is also intended as an introduction to and survey of current research
in the very active area of nonlinear wave equations. The key
ingredients throughout the survey are the use of the null structure
of the equations we consider and, intimately tied to it, bilinear
estimates.
\end{abstract}
\section{Introduction}
In this paper we undertake a systematic
review of results proved in
\cite{Kl-Ma1,Kl-Ma3,Kl-Ma5,Kl-Se,Kl-Ta} concerning local well-posedness of
the Cauchy problem for certain systems of nonlinear wave equations, with
minimal regularity assumptions on the initial data. Moreover we give a
vastly simplified and unified treatment
of these results and provide also complete proofs for large data. The key
ingredient throughout the survey is the use of
space-time bilinear estimates; they are intimately tied to the null
structure of the equations we consider. The simplest type of bilinear
estimates are of $L^2$ type; they transfom, by Plancherel's identity, to
bilinear $L^2$ convolution estimates in Fourier space. This leads
naturally to weighted $ L^2 $ spaces; in view of their similarity to the
Sobolev $H^s$ spaces we denote them $H^{s,\theta}$, and propose to call
them
Wave-Sobolev spaces\footnote{These spaces have apppeared before in PDE, in
connection with questions of propagation of singularities for nonlinear
wave equations. In the context of bilinear estimates and optimal
well-posedness of the Cauchy problem they appear first in \cite{B}, in the
study of periodic solution to KdV and nonlinear Schr\"odinger equations (see
also \cite{KPV}), and in \cite{Kl-Ma0.1} in connection with semilinear wave
equations satisfying the null condition. See section \ref{FurtherResults}
for more complete
historical remarks.}. Though these spaces play a fundamental role, in most
applications they need to be refined. In this survey we do this by taking
their
intersection with suitable weighted $\Mixed{q}{r}$ type spaces. These
spaces are described in detail in section \ref{MixedSpace}.
The main bilinear estimates
are summarized in section \ref{WaveEstimates}, the spaces $H^{s,\theta}$
are discussed in
section \ref{WaveSobolevSpaces}. The main nonlinear results are
stated below and proved in
sections \ref{WMProof}--\ref{WMMProof}. In section \ref{FurtherResults}
we discuss some of the main open problems and
provide some historical remarks.

On the Minkowski space-time $\R \times \R^{n} \simeq \R^{1+n}$ we use
coordinates $(t,x) = (x^{0},\dots,x^{n})$, and indices are raised and
lowered relative to the metric $m_{\mu\nu} = \diag(-1,1,\dots,1)$. The
summation convention is used in some sections. We write $\partial_{\mu} =
\partial_{x^{\mu}}$ and $\partial_{t} = \partial_{0}$.

We are interested in the Cauchy problem for systems of the type $$
\square u = \mathcal N(u),
$$
where $\square =
-\partial_{t}^{2} + \Delta$ is the standard wave operator on $\R \times
\R^{n}$, $\Delta = \sum_{1}^{n} \partial_{j}^{2}$ is the Laplacian on
$\R^{n}$, $u = u(t,x)$ takes values in $\R^{N}$ for some $N \ge 1$ and
$\mathcal N$ is an operator which (i) is local in time, in the sense
that for any open interval $I \subseteq \R$, the values of
$\mathcal N(u)$ on $I \times \R^n$ only depends on the values of $u$
on the same region; (ii) is time-translation invariant, in the sense
that $\mathcal N\bigl( u(\cdot + t, \cdot) \bigr) = \mathcal N(u)
(\cdot + t, \cdot)$ for all $t \in \R$; and (iii) satisfies $\mathcal N(0) = 0$.

Cauchy data are prescribed on the initial hypersurface $\{0\} \times
\R^{n} \simeq \R^{n}$:
$$
(u,\partial_{t} u) \init = (f,g) \in H^{s} \times H^{s-1} $$
where $H^s = \{ f : (I - \Delta)^{s/2} f \in L^2 \}$.
\subsection{Statement of Main Results}\label{BasicEquations}
We shall in fact concentrate on systems of the following types:
\begin{enumerate}
\item
\textit{Wave Maps Type:}
\begin{equation}\tag{WM}\label{WMType}
\square u^{I} + \sum_{J,K} \Gamma^{I}_{JK}(u) Q_{0}(u^{J},u^{K}) = 0.
\end{equation}
Here, $u^{I}$ denotes the $I$-th component function of $u$, the
$\Gamma^{I}_{JK}$ are smooth functions from $\R^{N}$ into $\R$ and $Q_{0}$
is the null form
$$
Q_{0}(\phi,\psi) = \sum_{\mu = 0}^{n} \partial_{\mu} \phi \partial^{\mu}
\psi = -\partial_{t} \phi \partial_{t} \psi + \sum_{j = 1}^{n}
\partial_{j} \phi \partial_{j} \psi. $$
\item
\textit{Yang-Mills Type:}
\begin{equation}\tag{``YM''}\label{YMType}
\square u = D^{-1} Q(u,u) + Q(D^{-1}u,u), \end{equation}
where $D^{\alpha} = (- \Delta)^{\alpha/2}$ and $Q$ stands for any bilinear
operator of the following type: Given vector-valued functions $u$ and $v$,
the $I$-th component function of $Q(u,v)$ is a linear combination, with
constant, real coefficients, of $Q_{ij}(u^{J},v^{K})$ for all $1 \le i < j
\le n$ and all $J,K$, where $Q_{ij}$ is the null form
$$
Q_{ij}(\phi,\psi) = \partial_{i} \phi \partial_{j} \psi - \partial_{i}
\phi \partial_{j} \psi.
$$
(The two $Q$'s on the right hand side of \eqref{YMType} may represent two
different such operators.) %
\item
\textit{Maxwell-Klein-Gordon Type:}
\begin{equation}\tag{``MKG''}\label{MKGType}
\begin{cases}
\square u = D^{-1} Q(v,v), \\
\square v = Q(D^{-1}u,v),
\end{cases}
\end{equation}
where $u = (u^{1},\dots,u^{N_{1}})$, $v = (v^{1},\dots,v^{N_{2}})$, $N =
N_{1} + N_{2}$ and $Q$ has the same meaning as before. Thus
\eqref{MKGType} is a special case of \eqref{YMType}. %
\item
\textit{Wave Maps Model Problem:}
\begin{equation}\tag{WMM}\label{CFWMType}
\square u^{I} = \sum_{J,K = 1}^{N} a^{I}_{JK} \widetilde Q(u^{J},u^{K}),
\end{equation}
where the $a^{I}_{JK}$ are real constants, $$
\widetilde Q(\phi,\psi) = \sum_{j = 1}^{n} \partial_{j} \bigl( R_{0} R_{j}
\phi \cdot \psi - \phi \cdot R_{0} R_{j} \psi \bigr) $$
and $R_{\mu} = D^{-1} \partial_{\mu}$.
\end{enumerate}
The following theorem summarizes the main well-posedness results
proved\footnote{Strictly speaking most of these results were proved only
for sufficiently small data. Large data recquire some technical
considerations discussed in this paper.} in
\cite{Kl-Ma1,Kl-Ma3,Kl-Ma5,Kl-Se,Kl-Ta}. %
\begin{MainTheorem}
{
\renewcommand{\theenumi}{\alph{enumi}}
\renewcommand{\labelenumi}{(\theenumi)}
\begin{enumerate}
\item\label{WMResult}
(\cite{Kl-Ma1, Kl-Se}.) If $n \ge 2$ and $s > \frac{n}{2}$, then
\eqref{WMType} is locally well-posed for initial data in $H^{s} \times
H^{s-1}$.
\item
(\cite{Kl-Ma5, Kl-Ta}.) If $n \ge 4$ and $s > \frac{n-2}{2}$, then
\eqref{MKGType} and \eqref{YMType} are locally well-posed for initial data
in $H^{s} \times H^{s-1}$.
\item
(\cite{Kl-Ma3}.) If $n \ge 3$ and
$s > \frac{n-2}{2}$, then \eqref{CFWMType} is locally well-posed for
initial data in $H^{s} \times H^{s-1}$.
\end{enumerate}
}
\end{MainTheorem}
By \emph{locally well-posed} we mean that for all $(f,g) \in H^{s} \times
H^{s-1}$ there exist $T > 0$ and
$$
u \in C([0,T],H^{s}) \cap C^{1}([0,T],H^{s-1}) $$
such that $u$ solves the equation on $(0,T) \times \R^{n}$ in the sense of
distributions, and such that the initial condition is satisfied. Moreover,
$T$ is bounded below by a strictly positive and continuous function of
$\Sobnorm{f}{s} + \Sobnorm{g}{s-1}$, the map $(f,g) \mapsto u$ is locally
Lipschitz\footnote{In fact, the solution depends smoothly (or even
analytically in most of the above examples) on the data,
in the sense that if
$\varepsilon
\mapsto (f_{\varepsilon},g_{\varepsilon})$ is a smooth map into $H^{s}
\times H^{s-1}$ for $\abs{\varepsilon} < \varepsilon_{0}$, and if
$u_{\varepsilon}$ is the solution corresponding to the initial data
$(f_{\varepsilon},g_{\varepsilon})$, then $\varepsilon \mapsto
u_{\varepsilon}$ is a smooth map into $C([0,T],H^{s}) \cap
C^{1}([0,T],H^{s-1})$ for some $T > 0$.
This is because the solution is obtained by a Picard iteration procedure;
see \cite{Se2}.}, and $u$ is unique in some subspace of $C([0,T],H^{s})
\cap C^{1}([0,T],H^{s-1})$. Moreover, any additional regularity of the
initial data persists in time, but for simplicity we ignore this issue.
\subsection{Motivation of the Equations} %
With the exception of (\ref{WMType}), the equations we work with are model
problems derived from the actual Maxwell-Klein-Gordon, Yang-Mills and wave
maps equations. Here we review these equations and discuss how our model
problems relate to them.
\subsubsection{Wave Maps}
A \emph{wave map} from the Minkowski space-time into a Riemannian manifold
$(M,g)$ is a map $u : \R^{1+n} \to M$ which is a critical point with
respect to compactly supported variations of the Lagrangian $$
\mathcal L[u] = \half \int_{\R^{1+n}} \innerprod{du}{du} \, dt \, dx, $$
where $\innerprod{du}{du} = \sum_{\mu = 0}^{n} \sum_{a,b} g_{ab}
\partial_{\mu} u^{a} \partial^{\mu} u^{b}$ in local coordinates on $M$.
The Euler-Lagrange equation for this variational problem is exactly of the
form \eqref{WMType}, in local coordinates on $M$, with $\Gamma^{I}_{JK}$
the Christoffel symbols of $M$ in the local chart and $N = \dim M$ (see,
e.g., Shatah-Struwe \cite{Sh-St}). %
\subsubsection{Maxwell-Klein-Gordon Equations} %
In the following discussion, the summation convention is in effect. Greek
indices are summed from $0$ to $n$, roman indices from $1$ to $n$. Recall
that indices are raised and lowered relative to the Minkowski metric
$m_{\mu\nu} = \diag(-1,1,\dots,1)$. For example, $\square = \partial^{\mu}
\partial_{\mu}$ and $\Delta = \partial^{j} \partial_{j}$. We denote by $i$
the imaginary unit.

The unknowns of the equations are a one-form $A_{\mu} dx^{\mu}$ (the gauge
potential) and a scalar $\phi$, both defined on the Minkowski space-time:
\begin{align*}
A_{\mu} : \R^{1+n} &\to \R,
\\
\phi : \R^{1+n} &\to \C.
\end{align*}
The electromagnetic field is the two-form $F_{\mu \nu} = \partial_{\mu}
A_{\nu} - \partial_{\nu} A_{\mu}$. The covariant derivative relative to
the gauge potential is $$
D_{\mu} \phi = \partial_{\mu} \phi + i A_{\mu} \phi. $$

We are looking for critical points of the Lagrangian $$
\mathcal L[A_\mu,\phi] =
\int_{\R^{1+n}} \left( - \frac{1}{4} F_{\mu \nu} F^{\mu \nu} - \half
D_{\mu} \phi \overline{D^{\mu} \phi} \right) \, dt \, dx. $$
The corresponding Euler-Lagrange equations are \begin{align}
\label{MKGa}
\tag{MKGa}
\partial^\mu F_{\mu\nu} &= - \Im \bigl(\phi \overline{D_\nu \phi} \bigr),
\\
\label{MKGb}
\tag{MKGb}
D^\mu D_\mu \phi &= 0,
\end{align}
where $\Im z$ denotes the imaginary part of $z$.

Let $\chi$ be a real-valued function on $\R^{1+n}$, and consider the
transformation $(A_\mu , \phi) \to (\widetilde A_\mu , \widetilde \phi)$
given by \begin{align*}
\widetilde A_\mu &= A_\mu - \partial_\mu \chi, \\
\widetilde \phi &= e^{i\chi} \phi.
\end{align*}
Clearly, the electromagnetic field is left unchanged by the gauge
transformation $A_\mu \to \widetilde A_\mu$, and a simple calculation
reveals that if $(A_\mu,\phi)$ verifies (MKG), then so does $(\widetilde
A_\mu, \widetilde \phi)$ (keep in mind that $D_\mu$ depends on $A_\mu$).
This gives an equivalence relation on the set of pairs $(A_\mu,\phi)$
verifying (MKG), and by a \emph{solution} of the latter, we understand an
equivalence class of such pairs.

Thus, we have gauge freedom; i.e., we are free to choose any
representative of a given solution (equivalence class), and we may
stipulate a condition that the gauge potential should satisfy. The
traditional gauge conditions are: %
\begin{itemize}
\item
\emph{Lorentz:} $\partial^\mu A_\mu = 0$, %
\item
\emph{Coulomb:} $\partial^j A_j = 0$,
\item
\emph{Temporal:} $A_0 = 0$.
\end{itemize}
\paragraph{(MKG) in Lorentz gauge.}
Coupling the Lorentz condition with (MKG) yields the system
\begin{subequations}
\begin{align}
\label{MKGLorentzA}
\square A_{\mu} &= - \Im \bigl(\phi \overline{\partial_\mu \phi} \bigr) +
\abs{\phi}^2 A_\mu, \\
\label{MKGLorentzB}
\square \phi &= -2i A^\mu \partial_\mu \phi + A^\mu A_\mu \phi, \\
\label{MKGLorentzC}
\partial^\mu A_\mu &= 0.
\end{align}
\end{subequations}
Now observe that if $(A_\mu,\phi)$ satisfies \eqref{MKGLorentzA} and
\eqref{MKGLorentzB} with initial data
\begin{subequations}
\begin{alignat}{2}
\label{MKGInitialDataA}
A_\mu \init &= a_\mu, \qquad & \partial_t A_\mu \init &= b_\mu, \\
\label{MKGInitialDataB}
\phi \init &= \phi_0, \qquad & \partial_t \phi \init &= \phi_1
\end{alignat}
\end{subequations}
satisfying the constraints
\begin{equation}\label{LorentzConstraint}
b_0 = \partial^j a_j, \qquad \Delta a_0 - \abs{\phi_0}^2 a_0 = \partial^j
b_j - \Im (\phi_0 \overline \phi_1 ),
\end{equation}
then \eqref{MKGLorentzC} is automatically satisfied. For by
\eqref{MKGLorentzA} and \eqref{MKGLorentzB}, $u = \partial^\mu A_\mu$
solves $$
\square u = \abs{\phi}^2 u,
$$
and by \eqref{MKGInitialDataA} and \eqref{LorentzConstraint}, $u \init =
\partial_t u \init = 0$. By uniqueness of solutions, $u = 0$.

Thus, \eqref{MKGLorentzC} is equivalent to the constraint
\eqref{LorentzConstraint} on the initial data, so we are left with
\eqref{MKGLorentzA} and \eqref{MKGLorentzB}. Therefore, (MKG) in Lorentz
gauge is schematically of the form $\square u = u \partial u + u^3$.
Unfortunately\footnote{See our discussion concerning
the first iterate in section \ref{Motivation} below.}, generic equations of
this type
do not have good local regularity properties, so the Lorentz gauge is not
very useful for our purposes.
\paragraph{(MKG) in Coulomb gauge.}
Coupling the Coulomb condition with (MKG) gives \begin{subequations}
\begin{align}
\label{MKGCoulombA}
\Delta A_{0} &= - \Im \bigl(\phi \overline{\partial_t \phi} \bigr) +
\abs{\phi}^2 A_0, \\
\label{MKGCoulombB}
\square A_{j} &= - \Im \bigl(\phi \overline{\partial_j \phi} \bigr) +
\abs{\phi}^2 A_j - \partial_j \partial_t A_0,
\\
\label{MKGCoulombC}
\square \phi &= - 2i A^j \partial_j \phi + 2i A_0 \partial_t \phi + i
(\partial_t A_0) \phi + A^\mu A_\mu \phi, \\
\label{MKGCoulombD}
\partial^j A_j &= 0.
\end{align}
\end{subequations}
Here we have split the gauge potential into its time component $A_0$ and
its spatial component $A = A_j dx^j$. We prescribe initial data at time $t
= 0$: \begin{subequations}\label{MKGCoulombInitialData} \begin{alignat}{2}
\label{MKGCoulombInitialDataA}
A_j \init &= a_j, \qquad & \partial_t A_j \init &= b_j, \\
\label{MKGCoulombInitialDataB}
\phi \init &= \phi_0, \qquad & \partial_t \phi \init &= \phi_1.
\end{alignat}
\end{subequations}
No initial condition is imposed on $A_0$; if we set $a_0 = A_0 \init$,
then by \eqref{MKGCoulombA}, $\Delta a_0 - \abs{\phi_0}^2 a_0 = - \Im
(\phi_0 \overline \phi_1 )$.

Equation \eqref{MKGCoulombD} is automatically satisfied if the data are
divergence-free:
\begin{equation}\label{CoulombDataConstraint}
\partial^j a_j = \partial^j b_j = 0.
\end{equation}
For if $(A_0,A,\phi)$ satisfies \eqref{MKGCoulombA}--\eqref{MKGCoulombC},
then $u = \partial^j A_j$ solves $\square u = \abs{\phi}^2 u$, and if
\eqref{MKGCoulombInitialData} and \eqref{CoulombDataConstraint} are
satisfied, then $u \init = \partial_t u \init = 0$.

We are then left with the equations
\eqref{MKGCoulombA}--\eqref{MKGCoulombC}. The first of these, being an
elliptic equation, is relatively easy to handle, so we leave it out of our
model equations. The two remaining equations have terms of three types on
the right hand side: \begin{itemize}
\item ``Elliptic terms'' involving $A_0$; these are collectively denoted
by $\mathcal E$. %
\item Cubic terms in $A_j$ and $\phi$; these are collectively denoted by
$\mathcal C$. %
\item Quadratic terms with a null-form structure. \end{itemize}
The terms falling into the latter category are $- \Im \bigl(\phi
\overline{\partial_j \phi} \bigr)$ and $- 2i A^j \partial_j \phi$. We now
uncover the null-form structure inherent in these expressions (due to the
Coulomb condition).

Split $\phi$ into its real and imaginary parts: $\phi = u + iv$. Then $$
- \Im \bigl(\phi \overline{\partial_j \phi} \bigr) = u \partial_j v - v
\partial_j u, $$
so \eqref{MKGCoulombB} reads, as an equation of (time-dependent) one-forms
on $\R^n$: $$
\square A = u dv - v du + \mathcal C - d(\partial_t A_0). $$
Apply $d$ to both sides:
$$
\square (dA) = 2 du \wedge dv + d\mathcal C. $$
But
$$
du \wedge dv = \half Q_{jk} (u,v) dx^j \wedge dx^k, $$
whence
$$
\square F_{jk} = Q_{jk}(u,v) + \partial \mathcal C. $$
The Coulomb gauge condition implies that $\partial^k F_{jk} = - \Delta
A_j$, so we have $$
- \Delta \square A_j = \partial^k Q_{jk} (u,v) + \partial^2 \mathcal C. $$
Thus, modulo Riesz operators,
\begin{equation}\label{MKGCoulombBmodel}
\square A = D^{-1} Q(\Re \phi, \Im \phi) + \mathcal C, \end{equation}
where $Q$ is some linear combination of the null forms\footnote{To be
precise, the $j$-th component of $Q$ is $\sum_{k} R_{k} Q_{jk}$, where
$R_{k} = D^{-1} \partial^{k}$ is the $k$-th Riesz operator. Since we work
with norms which only depend on the size of the Fourier transform, we
ignore the Riesz operators.} $Q_{jk}$. Since the cubic term $\mathcal C$
is easier to estimate, we leave it out of our model problem.

Now consider equation \eqref{MKGCoulombC}. Separating real and imaginary
parts, we have \begin{align*}
\square u &= 2 A \cdot \nabla v + \mathcal C + \mathcal E, \\
\square v &= - 2 A \cdot \nabla u + \mathcal C + \mathcal E. \end{align*}
(Here we consider $A$ as a vector field by raising its indices; $\nabla$
denotes the gradient in the space variables.) We claim that
the terms $A \cdot \nabla u$ and $A \cdot \nabla v$ have a null-form
structure, due to the fact that $A$ is divergence-free (by the Coulomb
condition). Let $B_{jk}$ be the unique solution of
\begin{equation}\label{DefOfB}
\Delta B_{jk} = \partial_j A_k - \partial_k A_j \end{equation}
(with appropriate regularity assumptions). By the Coulomb condition,
\begin{equation}\label{PropertyOfB}
\Delta \partial^j B_{jk} = \Delta A_k, \quad \text{which implies} \quad
\partial^j B_{jk} = A_k.
\end{equation}
Thus,
$$
A \cdot \nabla u = \partial^j B_{jk} \partial^k u = \half Q_{jk}( u,
B^{jk}). $$
The above equations for $u = \Re \phi$ and $v = \Im \phi$ can therefore be
rewritten \begin{align*}
\square \Re \phi &= Q_{jk}( \Im \phi, B^{jk}) + \mathcal C + \mathcal E, \\
\square \Im \phi &= Q_{jk}( B^{jk}, \Re \phi) + \mathcal C + \mathcal E.
\end{align*}
But in view of \eqref{DefOfB}, $B$ is of the form $D^{-1} A$ modulo Riesz
operators. Combining this with \eqref{MKGCoulombBmodel} and discarding the
terms $\mathcal C$ and $\mathcal E$ throughout, we obtain a system of the
form (``MKG''), which is our model for (MKG). %
\subsubsection{Yang-Mills Equations}
Let $G$ be one of the classical, compact Lie groups of matrices (such as
$\mathrm{SO}(k,\R)$ or $\mathrm{SU}(k,\C)$), and let $\mathfrak g$ be its
Lie algebra. The unknown is a $\mathfrak g$-valued one-form $A_\mu dx^\mu$
on $\R^{1+n}$. The corresponding covariant derivative is $$
D_\mu H = \partial_\mu H + [A_\mu,H],
$$
where $H$ is any $\mathfrak g$-valued tensor field on $\R^{1+n}$ and
$[\cdot,\cdot]$ is the matrix commutator.

The curvature is the $\mathfrak g$-valued two-form $$
F_{\mu \nu} = \partial_\mu A_\nu - \partial_\nu A_\mu + [A_\mu,A_\nu]. $$
The Lagrangian is
$$
\mathcal L[A_\mu] = - \frac{1}{4} \int \innerprod{F_{\mu \nu}}{F^{\mu
\nu}} \, dt \, dx, $$
where $\innerprod{\cdot}{\cdot}$ is the inner product on $\mathfrak g$
inherited from the ambient space (e.g., $\mathnormal{SO}(k,\R)$ embeds in
$\R^{k^2}$, so its Lie algebra can be viewed as a subspace of the latter).
The Euler-Lagrange equations are
\begin{equation}\label{YM}\tag{YM}
D^\nu F_{\mu \nu} = 0.
\end{equation}

Let $O$ be a $G$-valued function on $\R^{1+n}$. Consider the gauge
transformation $A_\mu \to \widetilde A_\mu$, given by
$$
\widetilde A_\mu = O A_\mu O^{-1} - \partial_\mu O O^{-1}. $$
A calculation shows that the curvature then transforms into $$
\widetilde F_{\mu \nu} = O F_{\mu \nu} O^{-1}. $$
Denoting by $\widetilde D_\mu$ the covariant derivative corresponding to
$\widetilde A_\mu$, we then have
$$
\widetilde D^\nu \widetilde F_{\mu \nu} = O D^\nu F_{\mu \nu} O^{-1}, $$
so \eqref{YM} is invariant under gauge transformations. We therefore have
gauge freedom, and may impose a gauge condition on $A_\mu$. %
\paragraph{(YM) in Coulomb gauge.}
Relative to the Coulomb condition $\partial^j A_j = 0$, \eqref{YM} takes
the form (see \cite{Kl-Ma0.3})
\begin{subequations}\label{YMCoulomb}
\begin{align}
\label{YMCoulombA}
\Delta A_0 &= 2 [\partial^j A_0, A_j] + [A^j, \partial_t A_j] + [ A^j,
[A_0,A_j]],
\\
\label{YMCoulombB}
\square A_j + \partial_t \partial_j A_0 &= - 2 [A^k, \partial_k A_j] +
[A^k, \partial_j A_k] + [\partial_t A_0, A_j] + 2[A_0, \partial_tA_j] \\
\notag
& \qquad - [A_0, \partial_j A_0] - [A^k, [A_k,A_j]] + [A_0,[A_0,A_j]], \\
\label{YMCoulombC}
\partial^j A_j &= 0.
\end{align}
\end{subequations}
Unfortunately, assuming the existence of a global Coulomb gauge forces a
restrictive smallness assumption on the initial data. In \cite{Kl-Ma0.3}
this difficulty was resolved by using local arguments. Following
\cite{Kl-Ta}, we ignore this complication, and derive our model equation
from the system \eqref{YMCoulomb}.

As in the discussion of (MKG), \eqref{YMCoulombC} reduces to a constraint
on the initial data. The equation for $A_0$ is elliptic, so we ignore it.
As for \eqref{YMCoulombB}, we only retain the first two terms on the
right, since the other terms either involve $A_0$ (for which we expect to
have better estimates than for $A_j$), or are cubic.

Now write \eqref{YMCoulombB} as an equation of time-dependent, $\mathfrak
g$-valued one-forms on $\R^{n}$ (ignoring all but the first two terms on
the right):
$$
\square A + d(\partial_{t} A_{0}) = S + T, $$
where $A = A_{j} dx^{j}$, $S = - 2 [A^k, \partial_k A_j] dx^{j}$ and $T =
[A^k, \partial_j A_k] dx^{j}$.
Apply the exterior derivative $d$ to both sides: $$
\square dA = dS + dT.
$$
Let $B$ be the two-form (in this case $\mathfrak g$-valued) determined by
equation \eqref{DefOfB}. Thus
$$
\square \Delta B_{jk} = \partial_{j} S_{k} - \partial_{k} S_{j} +
\partial_{j} T_{k} - \partial_{k} T_{j}. $$
By \eqref{PropertyOfB}, it follows that
$$
- \square \Delta A_{j} = \partial^{k} \left( \partial_{j} S_{k} -
\partial_{k} S_{j} + \partial_{j} T_{k} - \partial_{k} T_{j} \right), $$
so for the purposes of estimates in frequency space, we may replace
\eqref{YMCoulombB} by
\begin{equation}\label{YMCoulombModel}
\square A = S + D^{-1} dT.
\end{equation}

It remains to identify the null form structure hidden in $S$ and $dT$. To
begin with, we have
\begin{align*}
S_{j} &= [ \partial_{k} B^{kl}, \partial_{l} A_{j}] = \half [ \partial_{k}
B^{kl} - \partial_{k} B^{lk}, \partial_{l} A_{j}] \\
&= \half \left( [ \partial_{k} B^{kl} , \partial_{l} A_{j}] - [
\partial_{l} B^{kl} , \partial_{k} A_{j}] \right) \\
&= \half \left( \partial_{k} B^{kl} \partial_{l} A_{j} - \partial_{l}
B^{kl} \partial_{k} A_{j} + \partial_{k} A_{j} \partial_{l} B^{kl} -
\partial_{l} A_{j} \partial_{k} B^{kl} \right),
\end{align*}
so each entry of the matrix $S_{j}$ is a linear combination of terms of
the form $Q_{kl}(A,B)$, where $A$ and $B$ stand for any two entries of
$A_{j}$ and $B_{kl}$. But by \eqref{PropertyOfB}, we may replace $B$ by
$D^{-1} A$. Schematically,
\begin{equation}\label{YMCoulombFirstNullForm}
S \sim Q(A,D^{-1}A).
\end{equation}

Now consider the one-form $T$. We calculate: \begin{align*}
(dT)_{jk} &= \partial_{j} [A^l, \partial_k A_l] - \partial_{k} [A^l,
\partial_j A_l] \\
&= [\partial_{j} A^l, \partial_k A_l] - [\partial_{k} A^l, \partial_j A_l].
\end{align*}
Thus, each entry of the matrix $(dT)_{jk}$ is a linear combination of
terms of the form $Q_{jk}(A,A')$, where $A$ and $A'$ stand for any two
entries of $A_{l}$, $1 \le l \le n$.
Combining this with \eqref{YMCoulombFirstNullForm} and
\eqref{YMCoulombModel}, we arrive at the model \eqref{YMType} for the
Yang-Mills
equations.
\subsubsection{Wave Maps Model Problem}
The (WMM) equation arises from a simple reformulation of Wave Maps whose
target manifold has a bi-invariant Lie group structure.
Let $G$ be a Lie group, and let $\mathfrak g$ be its Lie algebra,
identified with the tangent space $T_e G$, where $e$ is the unit in $G$.
For any $a \in G$, we denote by $L_a$ and $R_a$ the left and right
translation operators on $G$, given by $L_a(g) = ag$ and $R_a(g) = ga$.
Their derivatives are denoted by $L_{a*}$ and $R_{a*}$ respectively.

Assume that $G$ is endowed with a Riemannian metric $h$ which is
bi-invariant; i.e., $h(L_{a*} X, L_{a*} Y) = h(X,Y)$ and $h(R_{a*} X,
R_{a*} Y) = h(X,Y)$ for all $a \in G$ and all tangent vector fields $X$
and $Y$.

Let $u : \R^{1+n} \to G$. Then for all $0 \le \mu \le n$ and $(t,x) \in
\R^{1+n}$, $\partial_\mu u(t,x)$ is a vector
in the tangent space $T_{u(t,x)}G$, and we move this vector into the Lie
algebra $T_{e}G$ by left translation. More precisely, we define a
$\mathfrak g$-valued one-form $A_\mu dx^\mu$ by
$$
A_\mu = L_{u^{-1}*}(\partial_\mu u),
$$
where $u^{-1}$ denotes the group inverse.

It turns out that $u$ is a wave map if and only if $A_\mu dx^\mu$
satisfies \begin{equation}\label{WaveMapsModelA}
\begin{split}
\partial^\mu A_\mu &= 0,
\\
\partial_\mu A_\nu - \partial_\nu A_\mu &= - [A_\mu,A_\nu], \end{split}
\end{equation}
where $[\cdot,\cdot]$ is the Lie bracket. The advantage of this
formulation of the wave maps problem is that it avoids the use of local
charts in the target manifold. See Christodoulou and Tahvildar-Zadeh
\cite{Ch-Ta} for an application of this system to prove global regularity
of spherically symmetric wave maps for $n = 2$.

First, let us see how the model equation \eqref{CFWMType} arises from this
system. We start by transforming the variables, using the nonlocal
operators $R_{\mu} = D^{-1} \partial_{\mu}$. We assume that $\mathfrak g$
is a Lie algebra of matrices, and that $[\cdot,\cdot]$ is the usual matrix
commutator. Set
$$
\bar A_i = A_i + R_0 R_i A_0.
$$
Then it follows from \eqref{WaveMapsModelA} that the one-form $A_0 dx^0 +
\bar A_i dx^i$ satisfies
\begin{align*}
\square A_0 &= \partial^i [ A_0 , \bar A_i - R_0 R_i A_0 ] \\ \partial^i
\bar A_i &= 0 \\
\partial_i \bar A_j - \partial_j \bar A_i &= [ \bar A_j - R_0 R_j A_0,
\bar A_i - R_0 R_i A_0 ].
\end{align*}
Since the spatial part $\bar A_i$ satisfies an elliptic Hodge system, it
is easier to estimate than the temporal part $A_0$, and therefore we
ignore it. In other words, we set $\bar A_i = 0$ in the equation for
$A_0$, which gives the model problem \eqref{CFWMType}.

We remark that if we set $A_{0} = 0$, then the above system describes a
time-independent wave map (a \emph{harmonic map}) $u : \R^{n} \to G$. This
formulation of the harmonic map problem was used by F. H\'elein \cite{He}
to prove regularity of weakly harmonic maps in dimension $n = 2$.

We now outline the derivation of the system \eqref{WaveMapsModelA}.
Following \cite[Section 3.1]{Ch-Ta}, we first choose an orthonormal basis
$\Omega_I$ of $\mathfrak g$, and we let $\omega^I$ be the dual basis of
left-invariant one-forms on $G$. Let $e^I_{JK}$ be the structure
constants, defined by
$$
[\Omega_J, \Omega_K] = e^I_{JK} \Omega_I. $$
Express $A_\mu$ relative to the basis:
$$
A_\mu = \psi^I_\mu \Omega_I.
$$
Since $\partial_\mu u = L_{u*} A_\mu$, it follows that
\begin{equation}\label{EquationForPsi}
\psi^I_\mu(t,x) = \omega^I_{u(t,x)} \bigl(\partial_\mu u(t,x) \bigr),
\end{equation}
which gives the precise dependence of $\psi^I$ on $u$ and $\partial u$.

Recall that the wave map Lagrangian is $\mathcal L[u] = \int L \, dt \,
dx$, where
$$
L(u,\partial u) = \half \innerprod{du}{du} = \half h( \partial_\mu u,
\partial^\mu u) = \half h(A_\mu, A^\mu)
= \half \sum_{\mu,I} (\psi^I_\mu)^2.
$$
Here we used the left invariance of $h$. Using the last expression for
$L(u,\partial u)$, together with \eqref{EquationForPsi} and the Cartan
structure equations $$
d \omega^I = -\half e^I_{JK} \omega^J \wedge \omega^K, $$
a calculation reveals (see \cite{Ch-Ta} for the details) that the
Euler-Lagrange equation takes the form
\begin{align}
\label{WaveMapsModelB}
\partial^\mu \psi^I_\mu &= \sum_{J,K} e^J_{KI} \psi^{J\mu} \psi^K_\mu. \\
\intertext{A direct calculation also gives}
\notag
\partial_\mu \psi^I_\nu - \partial_\nu \psi^I_\mu &= - \sum_{J,K} e^I_{JK}
\psi^J_\mu \psi^K_\nu. \end{align}

Observe that the last equation is equivalent to the second equation in
\eqref{WaveMapsModelA}. We claim that \eqref{WaveMapsModelB} is equivalent
to
$$
\partial^\mu A_\mu = [A^\mu,A_\mu].
$$
Since the right hand side vanishes, we obtain the first equation in
\eqref{WaveMapsModelA}.

To prove the claim, we only have to note that, because of the assumption
that the metric on $G$ is bi-invariant, the structure constants satisfy
$$
e^I_{JK} = e^J_{KI}.
$$
Equivalently,
$$
h\bigl( \Omega_I, [\Omega_J, \Omega_K] \bigr) + h\bigl( \Omega_J,
[\Omega_I, \Omega_K] \bigr) = 0. $$
To see this when $G$ is a group of matrices, let $e^{X}$ denote the
exponential map, where $X \in \mathfrak g$. Fix $X,Y,Z \in \mathfrak g$.
By the bi-invariance of $h$,
$$
h \bigl( e^{tX} Y e^{-tX}, e^{tX} Z e^{-tX} \bigr) = h (Y,Z).
$$
Since
$$
\frac{d}{dt} \bigl( e^{tX} Y e^{-tX} \bigr) \init = XY - YX = [X,Y], $$
it follows that
$$
\frac{d}{dt} h \bigl( e^{tX} Y e^{-tX}, e^{tX} Z e^{-tX} \bigr) \init = h(
[X,Y] , Z) + h( Y, [X,Z] ) = 0,
$$
which proves the claim.
\subsection{Motivation of the Main Theorem}\label{Motivation}
Consider the system
\begin{equation}\label{ClassicalSystem}
\square u = F(u,\partial u),
\end{equation}
where $u : \R^{1+n} \to \R^N$,
$\partial u = (\partial_t u, \partial_1 u, \dots, \partial_n u)$ and $F$
is a smooth $\R^N$-valued function satisfying $F(0) = 0$. For this
equation one has the following standard existence and uniqueness result
(concerning the proof, see Example \ref{ClassicalLocalExistenceExample}). %
\begin{ClassicalLocalEx}
Equation \eqref{ClassicalSystem} is locally well-posed for initial data in
$H^s \times H^{s-1}(\R^n)$ for all $s > \frac{n}{2} + 1$.
\end{ClassicalLocalEx}
This result is far from being sharp insofar as the regularity assumption
on the initial data is concerned.

To understand better the issue of optimal local well-posedness, in the
context of our examples (wave maps, Maxwell-Klein-Gordon and Yang-Mills
equations), we need to define the critical well-posedness (henceforth
abbreviated WP) exponent $s_{c}$. All our equations have a natural scaling
associated to them, and $s_{c}$ is the unique value of $s$ for which the
$\dot H^{s} \times \dot H^{s-1}$-norm of the initial data is invariant
under this scaling. For example, if $u$ solves \eqref{WMType}, then so
does $$
u_{\lambda}(t,x) = u(\lambda t, \lambda x), $$
for any $\lambda > 0$. Since $\Sobdotnorm{u_{\lambda}(t)}{s} =
\lambda^{\frac{n}{2}-s} \Sobdotnorm{u(\lambda t)}{s}$, the critical WP
exponent for \eqref{WMType} is $s_{c} = \frac{n}{2}$.

The same principle works for (MKG), (YM) and \eqref{CFWMType}. In fact,
they all have critical WP exponent $s_{c} = \frac{n-2}{2}$.

With this definition we formulate the following, taken from \cite{Kl}: %
\begin{GeneralWPConjecture}
\begin{enumerate}
\item
For all basic field theories the initial value problem is locally well
posed for initial data in $H^{s} \times H^{s-1}$, $s > s_{c}$. %
\item
The basic field theories are weakly\footnote{The solutions may fail to
depend smoothly (analytically) on the data.} globally well-posed for all
initial
data with small $H^{s_{c}} \times H^{s_{c}-1}$-norm. %
\item
The basic field theories are ill posed for initial data in $H^{s} \times
H^{s-1}$, $s < s_{c}$.
\end{enumerate}
\end{GeneralWPConjecture}
Our Main Theorem establishes part (i) of this conjecture for the equations
in section \ref{BasicEquations}. We prove local existence by Picard
iteration in a suitable Banach space, as discussed in section
\ref{TheIterationSpace}. The $0$-th iterate $u_{0}$ corresponding to a
Cauchy problem
$$
\square u = \mathcal N(u), \qquad (u,\partial_t u) \init = (f,g) $$
is just the homogeneous part of the solution: $$
\square u_{0} = 0, \qquad (u,\partial_t u) \init = (f,g). $$
The subsequent iterates are given inductively by $$
u_{j+1} = u_0 + \square^{-1} \mathcal N(u_j) $$
for $j \ge 0$, where $\square^{-1}$ is the operator which to any
sufficiently regular $F$ assigns the solution $v$ of $\square v = F$ with
$(v,\partial_t v) \init = 0$.

If we are to prove existence of a local solution of $\square u = \mathcal
N(u)$ with initial data in $H^{s} \times H^{s-1}$ by iteration, we must be
able to prove that the iterates remain in the data space:
\begin{equation}\label{IteratesCondition}
f \in H^{s}, \quad g \in H^{s-1} \implies u_{j}(t) \in H^{s},\quad
\partial_{t} u_{j}(t) \in H^{s-1}
\end{equation}
for all $j \ge 0$ and all $t$ in some interval $(0,T)$. For $j = 0$,
\eqref{IteratesCondition} is trivial, but the case $j = 1$ already offers
valuable insights. We will say that the first iterate is \emph{WP for
initial data in $H^{s}$} if \eqref{IteratesCondition} holds for $j = 1$
and all $(f,g) \in H^{s} \times H^{s-1}$.
\begin{example}\label{FirstIterateExample1} Consider the model problem $$
\square u = (\partial_{t} u)^{2},
$$
where $u$ is real-valued. This equation has the same scaling properties as
\eqref{WMType}, hence the WP-exponent is $s_{c} = \frac{n}{2}$. We want to
find the lower bound for the set of $s$ such that the first iterate
$u_{1}$ is WP for initial data in $H^{s}$. A simple calculation involving
Duhamel's principle, done in the Appendix, shows that this reduces to
proving an estimate of the type \begin{equation}\label{IntegralEstimate}
\int_{\R^{n} \times \R^{n}} K(\xi,\eta) f(\xi) g(\eta) h(\xi + \eta) \,
d\xi \, d\eta
\lesssim \twonorm{f}{} \twonorm{g}{} \twonorm{h}{} \end{equation}
for all $f,g, h \in L^{2}(\R^{n})$, where \begin{gather}
\label{KernelA}
K(\xi,\eta) = \frac{\elwt{\xi+\eta}^{s-1}}{\elwt{\xi}^{s-1}
\elwt{\eta}^{s-1} \bigl(1 + \Delta_{\pm}(\xi,\eta) \bigr)}, \\
\label{DeltaPlusMinus}
\Delta_{+} = \abs{\xi} + \abs{\eta} - \abs{\xi + \eta}, \qquad \Delta_{-}
= \abs{\xi + \eta} - \hypwt{\xi}{\eta}. \end{gather}
Here we use the notation $\elwt{\cdot} = 1 + \abs{\cdot}$.

In the Appendix we prove the following result concerning integral
estimates of the type \begin{equation}\label{SpecificIntegralEstimate}
\int_{\R^{n} \times \R^{n}} \frac{f(\xi) g(\eta) h(\xi + \eta)}
{\elwt{\xi}^{a} \elwt{\eta}^{b} \bigl(1 + \Delta_{\pm}(\xi,\eta)
\bigr)^{c}}
\, d\xi \, d\eta
\lesssim \twonorm{f}{} \twonorm{g}{} \twonorm{h}{}, \end{equation}
where $\Delta_{\pm}$ are given by \eqref{DeltaPlusMinus}. %
\begin{proposition}\label{IntegralEstimateProposition} Let $a,b,c \ge 0$.
Then \eqref{SpecificIntegralEstimate} holds if $a + b + c > \frac{n}{2}$
and $c < \frac{n-1}{4}$. \end{proposition}
It should be remarked that the estimate fails if $a + b + c < \frac{n}{2}$
or $c \ge \frac{n-1}{4}$, although we do not prove this here.

The kernel \eqref{KernelA} clearly satisfies (assuming $s \ge 1$) $$
K \lesssim \left( \elwt{\xi}^{1-s} + \elwt{\eta}^{1-s} \right) \bigl(1 +
\Delta_{\pm}(\xi,\eta) \bigr)^{-\alpha} $$
for any $0 \le \alpha \le 1$. In view of Proposition
\ref{IntegralEstimateProposition}, we must have $\alpha < \frac{n-1}{4}$.
In order for the other hypothesis of Proposition
\ref{IntegralEstimateProposition} to be satisfied, we need
$$
s-1 + \min \left(1,\frac{n-1}{4} \right) > \frac{n}{2}, $$
i.e., $s > \max \left( \frac{n}{2}, \frac{n+5}{4} \right)$. \end{example}
Thus, for the model equation $\square u = (\partial_{t} u)^{2}$ in
dimension $n = 3$, the above example shows that the first iterate is WP
for initial data in $H^{s}$ if $s > 2 = s_{c} + \half$; in fact, one can
show that this fails to be true if $s \le 2$. This should be compared to
the counterexamples of Lindblad \cite{Li1} in dimension $n = 3$, which
show that there are equations of the type $\square u = q(\partial u)$,
where $q$ is a quadratic form on $\R^{4}$, which are ill posed for data in
$H^{2} \times H^{1}(\R^{3})$. However, if the quadratic form $q$ is of
null form type, one can go almost all the way to the critical WP-exponent
$s_{c}$. The next two examples verify this at the level of the first
iterate.
\begin{example}\label{FirstIterateExample2} Consider the
equation\footnote{The equation below can in fact be trivially solved and
analyzed, see the first page in the introduction of \cite{Kl-Ma0.1}.}
$$
\square u = Q_{0}(u,u),
$$
where $u$ is real-valued. Again the question of WP of the first iterate
leads to the problem of proving an estimate of the type
\eqref{IntegralEstimate}, but because of the special null structure of the
operator $Q_{0}$, the singular factors $\Delta_{\pm}$ cancel out
completely from the denominator of the kernel. In fact, $K$ is given by $$
K(\xi,\eta) = \elwt{\xi}^{-s} + \elwt{\eta}^{-s}, $$
so by Proposition \ref{IntegralEstimateProposition}, the first iterate is
WP for data in $H^{s}$, $s > s_{c} = \frac{n}{2}$.
\end{example}
\begin{example}\label{FirstIterateExample3} Consider the equation
$$
\square u = Q(u,u),
$$
where $u$ is vector-valued and $Q(u,u)$ is a vector whose $I$-th component
is a linear combination of $Q_{ij}(u^{J},u^{K})$ for all $i,j,J$ and $K$.
As in the preceding example, there is a cancellation due to the null
strucure of $Q_{ij}$, but in this case we only get rid of half a power of
$\Delta_{\pm}$. In fact, $K$ is now given by
$$
K(\xi,\eta) = \left( \elwt{\xi}^{-s+\half} + \elwt{\eta}^{-s+\half} \right)
\bigl(1 + \Delta_{\pm}(\xi,\eta) \bigr)^{-\half} $$
so the first iterate is WP for data in $H^{s}$, $s > \max \left(
\frac{n}{2}, \frac{n+3}{4} \right)$.

By an obvious modification, if we consider instead the equation $$
\square u = Q(D^{-1}u,u),
$$
we find that the first iterate is WP for data in $H^{s}$, $s > \max \left(
\frac{n-2}{2}, \frac{n-1}{4} \right)$. \end{example}
The preceding examples are worked out in more detail in the Appendix. %
\subsection{Notation}\label{Notation}
Throughout the paper, $p \lesssim q$ means that $p \le C q$ for some
positive constant $C$. Similarly, $\simeq$ means $=$ modulo a positive
constant. The notation $p \sim q$ means $p \lesssim q \lesssim p$.

If $\X$ is a separable Banach space, $L^{p}(\R^{k},\X)$ denotes the usual
$L^{p}$ space, relative to Lebesgue measure on $\R^{k}$, and we write
$L^{p}(\R^{k}) = L^{p}(\R^{k},\C)$. If $(\alpha,\beta) \in \R^{k} \times
\R^{l}$, we define the mixed norm
$\norm{f(\alpha,\beta)}_{L^{q}_{\alpha}(L^{r}_{\beta})}$ by first taking
the $L^{r}(\R^{l})$-norm in $\beta$, followed by the $L^{q}(\R^{k})$-norm
in $\alpha$. Thus, $L^{q}_{\alpha}(L^{r}_{\beta}) =
L^{q}\bigl(\R^{k},L^{r}(\R^{l})\bigr)$ if $r < \infty$.

The space of Schwartz functions on $\R^{k}$ is denoted by
$\Schwartz(\R^{k})$, and its dual, the space of tempered distributions, is
written $\Schwartz'(\R^{k})$. If $u \in \Schwartz'(\R^{1+n})$ and it makes
sense to restrict $u$ to any time-slice $\{t\} \times \R^{n}$, we write
$u(t)$ instead of $u(t,\cdot)$. The Fourier transform of a tempered
distribution $u$ is denoted by $\Fourier u$ or $\widehat u$, in any
space-dimension. In frequency space we use coordinates $(\tau,\xi) = \Xi =
(\Xi^{0},\dots,\Xi^{1})$, where $\tau \in \R$ and $\xi \in \R^{n}$
correspond to the time variable $t$ and the space variable $x$
respectively. The Lorentzian inner product on $\R^{1+n}$ is denoted by
$\biginnerprod{\Xi}{\widetilde \Xi}$. Thus
$$
\biginnerprod{\Xi}{\widetilde \Xi} = \sum_{\mu = 0}^n \Xi_{\mu} \widetilde
\Xi^{\mu} = - \Xi^0 \widetilde \Xi^0
+ \sum_{j = 1}^n \Xi^j \widetilde \Xi^j, $$
and the symbol of the wave operator $\square$ is $- \innerprod{\Xi}{\Xi} =
\tau^2 - \abs{\xi}^2$. By $\abs{\Xi}$ we always mean the Euclidean norm.

Let $\Lambda^{\alpha}$, $\Lambda_{+}^{\alpha}$ and $\Lambda_{-}^{\alpha}$
be the multipliers given by \begin{align*}
\widehat{\Lambda^{\alpha} f}(\xi) &= \bigl( 1 + \abs{\xi}^{2}
\bigr)^{\alpha/2} \widehat{f}(\xi),
\\
\widehat{\Lambda_{+}^{\alpha} u}(\Xi) &= \bigl( 1 + \abs{\Xi}^{2}
\bigr)^{\alpha/2} \widehat{u}(\Xi),
\\
\widehat{\Lambda_{-}^{\alpha} u}(\Xi) &= \left( 1 + \frac{
\innerprod{\Xi}{\Xi}^{2} } {1 + \abs{\Xi}^{2} } \right)^{\alpha/2}
\widehat{u}(\Xi). \end{align*}
Observe that these operators are isomorphisms of $\Schwartz(\R^{1+n})$ as
well as $\Schwartz'(\R^{1+n})$. Moreover, $\Lambda^{\alpha}$ may also be
regarded as an isomorphism of $\Schwartz(\R^{n})$ and
$\Schwartz'(\R^{n})$, since it only acts in the space variable.

We also need homogeneous versions of these operators: Let $D^{\alpha}$,
$D_{+}^{\alpha}$ and $D_{-}^{\alpha}$ be the multipliers with symbols $$
\abs{\xi}^\alpha,
\quad (\abs{\tau} + \abs{\xi})^\alpha,
\quad \hypwt{\tau}{\xi}^\alpha
$$
respectively.

If $u, v \in \Schwartz'$ and $\widehat u, \widehat v$ are tempered
functions, we write $u \preceq v$ iff $\abs{\widehat u} \le \widehat v$,
and $\precsim$ means $\preceq$ up to a constant. If $u = (u^1,\dots,u^N)$
and $v = (v^1,\dots,v^N)$, then $u \preceq v$ (resp. $u \precsim v$) means
$u^I \preceq v^I$ (resp. $u^I \precsim v^I$) for $I = 1, \dots, N$.

If $\X$ is a normed vector space of tempered distributions such that
$\widehat u$ is a tempered function whenever $u \in \X$, then we say that
the norm on $\X$ depends only on the size of the Fourier transform if
$\norm{u} = \norm{v}$ whenever $\abs{\widehat u} = \abs{\widehat{v}}$, and
we say that the norm is compatible with the relation $\preceq$ if
$\norm{u} \le \norm{v}$ whenever $u \preceq v$.

The solution of the homogeneous wave equation $\square u = 0$ with initial
data $(u,\partial_{t} u) \init = (f,g)$ can be decomposed into \emph{half
waves}: $u = u_{+} + u_{-}$, where $u_{\pm}(t) = e^{\pm itD} \half ( f \pm
i^{-1} D^{-1} g )$.
We shall often restrict ourselves to the reduced initial value problem
$\square u = 0$ with data $(u,\partial_t u) \init = (f,0)$; the general
case can easily be reduced to this.

The symbol $\hookrightarrow$ means continuous inclusion. For example, we
have the Sobolev embeddings
\begin{align}
\label{SobolevEmbedding1}
\dot H^{\frac{n}{2}-\frac{n}{p}}(\R^{n}) &\hookrightarrow L^{p}(\R^{n})
\quad \text{iff} \quad 2 \le p < \infty, \\
\label{InftySobolevEmbedding}
H^s &\hookrightarrow L^\infty(\R^n)
\quad \text{iff} \quad s > \frac{n}{2}.
\end{align}
If $\X,\Y,\mathcal Z$ are normed function spaces, then $\X \cdot \Y
\hookrightarrow \mathcal Z$ means that $\norm{uv}_\mathcal Z \lesssim
\norm{u}_\X \norm{v}_\Y$ for all $(u,v) \in \X \times \Y$. More generally,
if $B(u,v)$ is some bilinear operator on $\X \times \Y$, we shall write
$B(\X,\Y) \hookrightarrow \mathcal Z$ to mean that $B$ is bounded from $\X
\times \Y$ into $\mathcal Z$.
\section{Estimates for the Wave Equation}\label{WaveEstimates}
Here we review some of the well known estimates for solutions of the
homogeneous wave equation $\square u = 0$ which will be needed throughout
the paper.

Without loss of generality, we restrict ourselves to the reduced initial
value problem
\begin{equation}\label{HomogeneousEquation}
\square u = 0, \qquad (u,\partial_{t} u) \init = (f,0). \end{equation}
Estimates for the general case $(u,\partial_t u) \init = (f,g)$ can easily
be deduced from this.

We start by recalling the Strichartz type estimates
\begin{equation}\label{HomogeneousLinearStrichartz}
\mixednorm{u}{q}{r}
\lesssim \Sobdotnorm{f}{s},
\end{equation}
where $u$ solves \eqref{HomogeneousEquation} and $\Sobdotnorm{f}{s} =
\twonorm{D^{s} f}{}$. Scaling considerations impose the condition
\begin{equation}\label{ScalingCondition}
s = \frac{n}{2} - \frac{n}{r} - \frac{1}{q}. \end{equation}
The pair $(q,r)$ is said to be \emph{wave admissible} if
\begin{equation}\label{WaveAdmissible}
2 \le q \le \infty, \quad
2 \le r < \infty, \quad
\frac{2}{q} \le (n-1) \left(\half - \frac{1}{r} \right). \end{equation}
For the proof of the next result, and further references, see
\cite{Ke-Ta}. The case $q,r = 4$, $n = 3$ corresponds to the original
inequality of Strichartz \cite{Str}. %
\begin{ThmA}
The estimate \eqref{HomogeneousLinearStrichartz} is satisfied by the
solution of \eqref{HomogeneousEquation} for all $f \in \dot H^{s}$ iff
\eqref{ScalingCondition} holds and $(q,r)$ is wave admissible. \end{ThmA}
The next result is a generalization of Theorem A to bilinear estimates of
the type \begin{equation}\label{BilinearEstimate}
\mixednorm{D^{-\sigma}(uv)}{q/2}{r/2}
\lesssim \Sobdotnorm{f}{s_{1}} \Sobdotnorm{g}{s_{2}}, \end{equation}
where $u$ and $v$ solve
\begin{equation}\label{TwoHomogeneousEquations}
\square u = \square v = 0, \qquad (u,\partial_{t}u) \init = (f,0), \qquad
(v,\partial_{t}v) \init = (g,0).
\end{equation}
Note that if $\sigma = 0$ and $s_{1} = s_{2} = s$, then
\eqref{BilinearEstimate} reduces to \eqref{HomogeneousLinearStrichartz} by
H\"older's inequality.

The following theorem was first proved by Klainerman-Machedon
\cite{Kl-Ma2} in the case $q,r = 4$; the general statement was proved by
Klainerman-Tataru \cite{Kl-Ta}.
\begin{ThmB}
Let $n \ge 2$, and let $(q,r)$ be a wave admissible pair:
$$
2 \le q \le \infty, \quad
2 \le r < \infty, \quad
\frac{2}{q} \le (n-1) \left(\half - \frac{1}{r} \right)
$$
Assume that
\begin{gather*}
0 < \sigma < n - \frac{2n}{r} - \frac{4}{q}, \\
s_1, s_2 < \frac{n}{2} - \frac{n}{r} - \frac{1}{q}, \\
s_1 + s_2 + \sigma = n - \frac{2n}{r} - \frac{2}{q}. \end{gather*}
Then \eqref{BilinearEstimate} holds for all solutions of
\eqref{TwoHomogeneousEquations}.
\end{ThmB}
In practically all our applications of the above theorem, $s_{1} = s_{2}$.
It should be remarked that in the asymmetric case $s_{1} \neq s_{2}$, the
above conditions on $s_{1},s_{2}$ are not optimal (cf. the proof of
\eqref{AsymmetricYMEstimate} in section \ref{YMInformalProof}).

Now consider more general bilinear estimates, of the form
\begin{equation}\label{GeneralBilinearEstimate}
\mixednorm{D^{\gamma} D_{+}^{\gamma_{+}} D_{-}^{\gamma_{-}}(uv)}{q/2}{r/2}
\lesssim \Sobdotnorm{f}{s_{1}} \Sobdotnorm{g}{s_{2}}, \end{equation}
where $u$ and $v$ solve \eqref{TwoHomogeneousEquations}. In the case $q,r
= 4$ all such estimates are known. Special cases of the following theorem
have appeared first in \cite{Kl-Ma0.1} and later in \cite{Kl-Ma2, Kl-Se,
Kl-Ma5}. The complete solution was carried out recently by
Foschi-Klainerman
\cite{Fo-Kl}, see also \cite{Tao2}.
\begin{ThmC}
Let $n \ge 2$ and $\gamma, \gamma_{-}, \gamma_{+}, s_{1}, s_{2} \in \R$.
The estimate
$$
\twonorm{D^{\gamma} D_{+}^{\gamma_{+}}
D_{-}^{\gamma_{-}}(uv)}{(\R^{1+n})}
\lesssim \Sobdotnorm{f}{s_{1}} \Sobdotnorm{g}{s_{2}} $$
is satisfied by the solutions of
\eqref{TwoHomogeneousEquations} for all $f,g$ iff the following conditions
hold:
\begin{align*}
\gamma + \gamma_{+} + \gamma_{-} &= s_{1} + s_{2} - \frac{n-1}{2}, \\
\gamma_{-} &\ge - \frac{n-3}{4},
\\
\gamma &> - \frac{n-1}{2},
\\
s_{i} &\le \gamma_{-} + \frac{n-1}{2}, \quad i = 1,2, \\
s_{1} + s_{2} &\ge \half,
\\
(s_{i},\gamma_{-}) &\neq \left( \frac{n+1}{4}, -\frac{n-3}{4} \right),
\quad i = 1,2,
\\
(s_{1} + s_{2},\gamma_{-}) &\neq \left( \half, -\frac{n-3}{4} \right).
\end{align*}
\end{ThmC}
\section{Wave-Sobolev Spaces}\label{WaveSobolevSpaces}
We define the space $H^{s,\theta}$, which is adapted to the wave operator
on $\R^{1+n}$ in the same way that $H^s$ is adapted to the Laplacian on
$\R^n$, and we show that the estimates in Theorems A, B and C for
solutions of the homogeneous wave equation imply corresponding estimates
for elements of $H^{s,\theta}$. %
\begin{definition}
For $s,\theta \in \R$, define
$$
H^{s,\theta} = \left\{ u \in \Schwartz' : \Lambda^{s} \Lambda_{-}^{\theta}
u \in L^{2} \right\} $$
with norm $\spacetimenorm{u}{s}{\theta} = \twonorm{\Lambda^{s}
\Lambda_{-}^{\theta} u}{}$ (see section \ref{Notation} for the definition
of the operators $\Lambda^s$ and $\Lambda_-^\theta$).
\end{definition}
Since $\Lambda^{s} \Lambda_{-}^\theta(\Schwartz) = \Schwartz$ and
$\Schwartz$ is dense in $L^{2}$, it is immediate from the definition that
$\Schwartz$ is dense in $H^{s,\theta}$.

There is a remarkably simple connection between $H^{s,\theta}$ and the
space of solutions of the homogeneous wave equation with data in $H^{s}$.
In effect, every $u \in H^{s,\theta}$ is of the form
\begin{equation}\label{SimplifiedIntegralFormula}
u(t) = \frac{1}{2\pi} \int_{-\infty}^{\infty} \frac{e^{it\lambda}
u_{\lambda}(t)}{(1+\abs{\lambda})^{\theta}} d\lambda \qquad
\text{($H^{s}$-valued integral)}
\end{equation}
where $\{u_{\lambda}\}_{\lambda \in \R}$ is a one-parameter family of
solutions of \eqref{HomogeneousEquation} with data in $H^{s}$; i.e.,
$\square u_{\lambda} = 0$ and $(u_{\lambda},\partial_{t} u_{\lambda})
\init = (f_{\lambda},0)$, where $\lambda \mapsto f_{\lambda}$ belongs to
$L^{2}(\R,H^{s})$. Moreover, $\spacetimenorm{u}{s}{\theta}^{2} = \int
\Sobnorm{f_{\lambda}}{s}^{2} \, d\lambda$. This is a slight simplification
(a precise description is given below), but for most practical purposes it
will suffice.

An important consequence of \eqref{SimplifiedIntegralFormula} is the
following:
\begin{principle}\label{WaveEstimatesImpliesSpaceTimeEstimates} A linear
or multilinear space-time estimate for solutions of the homogeneous wave
equation with data in $H^s$ implies a corresponding estimate for elements
of $H^{s,\theta}$.
\end{principle}
This is made precise in Proposition \ref{EmbeddingCorollary} below. To
illustrate this principle, let us interpret Theorem A in terms of
$H^{s,\theta}$. Assume that the hypotheses of Theorem A are satisfied, and
take the $L_{t}^{q}(L_{x}^{r})$-norm in \eqref{SimplifiedIntegralFormula}.
By Minkowski's integral inequality,
$$
\mixednorm{u}{q}{r} \lesssim \int_{-\infty}^{\infty}
\frac{\mixednorm{u_{\lambda}}{q}{r}}{(1+\abs{\lambda})^{\theta}} d\lambda,
$$
and by Theorem A, $\mixednorm{u_{\lambda}}{q}{r} \lesssim
\Sobnorm{f_{\lambda}}{s}$. Thus, if $\theta > 1/2$, $$
\int_{-\infty}^{\infty}
\frac{\Sobnorm{f_{\lambda}}{s}}{(1+\abs{\lambda})^{\theta}} d\lambda
\le C_{\theta} \left( \int \Sobnorm{f_{\lambda}}{s}^{2} d\lambda
\right)^{\half} = C_{\theta} \spacetimenorm{u}{s}{\theta}, $$
whence $\mixednorm{u}{q}{r} \lesssim \spacetimenorm{u}{s}{\theta}$. We
summarize: %
\begin{ThmD}
The embedding
$$
H^{\frac{n}{2}-\frac{n}{r}-\frac{1}{q},\theta} \hookrightarrow
L_{t}^{q}(L_{x}^{r})
$$
holds whenever $(q,r)$ is wave admissible and $\theta > 1/2$. \end{ThmD}
Theorem D may be viewed as an analog for $H^{s,\theta}$ of the Sobolev
embedding \eqref{SobolevEmbedding1} for the standard Sobolev spaces.

Just as in the linear case, via \eqref{SimplifiedIntegralFormula} we can
interpret the bilinear estimates of Theorem D in $H^{s,\theta}$. %
\begin{ThmE}
If $n \ge 2$ and $q,r,s_1,s_2$ and $\sigma$ satisfy the hypotheses of
Theorem D, then
$$
\mixednorm{D^{-\sigma}(uv)}{q/2}{r/2}
\lesssim
\spacetimenorm{u}{s_1}{\theta}
\spacetimenorm{v}{s_2}{\theta}
$$
provided $\theta > 1/2$.
\end{ThmE}
The crucial observation is that since $D^{-\sigma}$ does not involve the
time variable, the integral formula \eqref{SimplifiedIntegralFormula}
implies $$
D^{-\sigma}(u^2)(t) = \frac{1}{4\pi^2} \int \!\! \int
\frac{e^{it(\lambda_1+\lambda_2)} D^{-\sigma} ( u_{\lambda_1}
u_{\lambda_2}) (t)} {(1+\abs{\lambda_1})^{\theta}
(1+\abs{\lambda_2})^{\theta}} d\lambda_1 \, d\lambda_2.
$$
Take the $L_t^{q/2}(L_x^{r/2})$-norm, use Minkowski's integral inequality,
Theorem D and finally the Cauchy-Schwarz inequality to obtain the estimate
in Theorem E (see Proposition \ref{EmbeddingCorollary} and Remark
\ref{HalfWaveRemark} for the details).

Theorem C, in contrast to Theorem B, does not have an obvious
interpretation in terms of $H^{s,\theta}$ via the integral representation
\eqref{SimplifiedIntegralFormula}, since the operator $(u,v) \mapsto
D^{\gamma} D_{+}^{\gamma_{+}} D_{-}^{\gamma_{-}}(uv)$ acts in both space
and time. Nevertheless, if we set $\gamma_{+} = 0$, Theorem C does have an
$H^{s,\theta}$-analog, but with $D_{-}^{\gamma_{-}}$ replaced by the
operator $R^{\gamma_{-}}$ appearing in the following lemma. %
\begin{lemma}\label{HyperbolicLambdaEstimates} If $\alpha > 0$, then $$
D_{-}^{\alpha} (uv)
\precsim (D_{-}^{\alpha}u) v + u D_{-}^{\alpha}v + R^{\alpha}( u , v ),
$$
for all $u$ and $v$ with nonnegative Fourier transforms, where
$R^{\alpha}$ is the symmetric bilinear operator given by \begin{align*}
\Fourier R^{\alpha}(u,v) (\Xi)
&= \int_{\R^{1+n}} \bigl[ r(\Xi-\Xi';\Xi') \bigr]^{\alpha}
\widehat{u}(\Xi-\Xi') \widehat{v}(\Xi')
\, d\Xi',
\\
r(\tau,\xi;\lambda,\eta) &=
\begin{cases}
\abs{\xi} + \abs{\eta} - \abs{\xi+\eta} &\text{if} \quad \tau\lambda \ge
0, \\
\abs{\xi+\eta} - \bigl\vert \abs{\xi} - \abs{\eta} \bigr\vert &\text{if}
\quad \tau\lambda < 0.
\end{cases}
\end{align*}
Moreover, the same estimate holds with $D_{-}^{\alpha}$ replaced by
$\Lambda_{-}^{\alpha}$.
\end{lemma}
\begin{proof} It is enough to show that
$$
\hypwt{\tau+\lambda}{\xi+\eta}
\le \hypwt{\tau}{\xi} + \hypwt{\lambda}{\eta} + r(\tau,\xi,\lambda,\eta).
$$
The proof splits into four cases, corresponding to the four quadrants of
the $(\tau,\lambda)$-plane. For example, if $\tau, \lambda \ge 0$, then
\begin{align*}
\hypwt{\tau+\lambda}{\xi+\eta}
&= \bigabs{ \tau - \abs{\xi} + \lambda - \abs{\eta} + \abs{\xi} +
\abs{\eta} - \abs{\xi + \eta} } \\ &\le \bigabs{ \tau - \abs{\xi} } +
\bigabs{ \lambda - \abs{\eta} } + \abs{\xi} + \abs{\eta} - \abs{\xi +
\eta}, \end{align*}
and the remaining cases are similar. Finally, note that the symbol of
$\Lambda_-$ is comparable to $1 + \hypwt{\tau}{\xi}$.
\end{proof}
We now state the $H^{s,\theta}$-version of Theorem C. %
\begin{ThmF} If $n \ge 2$ and $\gamma,\gamma_{-},s_{1}$ and $s_{2}$
satisfy the hypotheses of Theorem C with $\gamma_{+} = 0$, then $$
\twonorm{D^{\gamma} R^{\gamma_{-}}(u,v)}{(\R^{1+n})} \lesssim
\spacetimenorm{D^{s_{1}}u}{0}{\theta} \spacetimenorm{D^{s_{2}}v}{0}{\theta}
$$
provided $\theta > 1/2$.
\end{ThmF}
To prove this result we need the precise version of the integral
representation \eqref{SimplifiedIntegralFormula}. First note that any $u
\in H^{s,\theta}$ has a unique decomposition
\begin{equation}\label{PlusMinusDecomposition}
u = u_{+} + u_{-}
\end{equation}
where $u_{+}$ and $u_{-}$ belong to $H^{s,\theta}$ and have Fourier
transforms supported in $[0,\infty) \times \R^{n}$ and $(-\infty,0] \times
\R^{n}$ respectively. Moreover,
$$
\spacetimenorm{u}{s}{\theta}^{2}
= \spacetimenorm{u_{+}}{s}{\theta}^{2} +
\spacetimenorm{u_{-}}{s}{\theta}^{2}.
$$
The notation in the decomposition \eqref{PlusMinusDecomposition} is
intentionally the same as the one used in the decomposition of a solution
of the homogeneous wave equation into half-waves (see section
\ref{Notation}).
\begin{proposition}\label{IntegralRepresentationProposition} If $u \in
H^{s,\theta}$, there exist $f_{+}, f_{-} \in L^{2}( \R , H^{s} )$ such that
\begin{equation}\label{IntegralRepresentation}
u_{\pm}(t) = \frac{1}{2\pi}\int_{-\infty}^{\infty} \frac{e^{it(\lambda \pm
D)} f_{\pm}(\lambda)} {( 1 + \abs{\lambda} )^{\theta} } \, d\lambda \qquad
\text{($H^{s}$-valued)}
\end{equation}
and $\spacetimenorm{u_{\pm}}{s}{\theta} = \norm{f_{\pm}}_{L^{2}( \R ,
H^{s} )}$. \end{proposition}
Thus, elements of $H^{s,\theta}$ may be thought of as superpositions of
half-waves with data in $H^{s}$. In fact, $f_{+}$ is given by $$
\Fourier\{ f_{+}(\lambda) \}(\xi) =
\begin{cases}
(1+\abs{\lambda})^{\theta}
\widehat{u}(\lambda+\abs{\xi},\xi) &\text{if $\lambda + \abs{\xi} \ge 0$,}
\\
0 &\text{otherwise,}
\end{cases}
$$
and $f_{-}$ has a similar definition. It is then easy to verify
\eqref{IntegralRepresentation} by applying the spacetime Fourier transform
(see \cite{Se}).

The following is the precise statement of Principle
\ref{WaveEstimatesImpliesSpaceTimeEstimates}. %
\begin{proposition}\label{EmbeddingCorollary} Assume that
$T : H^{s_{1}}(\R^{n}) \times \cdots \times H^{s_{k}}(\R^{n})
\longrightarrow H^{\sigma}(\R^{n})$
is k-linear, and let $\theta > 1/2$.
{\renewcommand{\theenumi}{\alph{enumi}}
\renewcommand{\labelenumi}{(\theenumi)}
\begin{enumerate}
\item Fix a $k$-tuple $\varepsilon \in \{-1,1\}^{k}$. If
\begin{equation}\label{HalfWaveEstimate}
\mixednorm{ T(e^{\varepsilon_{1} itD} f_{1},\dots, e^{\varepsilon_{k} itD}
f_{k}) }{q}{r} \lesssim \Sobnorm{f_{1}}{s_{1}} \cdots
\Sobnorm{f_{k}}{s_{k}}, \end{equation}
for all $(f_{1},\dots,f_{k}) \in H^{s_{1}} \times \cdots \times
H^{s_{k}}$, then \begin{equation}\label{MultilinearSpaceTimeEstimate}
\mixednorm{ T\bigl(u_{1}(t),\dots,u_{k}(t)\bigr)(x) }{q}{r} \le C_{\theta}
\spacetimenorm{u_{1}}{s_{1}}{\theta}
\cdots \spacetimenorm{u_{k}}{s_{k}}{\theta} \end{equation}
for all $(u_{1},\dots,u_{k}) \in
H^{s_{1},\theta}\times \cdots \times H^{s_{k},\theta}$ such that
\begin{equation}\label{EmbeddingCorollarySupport}
\supp \widehat{u_{j}} \subseteq
\begin{cases}
[0,\infty) \times \R^{n} &\text{if} \quad \varepsilon_{j} = 1, \\
(-\infty,0] \times \R^{n} &\text{if} \quad \varepsilon_{j} = -1.
\end{cases}
\end{equation}
\item If \eqref{HalfWaveEstimate} holds for all $\varepsilon \in
\{-1,1\}^{k}$ and all $(f_{1},\dots,f_{k}) \in H^{s_{1}} \times \cdots
\times H^{s_{k}}$,
then \eqref{MultilinearSpaceTimeEstimate} holds for all
$(u_{1},\dots,u_{k}) \in
H^{s_{1},\theta}\times \cdots \times H^{s_{k},\theta}$. \end{enumerate}}
\end{proposition}
\begin{proof}
By proposition \ref{IntegralRepresentationProposition} and the condition
\eqref{EmbeddingCorollarySupport}, which is equivalent to $$
u_{j} =
\begin{cases}
u_{j+} &\text{if} \quad \varepsilon_{j} = 1, \\ u_{j-} &\text{if} \quad
\varepsilon_{j} = -1, \end{cases}
$$
there exist $f_{j} \in L^{2}( \R , H^{s_{j}} )$ for $j = 1,\dots,k$ such
that $$
u_{j} =\int_{-\infty}^{\infty}
\frac{e^{it\lambda} e^{\varepsilon_{j} itD}f_{j}(\lambda)} {( 1 +
\abs{\lambda} )^{\theta} } \, d\lambda $$
and $\spacetimenorm{u_{j}}{s_{j}}{\theta} =
\norm{f_{j}}_{L^{2}(\R,H^{s_{j}})}$. By linearity, \begin{align*}
&T\bigl( u_{1}(t),\dots,u_{k}(t) \bigr) \\ &\qquad
= \int_{-\infty}^{\infty} \! \cdots \!
\int_{-\infty}^{\infty}
\frac{ e^{it(\lambda_{1} + \cdots + \lambda_{k})} T\bigl(
e^{\varepsilon_{1} itD} f_{1}(\lambda_{1}) , \dots , e^{\varepsilon_{k}
itD} f_{k}(\lambda_{k}) \bigr)} {( 1 + \abs{\lambda_{1}} )^{\theta} \cdots
( 1 + \abs{\lambda_{k}} )^{\theta}} d\lambda_{1} \cdots d\lambda_{k},
\end{align*}
so by Minkowski's integral inequality,
\eqref{HalfWaveEstimate} and the Cauchy-Schwarz inequality, \begin{align*}
&\mixednorm{ T(u_{1},\dots,u_{k}) }{q}{r} \\ &\qquad
\le \int_{-\infty}^{\infty} \! \cdots \! \int_{-\infty}^{\infty}
\frac{ C \Sobnorm{f_{1}(\lambda_{1})}{s_{1}} \cdots
\Sobnorm{f_{k}(\lambda_{k})}{s_{k}} }
{( 1 + \abs{\lambda_{1}} )^{\theta} \cdots ( 1 + \abs{\lambda_{k}}
)^{\theta}} d\lambda_{1} \cdots d\lambda_{k} \\ &\qquad
\le C \norm{f_{1}}_{L^{2}(\R,H^{s_{1}})} \cdots
\norm{f_{k}}_{L^{2}(\R,H^{s_{k}})}. \end{align*}
This concludes the proof of part (a), and to prove part (b) we simply
write $u_{j} = u_{j+} + u_{j-}$, use the multilinearity of $T$, and apply
part (a).
\end{proof}
\begin{remark}\label{HalfWaveRemark} Theorems A,B and C remain true with
$u$ and $v$ replaced by any of their half-waves (in fact, this is how the
estimates are proved). Thus, part (b) of Proposition
\ref{EmbeddingCorollary}, applied to Theorems A and B, proves Theorems D
and E respectively. Notice also that Theorems D and E remain true when
$H^{s,\theta}$ is replaced by the space with norm $\twonorm{D^{s}
\Lambda_{-}^{\theta} u}{}$. The reason is that any estimate of the form
\eqref{GeneralBilinearEstimate} may be rewritten as follows:
$$
\mixednorm{D^{\gamma} D_{+}^{\gamma_{+}} D_{-}^{\gamma_{-}}(D^{-s_{1}}u
D^{-s_{2}}v)}{q/2}{r/2} \lesssim \twonorm{f}{} \twonorm{g}{}.
$$
This fact will be used freely in the rest of the paper. \end{remark}
\begin{definition}\label{SDefinition}
Let $S_+^\alpha$ and $S_{-}^{\alpha}$ be the symmetric bilinear operators
given by \begin{align*}
\Fourier S_{+}^{\alpha} (f,g) (\xi)
&= \int_{\R^{n}} \bigl( \abs{\xi-\eta} + \abs{\eta}
- \abs{\xi} \bigr)^{\alpha} \widehat{f}(\xi-\eta) \widehat{g}(\eta) \,
d\eta,
\\
\Fourier S_{-}^{\alpha} (f,g) (\xi)
&= \int_{\R^{n}} \bigl( \abs{\xi} - \bigabs{\abs{\xi-\eta} - \abs{\eta}}
\bigr)^{\alpha} \widehat{f}(\xi-\eta)
\widehat{g}(\eta) \, d\eta.
\end{align*}
\end{definition}
The relation between the above operators and $D_{-}^{\alpha}$ and
$R^{\alpha}$ is given in the following lemma. Keep in mind that the latter
two operators act on functions defined on the space-time. If $u,v$ are two
such functions, then $S(u,v)$ denotes the function $(t,x) \mapsto S\bigl(
u(t), v(t) \bigr) (x)$, where as usual we write $u(t) = u(t,\cdot)$.
\begin{lemma}\label{LemmaForS} Let $f,g$ be defined on $\R^{n}$ and $u,v$
on $\R^{1+n}$. Then
{\renewcommand{\theenumi}{\alph{enumi}}
\renewcommand{\labelenumi}{(\theenumi)}
\begin{enumerate}
\item
$D_{-}^{\alpha} (e^{itD} f \cdot e^{\pm itD} g) = S_{\pm}^{\alpha}
(e^{itD} f,e^{\pm itD} g)$, %
\item
$R^{\alpha}(u,v) = S_{+}^{\alpha}(u_{+},v_{+}) +
S_{-}^{\alpha}(u_{+},v_{-}) + S_{-}^{\alpha}(u_{-},v_{+}) +
S_{+}^{\alpha}(u_{-},v_{-})$.
\end{enumerate}
}
\end{lemma}
\begin{proof} We prove (a) for $S_{+}$; the proof for $S_{-}$ is similar.
We have
$$
\Fourier (e^{itD} f \cdot e^{\pm itD} g) (\tau,\xi) = \int_{\R^{n}}
\delta( \tau - \abs{\xi-\eta} - \abs{\eta} ) \widehat{f}(\xi-\eta)
\widehat{g}(\eta) \, d\eta, $$
whence
\begin{align*}
&\Fourier \bigl\{ D_{-}^{\alpha} (e^{itD} f \cdot e^{\pm itD} g) \bigr\}
(\tau,\xi)
\\
&\qquad \quad = \hypwt{\tau}{\xi}^{\alpha} \int_{\R^{n}} \delta( \tau -
\abs{\xi-\eta} - \abs{\eta} ) \widehat{f}(\xi-\eta) \widehat{g}(\eta) \,
d\eta \\
&\qquad \quad = \int_{\R^{n}} \delta( \tau - \abs{\xi-\eta} - \abs{\eta} )
\bigl( \abs{\xi-\eta} + \abs{\eta} - \abs{\xi} \bigr)^{\alpha}
\widehat{f}(\xi-\eta) \widehat{g}(\eta) \, d\eta. \end{align*}
But the last expression equals
$$
\Fourier \bigl\{S_{\pm}^{\alpha} (e^{itD} f,e^{\pm itD} g) \bigr\}
(\tau,\xi). $$

As for (b), using the decomposition \eqref{PlusMinusDecomposition} and the
bilinearity of $R^{\alpha}$, we have $$
R^{\alpha}(u,v) =
R^{\alpha}(u_{+},v_{+}) + R^{\alpha}(u_{+},v_{-}) +
R^{\alpha}(u_{-},v_{+}) + R^{\alpha}(u_{-},v_{-}), $$
so it suffices to prove $R^{\alpha}(u_{+},v_{\pm}) =
S_{\pm}^{\alpha}(u_{+},v_{\pm})$ and
$R^{\alpha}(u_{-},v_{\mp}) = S_{\pm}^{\alpha}(u_{-},v_{\mp})$. We will
only prove the case
\begin{equation}\label{LemmaForSEqA}
R^{\alpha}(u_{+},v_{+}) =
S_{+}^{\alpha}(u_{+},v_{+}).
\end{equation}
The Fourier transform of the right hand side, at fixed $t$, is $$
\int_{\R^{n}} \bigl( \abs{\xi-\eta} + \abs{\eta} - \abs{\xi}
\bigr)^{\alpha} \widehat{u_{+}(t)}(\xi-\eta) \widehat{v_{+}(t)}(\eta) \,
d\eta,
$$
and applying the Fourier transform in $t$ yields the following expression
for the space-time Fourier transform at $(\tau,\xi)$ of the right hand side
of \eqref{LemmaForSEqA}: $$
\int_{\R^{n}} \! \! \int_{\R}
\bigl( \abs{\xi-\eta} + \abs{\eta} - \abs{\xi} \bigr)^{\alpha}
\widehat{u_{+}}(\tau-\lambda,\xi-\eta) \widehat{v_{+}}(\lambda,\eta) \,
d\lambda \, d\eta.
$$
But the latter is, by the definition of $R^{\alpha}$, equal to the
space-time Fourier transform of the left hand side of \eqref{LemmaForSEqA}.
The remaining cases are proved in a similar manner. \end{proof}
We can now prove Theorem F. If the hypotheses of the theorem are
satisfied, then by Lemma \ref{LemmaForS}(a) and Theorem C, we have (cf.
remark \ref{HalfWaveRemark}) $$
\bigtwonorm{D^{\gamma} S_{\pm}^{\gamma_{-}} (e^{itD} D^{-s_{1}} f, e^{\pm
itD} D^{-s_{2}} g)}{} \lesssim \twonorm{f}{} \twonorm{g}{}, $$
for all $f,g \in L^{2}(\R^{n})$, so by Proposition
\ref{EmbeddingCorollary}(a), \begin{equation}\label{ThmFProofA}
\bigtwonorm{D^{\gamma} S_{\pm}^{\gamma_{-}} ( D^{-s_{1}} u_{+}, D^{-s_{2}}
v_{\pm})}{} \lesssim \spacetimenorm{u_{+}}{0}{\theta}
\spacetimenorm{v_{\pm}}{0}{\theta}
\end{equation}
for all $u,v \in H^{0,\theta}$, where $\theta > 1/2$. A similar argument
gives
\begin{equation}\label{ThmFProofB}
\bigtwonorm{D^{\gamma} S_{\pm}^{\gamma_{-}} ( D^{-s_{1}} u_{-}, D^{-s_{2}}
v_{\mp})}{} \lesssim \spacetimenorm{u_{-}}{0}{\theta}
\spacetimenorm{v_{\mp}}{0}{\theta}.
\end{equation}
By \eqref{ThmFProofA}, \eqref{ThmFProofB} and Lemma \ref{LemmaForS}(b), we
conclude that
$$
\bigtwonorm{D^{\gamma} R^{\gamma_{-}} (u,v)}{} \lesssim
\spacetimenorm{D^{s_{1}} u}{0}{\theta}
\spacetimenorm{D^{s_{2}} v}{0}{\theta},
$$
so Theorem F is proved.

The following embedding is an easy consequence of the integral formula
\eqref{IntegralRepresentation} and the dominated convergence theorem for
$H^{s}$-valued integrals. We omit the details. Here $C_{b}(\R,H^{s})$
denotes the space of bounded, continuous maps from $\R$ to $H^{s}$, with
the supremum norm.
\begin{proposition}\label{EnergyEmbeddingProposition}
$H^{s,\theta} \hookrightarrow C_b(\R,H^s)$ if $\theta > 1/2$.
\end{proposition}
Associated to $H^{s,\theta}$ we have the following space. %
\begin{definition} For $s, \theta \in \R$, define $$
\scrH^{s,\theta} = \left\{ u : \text{$u \in H^{s,\theta}$ and
$\partial_{t}u \in H^{s-1,\theta}$} \right\} $$
with norm $\Spacetimenorm{u}{s}{\theta} = \spacetimenorm{u}{s}{\theta} +
\spacetimenorm{\partial_{t}u}{s-1}{\theta}$. \end{definition}
\begin{remark} An equivalent, but less intuitive definition is $$
\scrH^{s,\theta} = \left\{ u \in \Schwartz' : \Lambda^{s-1} \Lambda_{+}
\Lambda_{-}^{\theta}u \in L^{2} \right\} $$
with norm $\Spacetimenorm{u}{s}{\theta} = \twonorm{\Lambda^{s-1}
\Lambda_{+} \Lambda_{-}^{\theta}u}{}$. We will use these two definitions
of $\scrH^{s,\theta}$ and its norm interchangeably. \end{remark}
The following embedding is a corollary to Proposition
\ref{EnergyEmbeddingProposition}.
\begin{proposition}\label{SecondEnergyEmbeddingProposition}
$\scrH^{s,\theta} \hookrightarrow C_b(\R,H^s) \cap C^1_b(\R,H^{s-1})$ if
$\theta > \half$.
\end{proposition}
\section{The Space $\Mixed{q}{r}$}\label{MixedSpace}
Optimal local well-posedness for \eqref{WMType} (part (a) of the Main
Theorem) will be proved by iteration in the space $\scrH^{s,\theta}$,
defined in the previous section. An attempt to prove the corresponding
results for \eqref{MKGType}/\eqref{YMType} (part (b) of the Main Theorem)
and \eqref{CFWMType} (part (c) of the Main Theorem) by iteration in the
same space, leads to estimates which are false. It turns out that the
iteration works out if we replace $\scrH^{s,\theta}$ with the subspace
defined by a norm $$
\norm{u} = \Spacetimenorm{u}{s}{\theta} + \mixednorm{\Lambda^{\gamma}
\Lambda_{-}^{\gamma_{-}} u}{q}{r}, $$
where the choice of exponents $\gamma, \gamma_{-}, q$ and $r$ is dictated
by the specific equation under consideration.

However, since we want a space whose norm only depends on the size of the
Fourier transform, the space $L_{t}^{q}(L_{x}^{r})$ must be modified. To
motivate the following definition, recall that if $u \in \Schwartz'$ and
$1 \le q,r \le \infty$, then
$$
\mixednorm{u}{q}{r} = \sup_{v} \abs{\innerprod{u}{v}} = \sup_{v}
\abs{\innerprod{\widehat u}{\widehat v}}, $$
where the supremum is over all $v \in \Schwartz$ such that
$\mixednorm{v}{q'}{r'} = 1$ ($q'$ and $r'$ being the dual exponents of $q$
and $r$ respectively). Of course, if $\widehat u$ is a tempered function,
then $\innerprod{\widehat u}{\widehat v} = \int_{\R^{1+n}} \widehat u
\widehat v$. %
\begin{definition}\label{ModifiedMixednormDefinition} If $1 \le q, r \le
\infty$, $u \in \Schwartz'$ and $\widehat u$ is a tempered function, set $$
\Mixednorm{u}{q}{r} = \sup \left\{ \int_{\R^{1+n}} \abs{ \widehat u (\Xi)
} \widehat v(\Xi) \, d\Xi : \text{$v \in \Schwartz$, $\widehat v \ge 0$
and $\mixednorm{v}{q'}{r'} = 1$} \right\}, $$
where $q'$ and $r'$ are the conjugate exponents of $q$ and $r$
respectively; i.e., $1 = \frac{1}{q} + \frac{1}{q'}$ and $1 = \frac{1}{r}
+ \frac{1}{r'}$. Let $\Mixed{q}{r}$ be the corresponding subspace of
$\Schwartz'$. \end{definition}
Clearly, $\Mixednorm{\cdot}{q}{r}$ is a translation invariant norm on
$\Mixed{q}{r}$, it is compatible with the relation $\preceq$, and it only
depends on the size of the Fourier transform. Note that $\Mixed{2}{2} =
L^2(\R^{1+n})$. Observe also that
\begin{equation}\label{ModifiedMixednormBound}
\Mixednorm{u}{q}{r} \le \mixednorm{u}{q}{r} \quad \text{whenever} \quad
\widehat u \ge 0. \end{equation}

The above definition is inspired by the norms introduced in \cite{Kl-Ma4,
Kl-Ma5}. Another way of modifying the norm on $L_t^q(L_x^r)$ so that it
only depends on the size of the Fourier transform can be found in
\cite{Kl-Ta}. %
\begin{proposition}\label{CompletenessOfXs} Let
$s,\gamma,\gamma_+,\gamma_- \in \R$, $\theta > \half$ and $1 \le q,r \le
\infty$. Define
$$
\X^s = \bigl\{ u : \norm{u} < \infty \bigr\}, $$
where
$$
\norm{u} = \Spacetimenorm{u}{s}{\theta} + \Mixednorm{\Lambda^\gamma
\Lambda_+^{\gamma_+} \Lambda_-^{\gamma_-} u}{q}{r}.
$$
Then $\X^s$ is a Banach space.
\end{proposition}
\begin{proof}
Assume that $(u_{j})$ is a Cauchy sequence in $\X^{s}$. Then $(u_{j})$ is
Cauchy in $\scrH^{s,\theta}$, so it converges in the latter space to some
limit $u$. It remains to prove that $\norm{u_{j}-u} \to 0$ as $j \to
\infty$. Fix $\varepsilon > 0$. There exists $M \in \N$ such that $$
\bigMixednorm{\Lambda^\gamma \Lambda_{+}^{\gamma_{+}}
\Lambda_{-}^{\gamma_-} (u_{j}-u_{k})}{q}{r} \le \varepsilon
$$
for all $j,k \ge M$. We claim that
$$
\bigMixednorm{\Lambda^\gamma \Lambda_{+}^{\gamma_{+}}
\Lambda_{-}^{\gamma_-} (u_{j}-u)}{q}{r} \le \varepsilon
$$
for all $j \ge M$. To see this, fix $v \in \Schwartz$ such that
$\widehat{v} \ge 0$ and $\mixednorm{v}{q'}{r'} = 1$, where $q', r'$ are
conjugate to $q, r$. Then
$$
\int \abs{ \Fourier \left\{ \Lambda^\gamma \Lambda_{+}^{\gamma_{+}}
\Lambda_{-}^{\gamma_-} (u_{j}-u_{k}) \right\}(\Xi) } \widehat v(\Xi) \,
d\Xi \le \varepsilon $$
for all $j,k \ge M$, so it suffices to prove that
\begin{multline}\label{CompletenessProofA}
\lim_{k \to \infty} \int \abs{ \Fourier \left\{ \Lambda^\gamma
\Lambda_{+}^{\gamma_{+}} \Lambda_{-}^{\gamma_-} (u_{j}-u_{k})
\right\}(\Xi) } \widehat v(\Xi) \, d\Xi \\
= \int \abs{ \Fourier \left\{ \Lambda^\gamma \Lambda_{+}^{\gamma_{+}}
\Lambda_{-}^{\gamma_-} (u_{j}-u) \right\}(\Xi) } \widehat v(\Xi) \, d\Xi
\end{multline}
for fixed $j$. To prove this, we write
\begin{multline*}
\int \abs{ \Fourier \left\{ \Lambda^\gamma \Lambda_{+}^{\gamma_{+}}
\Lambda_{-}^{\gamma_-} (u_{j}-u_{k}) \right\}(\Xi) } \widehat v(\Xi) \,
d\Xi \\
= \int \abs{ \Fourier \left\{ \Lambda^{s-1} \Lambda_{+}
\Lambda_{-}^{\theta} (u_{j}-u_{k}) \right\}(\Xi) } \widehat{v'}(\Xi) \,
d\Xi,
\end{multline*}
where $v' = \Lambda^{\gamma+1-s} \Lambda_{+}^{\gamma_{+}-1}
\Lambda_{-}^{\gamma_- - \theta} v \in \Schwartz$. Since $\Lambda^{s-1}
\Lambda_{+} \Lambda_{-}^{\theta} (u_{j}-u_{k})$ converges to
$\Lambda^{s-1} \Lambda_{+} \Lambda_{-}^{\theta} (u_{j}-u)$ in $L^{2}$, we
conclude that \eqref{CompletenessProofA} holds.
\end{proof}
For later use, we mention some basic properties of $\Mixed{q}{r}$.

First, a version of H\"older's inequality holds. %
\begin{proposition}\label{MixednormHolder} Suppose $\frac{1}{q} =
\frac{1}{q_{1}} + \frac{1}{q_{2}}$ and $\frac{1}{r} = \frac{1}{r_{1}} +
\frac{1}{r_{2}}$, where the $q$'s and $r$'s all belong to $[1,\infty]$.
Then $$
\Mixednorm{uv}{q}{r} \le \Mixednorm{u}{q_{1}}{r_{1}}
\mixednorm{v}{q_{2}}{r_{2}}
$$
for all $u \in \Mixed{q_{1}}{r_{1}}$ and $v \in \Schwartz$ with $\widehat
v \ge 0$.
\end{proposition}
\begin{proof}
Since the norm $\Mixednorm{\cdot}{\cdot}{\cdot}$ is compatible with the
relation $\preceq$ and only depends on the size of the Fourier transform,
it suffices to prove the inequality when $\widehat u \ge 0$. Thus, we fix
$u$ and $v$ such that $v \in \Schwartz$ and $\widehat u, \widehat v \ge
0$. Let $q', r', q_{1}'$ etc. denote the dual exponents. Then
$$
\Mixednorm{uv}{q}{r}
= \sup_{w} \int uvw \, dt \, dx
\le \Mixednorm{u}{q_{1}}{r_{1}}
\mixednorm{vw}{q_{1}'}{r_{1}'},
$$
where the supremum is over all $w \in \Schwartz$ such that $\widehat w \ge
0$ and $\mixednorm{w}{q'}{r'} = 1$. But by H\"older's inequality,
$$
\mixednorm{vw}{q_{1}'}{r_{1}'}
\le \mixednorm{v}{q_{2}}{r_{2}} \mixednorm{w}{q'}{r'}, $$
finishing the proof.
\end{proof}
When applying the previous proposition, the following is useful: %
\begin{lemma}\label{DensityLemma}
$\{ u \in \Schwartz : \widehat u \ge 0 \}$ is dense in $\{ u \in
H^{a,\alpha} : \widehat u \ge 0 \}$ for all $a,\alpha \in \R$.
\end{lemma}
\begin{proof}
Since $\Lambda^{a} \Lambda_{-}^{\alpha}$ is an isomorphism of
$H^{a,\alpha}$ onto $L^{2}$, which preserves the Schwartz class and
positivity of the Fourier transform, we may take $a = \alpha = 0$. Since
the Fourier transform is an isomorphism of both $\Schwartz$ and $L^{2}$,
it then suffices to prove that the set $\{ v \in \Schwartz : v \ge 0 \}$
is dense in $\{ v \in L^{2} : v \ge 0 \}$. But the standard proof that
$\Schwartz$ is dense in $L^{2}$ shows this to be true.
\end{proof}
The following duality argument is fundamental to our approach. %
\begin{proposition}\label{DualityProposition} Let $1 \le a,b,q,r \le
\infty$.
{
\renewcommand{\theenumi}{\alph{enumi}}
\renewcommand{\labelenumi}{(\theenumi)}
\begin{enumerate}
\item If
\begin{equation}\label{DualityA}
\mixednorm{G}{a'}{b'} \lesssim \bigmixednorm{\Lambda^{\alpha}
\Lambda_{-}^{\beta} G}{q'}{r'}
\end{equation}
for all $G \in \Schwartz$, then
$$
\mixednorm{F}{q}{r} \lesssim \bigmixednorm{\Lambda^{\alpha}
\Lambda_{-}^{\beta} F}{a}{b}
$$
for all $F$.
\item If \eqref{DualityA} holds for all $G \in \Schwartz$ with $\widehat G
\ge 0$, then
$$
\Mixednorm{F}{q}{r} \lesssim \bigMixednorm{\Lambda^{\alpha}
\Lambda_{-}^{\beta} F}{a}{b}
$$
for all $F$.
\end{enumerate}
}
\end{proposition}
\begin{proof}
We have
\begin{align*}
\mixednorm{F}{q}{r} &= \sup_{G} \Bigabs{ \int FG} \\
&= \sup_{G} \Bigabs{ \int \Lambda^{\alpha}\Lambda_{-}^{\beta} F
\Lambda^{-\alpha}\Lambda_{-}^{-\beta} G} \\
&\le \bigmixednorm{\Lambda^{\alpha}\Lambda_{-}^{\beta} F}{a}{b} \sup_{G}
\bigmixednorm{\Lambda^{-\alpha}\Lambda_{-}^{-\beta} G}{a'}{b'},
\end{align*}
where the supremum is over all $G \in \Schwartz$ with
$\mixednorm{G}{q'}{r'} = 1$. Part (a) follows.

For part (b), we have
\begin{align*}
\Mixednorm{F}{q}{r} &= \sup_{G} \int \bigabs{\widehat F} \widehat G \\
&= \sup_{G} \int \bigabs{\Fourier \bigl(
\Lambda^{\alpha}\Lambda_{-}^{\beta} F \bigr)} \Fourier \bigl(
\Lambda^{-\alpha}\Lambda_{-}^{-\beta} G \bigr) \\
&\le \bigMixednorm{\Lambda^{\alpha}\Lambda_{-}^{\beta} F}{a}{b} \sup_{G}
\bigmixednorm{\Lambda^{-\alpha}\Lambda_{-}^{-\beta} G}{a'}{b'},
\end{align*}
where the supremum is over all $G \in \Schwartz$ such that $\widehat G \ge
0$ and $\mixednorm{G}{q'}{r'} = 1$. \end{proof}
\begin{corollary}\label{DualityCorollary} Let $1 \le a,b,q,r \le \infty$.
If
\begin{align*}
\Lambda^{-\alpha} \Lambda_{-}^{-\beta} \mixed{a}{b} &\hookrightarrow
\mixed{q}{r},
\intertext{then}
\Lambda^{-\alpha} \Lambda_{-}^{-\beta} \Mixed{a}{b} &\hookrightarrow
\Mixed{q}{r}.
\end{align*}
\end{corollary}
\begin{proof} We apply Proposition \ref{DualityProposition}. By part (a),
we get
\begin{align*}
\Lambda^{-\alpha} \Lambda_{-}^{-\beta} \mixed{q'}{r'} &\hookrightarrow
\mixed{a'}{b'}.
\intertext{Then, by part (b),}
\Lambda^{-\alpha} \Lambda_{-}^{-\beta} \Mixed{a}{b} &\hookrightarrow
\Mixed{q}{r}.
\end{align*}
\end{proof}
The next estimate is used in the proof of part (b) of the Main Theorem.
The proof is based on an idea from \cite{Kl-Ma5}. %
\begin{proposition}\label{KlMaEmbedding} If $n \ge 4$, $1 < p \le
\frac{2(n-1)}{n+1}$, $s = \frac{n}{p} - \frac{n}{2} - \half$ and $\theta >
\half$, then
$$
\mixednorm{D^{-s} \Lambda_{-}^{-\theta} u}{\infty}{2} \lesssim
\mixednorm{u}{\infty}{p}
$$
for all $u \in \Schwartz$ with nonnegative Fourier transform.
\end{proposition}
\begin{proof}
Set $U = D^{-s} \Lambda_{-}^{-\theta} u$. Since $\mixednorm{U}{\infty}{2}
\le \bigtwonorm{\int \widehat{U}(\tau,\xi) \, d\tau}{_{\xi}}$, it suffices
to show that
$$
\int_{\R^{1+n}} \abs{\xi}^{-s} \bigl(1 + \bigabs{\abs{\tau} - \abs{\xi}}
\bigr)^{-\theta}
\widehat{u}(\tau,\xi) \Check f(\xi) \, d\tau \, d\xi \lesssim
\mixednorm{u}{\infty}{p} \twonorm{f}{} $$
for all $f \in \Schwartz(\R^{n})$ whose inverse Fourier transform $\Check
f$ is nonnegative. The integral on the left hand side is dominated by
\begin{equation}\label{KlMaEmbeddingA}
\int_{\R^{1+n}} u (v_{+} + v_{-}) \, dt \, dx, \end{equation}
where
$$
\Fourier^{-1} v_\pm(\tau,\xi) = \abs{\xi}^{-s} \bigl(1 + \bigabs{\tau \mp
\abs{\xi}} \bigr)^{-\theta} \Check f(\xi). $$
By H\"older's inequality, \eqref{KlMaEmbeddingA} is bounded by
$\bigmixednorm{u}{\infty}{p}$ times $\mixednorm{v_{\pm}}{1}{r}$, where $1
= \frac{1}{r} + \frac{1}{p}$, so it suffices to show that
$$
\mixednorm{v_{\pm}}{1}{r} \lesssim \twonorm{f}{}. $$
But $v_{\pm}(t,\cdot) = c(t) D^{-s} e^{\pm itD} f$, where $\widehat
c(\tau) = (1+\abs{\tau})^{-\theta}$, and since $c \in L^{2}(\R)$, it
follows that
$$
\mixednorm{v_{\pm}}{1}{r} \lesssim \mixednorm{D^{-s} e^{\pm itD} f}{2}{r}.
$$
By Theorem A (see also Remark \ref{HalfWaveRemark}), the right hand side
is dominated by $\twonorm{f}{}$.
\end{proof}
The dual statement is as follows:
\begin{proposition}\label{DualKlMaEmbedding} If $\frac{2(n-1)}{n-3} \le r
< \infty$, $s = \frac{n}{2} - \frac{n}{r} - \half$ and $\theta > \half$,
then
$$
\Mixednorm{u}{1}{r}
\lesssim \Mixednorm{\Lambda^{s} \Lambda_{-}^{\theta} u}{1}{2} $$
for all $u$.
\end{proposition}
\begin{proof}
If $1 = 1/p + 1/r$, the hypotheses of Proposition \ref{KlMaEmbedding} are
satisfied, so the result follows by Proposition
\ref{DualityProposition}(b).
\end{proof}
\section{The Iteration Space}\label{TheIterationSpace}
The main point we want to make here is that proving local well-posedness
for a system $\square u = \mathcal N(u)$ with initial data in $H^s \times
H^{s-1}$ by iteration in some functional Banach space $\X^s$ (which should
satisfy certain conditions), reduces to proving estimates of the type
\begin{align}
\label{GenericEstimate}
\norm{\Lambda_+^{-1} \Lambda_-^{-1} \mathcal N (u)}_{\X^s} &\le A\bigl(
\norm{u}_{\X^s} \bigr),
\\
\label{GenericDifferenceEstimate}
\norm{\Lambda_+^{-1} \Lambda_-^{-1} \bigl( \mathcal N (u) - \mathcal N (v)
\bigr)}_{\X^s}
&\le A'\bigl( \max\{\norm{u}_{\X^s},\norm{v}_{\X^s} \} \bigr)
\norm{u-v}_{\X^s}, \end{align}
where $A$ and $A'$ are continuous functions and $A(0) = 0$.

Let us briefly describe the conditions that $\X^{s}$ should satisfy.

Firstly, we require that $\X^{s}$ embed in the continuation of the data
space $H^{s} \times H^{s-1}$, namely
\begin{equation}\label{GlobalEnergySpace}
C_b(\R,H^s) \cap C^1_b(\R,H^{s-1})
\end{equation}
with norm $\sup_{t \in \R} \bigl( \Sobnorm{u(t)}{s} +
\Sobnorm{\partial_{t} u(t)}{s-1} \bigr)$. In particular, this ensures that
restriction to any time-slab $S_{T} = [0,T] \times \R^{n}$ is well-defined
for elements of $\X^{s}$, and we denote by $\X^{s}_{T}$ the corresponding
restriction space.

Secondly, we limit our attention to spaces in which we have estimates for
the Cauchy problem for the \emph{linear} wave equation with data in $H^{s}
\times H^{s-1}$. Thus, we require that the solution of $\square u = F$
with initial data $(u, \partial_{t} u) \init = (f,g)$
satisfy\footnote{Here $\norm{u}_{\X^{s}_{T}}$ denotes the norm on the
restriction space $\X^{s}_{T}$. For example, if $\X^{s}$ is the data
continuation space \eqref{GlobalEnergySpace}, it is obvious which norm to
use on $\X^{s}_{T}$. It is less obvious if $\X^{s}$ is the space
$\scrH^{s,\theta}$, since the norm is then nonlocal in time. However, there is an
abstract way of defining a norm on $\X^{s}_{T}$ such that it becomes a
Banach space; see the statement of Theorem \ref{WellPosednessTheorem}.}
\begin{equation}\label{PreliminaryLinearEstimate}
\norm{u}_{\X^s_T} \lesssim \Sobnorm{f}{s} + \Sobnorm{g}{s-1} +
\norm{\Lambda_{+}^{-1} \Lambda_{-}^{-1} F}_{\X^{s}} \end{equation}
for all $0 < T < 1$, say. Let us write
$$
u = u_{0} + \square^{-1} F,
$$
where $u_{0}$ is the homogeneous part of the solution, i.e., $\square u_0
= 0$ with initial data $(f,g)$, and where $\square^{-1}$ is the operator
which to any sufficiently regular $F$ assigns the solution $v$ of $\square
v = F$ with $(v,\partial_t v) \init = 0$.

Thus, \eqref{PreliminaryLinearEstimate} splits into two estimates:
\begin{align*}
\norm{u_{0}}_{X^{s}_{T}} &\lesssim \Sobnorm{f}{s} + \Sobnorm{g}{s-1}, \\
\norm{\square^{-1} F}_{\X^s_T} &\lesssim \norm{\Lambda_{+}^{-1}
\Lambda_{-}^{-1} F}_{\X^{s}}. \end{align*}
The latter says that for the purpose of local-in-time estimates, we may
replace $\square^{-1}$ by the much nicer operator $\Lambda_{+}^{-1}
\Lambda_{-}^{-1}$.

The essential point is then the following: If we have a space $\X^{s}$
with the above properties (i.e., $\X^{s}$ embeds in the data continuation
space and \eqref{PreliminaryLinearEstimate} holds), and if the estimates
\eqref{GenericEstimate} and \eqref{GenericDifferenceEstimate} for the
nonlinearity $\mathcal N$ are true, then the equation $\square u =
\mathcal N(u)$ is locally well-posed for initial data in $H^{s} \times
H^{s-1}$. A precise statement is given in Theorem
\ref{SpecializedWellPosednessTheorem} below. %
\begin{remark}\label{PositivityRemark} If the norm on $\X^{s}$ only
depends on the size of the Fourier transform and is compatible with the
relation $\preceq$ (this terminology is defined in section
\ref{Notation}), and if there are operators $\mathcal N_{1}, \dots,
\mathcal N_{m}$ such that $$
u \preceq v \implies \mathcal N(u) \preceq \mathcal N_{1}(v) + \cdots +
\mathcal N_{m}(v),
$$
then in order to prove \eqref{GenericEstimate}, it suffices to prove $$
\norm{\Lambda_+^{-1} \Lambda_-^{-1} \mathcal N_j(v)}_{\X^s} \lesssim
A\bigl(\norm{v}_{\X^s}\bigr)
\qquad (1 \le j \le m)
$$
for all $v$ with $\widehat v \ge 0$.
\end{remark}
We start by discussing in fairly general terms how estimates imply local
well-posedness.
\subsection{Well-Posedness: A General Point of View} %
Consider again the generic Cauchy problem
\begin{equation}\label{GenericSystem}
\square u = \mathcal N(u), \qquad
(u,\partial_t u) \init = (f,g).
\end{equation}
Assume that $\mathcal N(0) = 0$ and $\mathcal N$ is local in time.

Associated to the Cauchy problem \eqref{GenericSystem} is the sequence
$(u_j)$ of iterates, defined inductively by setting $u_{-1} \equiv 0$ and
$$
\square u_{j} = \mathcal N(u_{j-1}), \qquad (u,\partial_t u) \init = (f,g)
$$
for $j \ge 0$. Thus $u_{0}$ is the homogeneous part of the solution, and
the subsequent iterates are given by \begin{equation}\label{InductiveStep}
u_{j+1} = u_0 + \square^{-1} \mathcal N(u_j) \end{equation}
for $j \ge 0$.

The strategy for proving local existence for \eqref{GenericSystem} in a
time slab $S_T = [0,T] \times \R^n$, for some $T > 0$ and for given data
$(f,g) \in H^s \times H^{s-1}$, is to find a Banach space
\begin{equation}\label{LocalEnergyEmbedding}
\X^s_T \hookrightarrow C([0,T],H^s) \cap C^1([0,T],H^{s-1}) \end{equation}
in which $\bigl( u_j \vert S_T \bigr)$ is Cauchy. The limit $u$ will then
be a solution of \eqref{GenericSystem} on $S_T$, provided that $\mathcal
N(u_j) \to \mathcal N(u)$ in the sense of distributions on $(0,T) \times
\R^n$ (this always follows from the estimates involving $\mathcal N$).

To prove that $(u_j)$ is Cauchy, we need estimates. Firstly, the inductive
step \eqref{InductiveStep} must be well-defined, so we need (for $0 < T <
1$, say) \begin{align}
\label{InductiveStepA}
\norm{u_0}_{\X^s_T} &\le C \bigl( \Sobnorm{f}{s} + \Sobnorm{g}{s-1}
\bigr), \\
\label{InductiveStepB}
\norm{\square^{-1} \mathcal N(u)}_{\X^s_T} &\le C_T A\bigl(
\norm{u}_{\X^s_T} \bigr), \\
\intertext{where $A$ is a continuous function vanishing at $0$. We may
always assume that $A$ is increasing. Secondly, we need estimates for the
difference of two iterates; i.e., we need}
\label{LocalDifferenceEstimate}
\norm{\square^{-1} \bigl( \mathcal N(u) - \mathcal N(v) \bigr) }_{\X^s_T}
&\le C_T A'\bigl( \max\{\norm{u}_{\X^s_T},\norm{v}_{\X^s_T} \} \bigr)
\norm{u-v}_{\X^s_T}, \end{align}
where $A'$ is continuous.
\begin{example}\label {ClassicalLocalExistenceExample} To prove the
Classical Local Existence Theorem (section \ref{Motivation}) for a system
$\square u = F(u,\partial u)$ , where $F$ is smooth and vanishes at the
origin, we set $$
\X^s_T = C([0,T],H^s) \cap C^1([0,T],H^{s-1}). $$
Then \eqref{InductiveStepA}, \eqref{InductiveStepB} and
\eqref{LocalDifferenceEstimate} hold for any $s > \frac{n}{2} + 1$, with
$C_T = O(T)$ as $T \to 0$. Indeed, by the energy inequality, \begin{align}
\notag
\norm{u_0}_{\X^s_T} &\le C \bigl( \Sobnorm{f}{s} + \Sobnorm{g}{s-1} \bigr)
\\
\intertext{and}
\label{ClassicalExampleA}
\norm{\square^{-1} F(u,\partial u)}_{\X^s_T} &\le C \int_{0}^{T}
\Sobnorm{F\bigl( u(t),\partial u(t) \bigr)}{s-1} \, dt
\end{align}
for all $0 \le T \le 1$, say.
Recall the Moser inequality, which says that if $\Gamma$ is smooth and
vanishes at the origin, and if $\sigma \ge 0$, then there exists a
continuous function $g : [0,\infty) \to [0,\infty)$ such that
\begin{equation}\label{MoserInequality}
\Sobnorm{\Gamma(f)}{\sigma} \le g( \inftynorm{f}{} ) \Sobnorm{f}{\sigma}
\end{equation}
for all $f \in H^\sigma \cap L^\infty$ ($f$ may be $\R^{N}$-valued). See,
e.g., Meyer \cite{Meyer}. By applying this, we get
\begin{multline}\label{ClassicalExampleB}
\Sobnorm{F\bigl( u(t),\partial u(t) \bigr)}{s-1} \\
\le g\bigl( \inftynorm{u(t)}{} + \inftynorm{\partial u(t)}{} \bigr) \bigl(
\Sobnorm{u(t)}{s} + \Sobnorm{\partial_{t} u(t)}{s-1} \bigr).
\end{multline}
Since $s-1 > \frac{n}{2}$, the $L^{\infty}$ Sobolev embedding
\eqref{InftySobolevEmbedding} gives
\begin{equation}\label{ClassicalExampleC}
\inftynorm{u(t)}{} + \inftynorm{\partial u(t)}{} \lesssim
\Sobnorm{u(t)}{s} + \Sobnorm{\partial_{t} u(t)}{s-1}. \end{equation}
By combining \eqref{ClassicalExampleA},\eqref{ClassicalExampleB} and
\eqref{ClassicalExampleC}, we obtain \eqref{InductiveStepB}. The proof of
the difference estimate \eqref{LocalDifferenceEstimate} is similar.
\end{example}
If \eqref{InductiveStepA} and \eqref{InductiveStepB} hold, and we let $R$
be twice the right hand side of \eqref{InductiveStepA}, then it follows by
induction that $\norm{u_j}_{\X^s_T} \le R$ provided $2 C_T A(R) \le R$
(keep in mind that $A$ is increasing).
There are two ways to ensure that the latter inequality holds: Take $T$
small (if $\lim_{T \to 0^+} C_T = 0$), or (since $A(0) = 0$) require that
$R$ be small, i.e., the data have small norm.

Then, by the difference estimate \eqref{LocalDifferenceEstimate}, we have
$$\norm{u_{j+1} - u_j}_{\X^s_T}
\le \half \norm{u_{j} - u_{j-1}}_{\X^s_T}$$ provided $2 C_T A'(R) \le R$
(so we need either $\lim_{T \to 0^+} C_T = 0$ or $A'(0) = 0$).

It follows that $(u_j)$ is Cauchy, establishing local existence. With a
little more work one can then prove uniqueness of local solutions in
$\X^s_T$ for any $T > 0$, and local Lipschitz continuity. For a fuller
discussion we refer to Selberg \cite{Se, Se2},
where the following is proved.
\begin{theorem}\label{WellPosednessTheorem} Let $\X^s$ be a Banach space
which embeds in \eqref{GlobalEnergySpace} and is time-translation
invariant:
$$
\norm{u(\cdot + T, \cdot)}_{\X^s} = \norm{u}_{\X^s} \quad (\forall T). $$
Also assume that for all $\phi \in C_{c}^{\infty}(\R)$, the multiplication
map $u \mapsto \phi(t) u(t,x)$ is bounded from $\X^s$ into itself.

For any $T > 0$, let $\X^s_T$ be the restriction of $\X^s$ to $[0,T]
\times \R^n$. (That is, we define an equivalence relation $\sim_T$ on
$\X^s$ by $$
u \sim_{T} v \iff
u(t) = v(t) \quad \forall 0 \le t \le T. $$
Since $\X^{s}$ embeds in \eqref{GlobalEnergySpace}, the equivalence
classes are closed sets in $\X^s$, so the quotient $\X^s_{T} = \X^s /
\!\!\sim_{T}$, with norm $$
\norm{u}_{\X^s_{T}} = \inf_{v \sim_{T} u} \norm{v}_{\X^s}, $$
is a Banach space.)

Consider the system \eqref{GenericSystem}. Assume that the estimates
\eqref{InductiveStepA},
\eqref{InductiveStepB} and \eqref{LocalDifferenceEstimate} hold for all $0
< T < 1$, $(f,g) \in H^s \times H^{s-1}$ and $u,v \in \X^s_T$. Moreover,
assume that
$$
\lim_{T \to 0^+} C_T = 0.
$$
Then \eqref{GenericSystem} is locally well-posed for initial data in $H^s
\times H^{s-1}$, in the following precise sense:
{ \renewcommand{\theenumi}{\alph{enumi}}
\renewcommand{\labelenumi}{(\theenumi)}
\begin{enumerate}
\item {\bfseries(Existence)}
For all $(f,g) \in H^{s} \times H^{s-1}$ there is a time $$
0 < T = T_{s}(\Sobnorm{f}{s} + \Sobnorm{g}{s-1}) $$
which depends continuously on the norm of the data, and there is a $u \in
\X^s_{T}$ which solves \eqref{GenericSystem} on $S_{T} = [0,T] \times
\R^{n}$.
\item {\bfseries(Uniqueness)} If $T > 0$ and $u,u' \in \X^s_T$ are two
solutions of \eqref{GenericSystem} on $S_{T}$ with the same data $(f,g)$,
then $u = u'$ in $\X^s_T$. %
\item {\bfseries(Continuous dependence on initial data)} If $u \in
\X^{s}_{T}$ solves \eqref{GenericSystem} on $S_{T}$ for some $T > 0$, then
for all $(f',g') \in H^{s} \times H^{s-1}$ sufficiently close to $(f,g)$
there exists a $u' \in \X^{s}_{T}$ which solves \eqref{GenericSystem} on
$S_{T}$ with data $(f',g')$, and
$$
\norm{u - u'}_{\X^{s}_{T}} \lesssim \Sobnorm{f-f'}{s} +
\Sobnorm{g-g'}{s-1}.
$$
\end{enumerate}
}
\end{theorem}
\subsection{Specialization to the Relevant Spaces} %
Here we specialize the preceding discussion to suit our present needs. The
following theorem gives a precise form to the idea outlined at the
beginning of this section. %
\begin{theorem}\label{SpecializedWellPosednessTheorem} Let $\X^s$ be a
Banach space with the following properties: {
\renewcommand{\theenumi}{\alph{enumi}}
\renewcommand{\labelenumi}{(\theenumi)}
\begin{enumerate}
\item $\X^{s}$ embeds in \eqref{GlobalEnergySpace}. %
\item The norm on $\X^{s}$ is invariant under time-translation. %
\item The estimate
$$
\norm{u_0}_{\X^s_T} \le C \bigl( \Sobnorm{f}{s} + \Sobnorm{g}{s-1} \bigr)
$$
holds for all $(f,g) \in H^{s} \times H^{s-1}$ and $0 < T < 1$, where $u_{0}$ is the
solution of the homogeneous wave equation with initial data $(f,g)$ and
$\X^s_T$ is the restriction space.
\item
For the purposes of local-in-time estimates, $\square^{-1}$ may be
replaced with $\Lambda_{+}^{-1} \Lambda_{-}^{-1}$. More precisely, assume
that for $\varepsilon > 0$ sufficiently small,
\begin{equation}\label{LocalInhomogeneousEstimate}
\norm{\square^{-1} F}_{\X^s_T}
\le C_{T,\varepsilon} \norm{\Lambda_+^{-1} \Lambda_-^{\varepsilon-1}
F}_{\X^s} \end{equation}
for all $F \in \Lambda_+ \Lambda_-^{1-\varepsilon} \X^s$ and $0 < T < 1$.
Furthermore, suppose
\begin{equation}\label{ContinuityAtZero}
\lim_{T \to 0^{+}} C_{T,\varepsilon} = 0 \end{equation}
for $\varepsilon > 0$.
\item For all $\phi \in C_{c}^{\infty}(\R)$, the multiplication map $u
\mapsto \phi(t) u(t,x)$ is bounded from $\X^s$ into itself. \end{enumerate}
}
Consider the Cauchy problem \eqref{GenericSystem}. Suppose
\begin{align}
\label{EpsilonGenericEstimate}
\norm{\Lambda_+^{-1} \Lambda_-^{\varepsilon-1} \mathcal N (u)}_{\X^s} &\le
A\bigl( \norm{u}_{\X^s} \bigr),
\\
\label{EpsilonGenericDifferenceEstimate} \norm{\Lambda_+^{-1}
\Lambda_-^{\varepsilon-1} \bigl( \mathcal N (u) - \mathcal N (v)
\bigr)}_{\X^s}
&\le A'\Bigl( \frac{1}{3}
\max\{\norm{u}_{\X^s},\norm{v}_{\X^s} \} \Bigr) \norm{u-v}_{\X^s}
\end{align}
for all $u,v \in \X^{s}$, where $A$ and $A'$ are continuous and $A(0) =
0$.

Then \eqref{GenericSystem} is locally well-posed for initial data in
$H^{s} \times H^{s-1}$, with uniqueness of solutions in $\X^{s}_{T}$ for
any $T > 0$. \end{theorem}
\begin{remarkNoLabel} In our applications of this theorem, $\varepsilon$
must be strictly positive to ensure that \eqref{ContinuityAtZero} holds.
The latter is not needed, however, if one imposes instead a smallness
assumption on the norms of the initial data. More precisely, if we take
$\varepsilon = 0$ in the above theorem, and if we replace
\eqref{ContinuityAtZero} with the assumption that
\begin{equation}\label{WeakerContinuityAtZero}
\limsup_{t \to 0^{+}}
\norm{u}_{\X^{s}_{t}} \le C( \Sobnorm{u(0)}{s} + \Sobnorm{\partial_{t}
u(0)}{s-1} )
\end{equation}
for all $u \in \X^{s}$, then the conclusion of the theorem still holds,
but we must require that the initial data satisfy $\Sobnorm{f}{s} +
\Sobnorm{g}{s-1} < \delta$ for some sufficiently small $\delta > 0$.
\end{remarkNoLabel}
\begin{proof}[Proof of Theorem \ref{SpecializedWellPosednessTheorem}] By
Theorem \ref{WellPosednessTheorem}, it suffices to prove the estimates
\eqref{InductiveStepB} and \eqref{LocalDifferenceEstimate} for $0 < T < 1$.

Fix $T$ and $u \in \X^s_T$. Let $u' \in \X^s$ be any extension of $u$
(meaning $u' \sim_T u$). Since $\mathcal N$ is local in time, we have
$\square^{-1} \mathcal N(u) = \square^{-1} \mathcal N(u')$ on $[0,T]
\times \R^{n}$, so by
\eqref{LocalInhomogeneousEstimate} and \eqref{EpsilonGenericEstimate}, $$
\norm{\square^{-1} \mathcal N(u)}_{\X^s_T} \le C_{T,\varepsilon} A\bigl(
\norm{u'}_{\X^s} \bigr). $$
Let $\norm{u'}_{\X^s} \to \norm{u}_{\X^s_T}$. Since $A$ is continuous,
\eqref{InductiveStepB} follows.

To prove \eqref{LocalDifferenceEstimate}, fix $u,v \in \X^s_T$, and let
$u',v' \in \X^s$ be any two extensions of $u$ and $v$. By
\eqref{LocalInhomogeneousEstimate} and
\eqref{EpsilonGenericDifferenceEstimate}, $$
\norm{\square^{-1} \bigl( \mathcal N(u) - \mathcal N(v) \bigr) }_{\X^s_T}
\le C_{T,\varepsilon} A'\Bigl( \frac{1}{3}
\max\{\norm{u'}_{\X^s},\norm{v'}_{\X^s} \} \Bigr) \norm{u'-v'}_{\X^s}. $$
Let $w = u' - v'$, and write
$$
\max\{\norm{u'}_{\X^s},\norm{v'}_{\X^s} \} \le
\max\{\norm{u'}_{\X^s},\norm{w}_{\X^s} + \norm{u'}_{\X^s} \} =
\norm{w}_{\X^s} + \norm{u'}_{\X^s}.
$$
We may assume that $A'$ is increasing. Thus $$
\norm{\square^{-1} \bigl( \mathcal N(u) - \mathcal N(v) \bigr) }_{\X^s_T}
\le C_{T,\varepsilon} A'\Bigl( \frac{1}{3} \{ \norm{w}_{\X^s} +
\norm{u'}_{\X^s} \} \Bigr) \norm{w}_{\X^s}. $$
Let $\norm{u'}_{\X^s} \to \norm{u}_{\X^s_T}$ and $\norm{w}_{\X^s} \to
\norm{u-v}_{\X^s_T}$. Since
$$
\norm{u-v}_{\X^s_T} + \norm{u}_{\X^s_T} \le 2 \norm{u}_{\X^s_T} +
\norm{v}_{\X^s_T} \le 3 \max\{\norm{u}_{\X^s_T},\norm{v}_{\X^s_T} \}, $$
we conclude that \eqref{LocalDifferenceEstimate} holds. \end{proof}
The next theorem gives sufficient conditions for $\X^s$ to satisfy
properties (c) and (d) of Theorem \ref{SpecializedWellPosednessTheorem}. %
\begin{theorem}\label{BasicConditionsTheorem} Let $\X^{s}$ be a Banach
space satisfying:
{ \renewcommand{\theenumi}{\alph{enumi}}
\renewcommand{\labelenumi}{(\theenumi)}
\begin{enumerate}
\item $\X^s \hookrightarrow \scrH^{s,\theta}$ for some $\half < \theta <
1$;
\item $\norm{u}_{\X^s} \le \norm{v}_{\X^s}$ whenever $u \preceq v$; %
\item There exists $\alpha \le \frac{5}{4} + \frac{\theta}{2}$ such that $$
\norm{u}_{\X^{s}} \lesssim \bignorm{ \Fourier \Lambda^{s-1} \Lambda_{+}
\Lambda_{-}^{\alpha} u (\tau,\xi) }_{L_{\xi}^{2}(L_{\tau}^{\infty})}
$$
for all $u \in \X^{s}$.
\end{enumerate}
}
Let $\X^s_T$ be the restriction, defined as in Theorem
\ref{WellPosednessTheorem}. Fix $0 \le \varepsilon < 1 - \theta$. Then the
solution of the linear Cauchy problem $$
\square u = F, \qquad (u,\partial_t u) \init = (f,g) $$
satisfies
$$
\norm{u}_{\X^s_T} \lesssim \Sobnorm{f}{s} + \Sobnorm{g}{s-1} +
T^{\varepsilon/2} \norm{\Lambda_{+}^{-1} \Lambda_{-}^{\varepsilon-1}
F}_{\X^{s}} $$
for all $0 < T \le 1$, $(f,g) \in H^s \times H^{s-1}$ and $F \in
\Lambda_{+} \Lambda_{-}^{1-\varepsilon} \X^{s}$. \end{theorem}
The proof can be found in \cite{Se2}.

Next, we verify that the iteration spaces used in the proof of the Main
Theorem satisfy the hypotheses of the previous theorem. The spaces in
parts (a)--(c) below are the iteration spaces used to prove parts
(a)--(c), respectively, of the Main Theorem. %
\begin{proposition}\label{BasicSpacesProposition} The hypotheses of
Theorem \ref{BasicConditionsTheorem} are satisfied by the following spaces:
{ \renewcommand{\theenumi}{\alph{enumi}}
\renewcommand{\labelenumi}{(\theenumi)}
\begin{enumerate}
\item $\X^{s} = \scrH^{s,\theta}$, provided $\half < \theta <
\frac{3}{2}$. %
\item The space $\X^s$ given by the norm $$
\norm{u} = \Spacetimenorm{u}{s}{\theta} + \bigMixednorm{\Lambda^{\gamma}
\Lambda_-^{\half} u}{1}{2n},
$$
where $n \ge 4$, $s > \frac{n-2}{2}$, $\half < \theta < \frac{3}{2}$ and
$0 < \gamma \le s - \frac{n-2}{2}$.
\item The space $\X^s$ given by the norm $$
\norm{u} = \Spacetimenorm{u}{s}{\theta} + \Mixednorm{\Lambda^{-1}
\Lambda_- u}{q}{\infty},
$$
where $n \ge 3$, $s > \frac{n-2}{2}$, $\half < \theta < \frac{3}{2}$ and
$1 \le q \le 2$.
\end{enumerate}
}
\end{proposition}
\begin{proof} For $\scrH^{s,\theta}$, we only have to note that $$
\Spacetimenorm{u}{s}{\theta} = \twonorm{ \Fourier \Lambda^{s-1}
\Lambda_{+} \Lambda_{-}^{\theta} u(\tau,\xi) }{_{\tau,\xi}} \le C_\delta
\bignorm{ \Fourier \Lambda^{s-1} \Lambda_{+} \Lambda_{-}^{\theta + \half +
\delta} u (\tau,\xi) }_{L_{\xi}^{2}(L_{\tau}^{\infty})}, $$
where
$$
C_\delta^2 \simeq \int_\R (1 + \abs{\lambda}^2)^{-\half-\delta} \,
d\lambda < \infty $$
for any $\delta > 0$. This proves part (a).

Let $\X^s$ be the space in part (b).
Since $n \ge 4$, Proposition \ref{DualKlMaEmbedding} gives $$
\bigMixednorm{ \Lambda^\gamma \Lambda_-^{\half} u}{1}{2n} \lesssim
\bigMixednorm{ \Lambda^{\gamma + \frac{n-2}{2}} \Lambda_-^{1 + \delta}
u}{1}{2}
$$
for any $\delta > 0$. We claim that
\begin{equation}\label{BasicSpacesClaim}
\Mixednorm{u}{1}{2} \lesssim
\norm{ \widehat u(\tau,\xi)}_{L_{\xi}^{2}(L_{\tau}^{\infty})}
\end{equation}
for all $u$; this finishes the proof of part (b), if we choose $\delta$
sufficiently small.

To prove the claim, notice that
$$
\Mixednorm{u}{1}{2} = \sup_{v} \int \abs{\widehat u(X)} \widehat v(X) \,
dX \le \norm{ \widehat u(\tau,\xi)}_{L_{\xi}^{2}(L_{\tau}^{\infty})}
\sup_{v} \norm{ \widehat v(\tau,\xi)}_{L_{\xi}^{2}(L_{\tau}^{1})}, $$
where the supremum is over all $v \in \Schwartz$ such that $\widehat v \ge
0$ and $\mixednorm{v}{\infty}{2} = 1$. But for such $v$, we have $\norm{
\widehat v(\tau,\xi)}_{L_{\xi}^{2}(L_{\tau}^{1})} =
(2\pi)^{-1}\twonorm{v(0,\cdot)}{(\R^{n})} \le (2\pi)^{-1}$ by Fourier
inversion.

Finally, we prove part (c). Since the norm only depends on the size of the
Fourier transform and is compatible with the relation $\preceq$, it
suffices, by
\eqref{ModifiedMixednormBound}, to prove that $$
\mixednorm{\Lambda^{-1} \Lambda_{-} u}{q}{\infty} \lesssim \bignorm{
\Fourier \Lambda^{s}
\Lambda_{-}^{\alpha} u (\tau,\xi) }_{L_{\xi}^{2}(L_{\tau}^{\infty})} $$
for some $\alpha \le \frac{5}{4} + \frac{\theta}{2}$. Let $\frac{1}{q'} =
1 - \frac{1}{q}$ and choose $\beta$ satisfying $\frac{1}{q'} < \beta \le
\frac{1}{4} + \frac{\theta}{2}$. By the inequalities of Hausdorff-Young,
Minkowski and H\"older, \begin{align*}
\mixednorm{v}{q}{\infty}
&\le \bignorm{ \widehat{v(t)}(\xi) }_{L_{t}^{q}(L_{\xi}^{1})} \le
\bignorm{ \widehat{v(t)}(\xi) }_{L_{\xi}^{1}(L_{t}^{q})} \le \bignorm{
\widehat{v}(\tau,\xi) }_{L_{\xi}^{1}(L_{\tau}^{q'})} \\
&\lesssim \bignorm{ (1 + \abs{\xi})^{-\frac{n}{2}-\varepsilon} \bigl(1 +
\hypwt{\tau}{\xi}\bigr)^{-\beta} \Fourier \Lambda^{\frac{n}{2} +
\varepsilon} \Lambda_{-}^{\beta} v (\tau,\xi)
}_{L_{\xi}^{1}(L_{\tau}^{q'})} \\
&\lesssim \bignorm{ \Fourier \Lambda^{\frac{n}{2} + \varepsilon}
\Lambda_{-}^{\beta} v (\tau,\xi) }_{L_{\xi}^{2}(L_{\tau}^{\infty})}
\end{align*}
for any $\varepsilon > 0$. This finishes the proof. \end{proof}
Finally, we check that condition (e) of Theorem
\ref{SpecializedWellPosednessTheorem} holds for the spaces that we use. %
\begin{proposition}\label{TimeCutOffProposition} Let
$\gamma,\gamma_+,\gamma_{-} \in \R$ and $1 \le q,r \le \infty$.
Let $\phi \in \Schwartz(\R)$, $M_{\phi}u(t,x)=\phi(t)u(t,x)$. Then $$
\Mixednorm{\Lambda^{\gamma} \Lambda_{+}^{\gamma_{+}}
\Lambda_{-}^{\gamma_{-}} M_{\phi}
u}{q}{r}
\le C_{\gamma_+,\gamma_{-},\phi}
\Mixednorm{\Lambda^{\gamma} \Lambda_{+}^{\gamma_{+}}
\Lambda_{-}^{\gamma_{-}} u}{q}{r}.
$$
\end{proposition}
\begin{proof}
Since $\Lambda^{\gamma} M_{\phi} u = M_{\phi} (\Lambda^{\gamma} u)$, we
may assume $\gamma = 0$.

Let $\chi$ and $v$ be defined by $\widehat \chi (\tau) = (1 +
\abs{\tau}^{2})^\frac{\gamma_{+} + \gamma_{-}}{2} \bigabs{\widehat \phi}$
and $\widehat v = \abs{\widehat u}$. Since \begin{gather*}
\hypwt{\tau}{\xi}
\le \hypwt{\tau}{\lambda} + \hypwt{\lambda}{\xi} \le \abs{\tau-\lambda} +
\hypwt{\lambda}{\xi}, \\
\abs{\tau} + \abs{\xi}
\le \abs{\tau-\lambda} + \abs{\lambda} + \abs{\xi}, \end{gather*}
we have
$$
\Lambda_{+}^{\gamma_{+}} \Lambda_{-}^{\gamma_{-}} M_{\phi} u \precsim
M_{\chi} (\Lambda_{+}^{\gamma_{+}} \Lambda_{-}^{\gamma_{-}} v). $$
Thus, it suffices to prove
$$
\Mixednorm{M_{\chi} v}{q}{r} \lesssim \Mixednorm{v}{q}{r} $$
for all $v$ such that $\widehat v \ge 0$.

Since $\chi \in H^{\alpha}(\R)$ for all $\alpha > 0$, we have $\chi \in
C^{\infty}$ and $\chi^{(j)} \in L^{\infty}$ for all $k \ge 0$, so
$M_{\chi}$ maps the set
$\{ w \in \Schwartz(\R^{1+n}) : \widehat w \ge 0 \}$ into itself. Thus, if
$\widehat v \ge 0$,
\begin{multline*}
\Mixednorm{M_{\chi} v}{q}{r}
= \sup_{w} \int (M_{\chi} v) w \, dt \, dx \\
= \sup_{w} \int v M_{\chi} w \, dt \, dx \le \Mixednorm{v}{q}{r} \sup_{w}
\mixednorm{M_{\chi} w}{q'}{r'}, \end{multline*}
where $q',r'$ denote the dual exponents of $q,r$ and the supremum is over
all $w \in \Schwartz$ such that $\widehat w \ge 0$ and
$\mixednorm{w}{q'}{r'} = 1$. But by H\"older's inequality and Sobolev
embedding,
$$
\mixednorm{M_{\chi} w}{q'}{r'}
\lesssim
\inftynorm{\chi}{(\R)} \mixednorm{w}{q'}{r'} \lesssim
\norm{\chi}_{H^{\half + \varepsilon}(\R)} \mixednorm{w}{q'}{r'}
$$
for any $\varepsilon > 0$, and we have
$\norm{\chi}_{H^{\half + \varepsilon}(\R)} = \norm{\phi}_{H^{\gamma_+ +
\gamma_{-} + \half + \varepsilon}(\R)}$. \end{proof}
\section{Some Special Embeddings}
We collect here some embeddings that are used repeatedly in the proof of
the Main Theorem.
\begin{alignat}{2}
\label{EnergyEmbedding}
H^{0,\theta} &\hookrightarrow \mixed{\infty}{2}, & \qquad & \theta > 1/2,
\\
\label{InftyEmbedding}
H^{s,\theta} &\hookrightarrow L^\infty(\R^{1+n}), && \text{$s > n/2$,
$\theta > 1/2$},
\\
\label{SpaceInftyEmbedding}
H^{s,0} &\hookrightarrow \mixed{2}{\infty}, && s > n/2,
\\
\label{SpecialEmbeddingA}
H^{1,0} &\hookrightarrow \mixed{2}{\frac{2n}{n-2}}, && n \ge 3,
\\
\label{SpecialEmbeddingB}
H^{1,\theta} &\hookrightarrow \mixed{2}{\frac{2n}{n-3}}, && \text{$n \ge
4$, $\theta > \half$},
\\
\label{SpecialEmbeddingC}
H^{\frac{n-1}{2},0} &\hookrightarrow \mixed{2}{2n}, &&
\\
\label{SpecialEmbeddingD}
H^{\frac{n-1}{2}+\varepsilon,\theta} &\hookrightarrow \mixed{2}{\infty},
&& \text{$n \ge 4$, $\varepsilon > 0$, $\theta > \half$}, \\
\label{SpecialEmbeddingE}
H^{\frac{n-2}{2},\theta} &\hookrightarrow \mixed{2}{2n}, && \text{$n \ge
4$, $\theta > \half$},
\\
\label{SpecialEmbeddingF}
H^{\frac{n-3}{2},0} &\hookrightarrow \mixed{2}{\frac{2n}{3}}, && n \ge 3
\\
\label{SpecialEmbeddingG}
H^{\frac{n-3}{2},\theta} &\hookrightarrow \mixed{2}{n}, && \text{$n \ge
5$, $\theta > \half$},
\\
\label{SpecialEmbeddingH}
H^{\frac{n-3}{2},\theta} &\hookrightarrow \mixed{\infty}{\frac{2n}{3}}, &&
\text{$n \ge 3$, $\theta > \half$},
\\
\label{SpecialEmbeddingI}
\Lambda^{-\half-\varepsilon} \mixed{1}{2n} &\hookrightarrow
\mixed{1}{\infty}, && \varepsilon > 0,
\\
\label{SpecialEmbeddingJ}
\Lambda^{-\half-\varepsilon} \Mixed{1}{2n} &\hookrightarrow
\Mixed{1}{\infty}, && \varepsilon > 0,
\\
\label{SpecialEmbeddingK}
\Mixed{q}{\infty} \cdot H^{0,\theta} &\hookrightarrow H^{0,-\alpha}, &&
\text{$1 \le q \le 2$, $\alpha > \frac{1}{q} - \half$, $\theta > \half$},
\\
\label{SpecialEmbeddingL}
\Lambda^{-\half-\varepsilon} \Mixed{1}{2n} \cdot H^{0,\theta}
&\hookrightarrow H^{0,-\theta},
&& \text{ $\varepsilon > 0$, $\theta > \half$}. \end{alignat}
\begin{remark}\label{DualityRemark} In view of Corollary
\ref{DualityCorollary}, the estimates
\eqref{EnergyEmbedding}--\eqref{SpecialEmbeddingH} remain valid if we
replace the $\mixed{q}{r}$-space on the right by the corresponding
$\Mixed{q}{r}$-space. \end{remark}
\paragraph{Proofs:}
\begin{enumerate}
\item \eqref{EnergyEmbedding} follows from Proposition
\ref{EnergyEmbeddingProposition}; then, since $$
\Lambda^{-s} L^{2}(\R^{n}) \hookrightarrow L^\infty(\R^{n}) $$
for $s > \frac{n}{2}$ by Sobolev embedding, and since $H^{s,\theta} =
\Lambda^{-s} H^{0,\theta}$ and $H^{s,0} = \Lambda^{-s} L^{2}$,
\eqref{InftyEmbedding} and \eqref{SpaceInftyEmbedding} follow. %
\item \eqref{SpecialEmbeddingB}, \eqref{SpecialEmbeddingE} and
\eqref{SpecialEmbeddingG} are special cases of Theorem D. %
\item \eqref{SpecialEmbeddingA}, \eqref{SpecialEmbeddingC} and
\eqref{SpecialEmbeddingF} hold by Sobolev embedding. %
\item \eqref{SpecialEmbeddingD} follows from \eqref{SpecialEmbeddingE} and
$$
\mixed{2}{2n} \hookrightarrow \Lambda^{\half + \varepsilon}
\mixed{2}{\infty}.
$$
The latter holds by Sobolev embedding.
\item \eqref{SpecialEmbeddingH} follows from \eqref{EnergyEmbedding} and
the Sobolev embedding $H^\frac{n-3}{2} \hookrightarrow
L^\frac{2n}{3}(\R^n)$, since these imply
$$
H^{\frac{n-3}{2},\theta} = \Lambda^{-\frac{n-3}{2}} H^{0,\theta}
\hookrightarrow \Lambda^{-\frac{n-3}{2}} \mixed{\infty}{2} \hookrightarrow
\mixed{\infty}{\frac{2n}{3}}. $$
\item \eqref{SpecialEmbeddingI} holds by Sobolev embedding;
\eqref{SpecialEmbeddingJ} follows by Corollary \ref{DualityCorollary}. %
\item \eqref{SpecialEmbeddingK} is proved as follows: By interpolation
between \eqref{EnergyEmbedding} and $L^{2} \hookrightarrow L^{2}$,
\begin{align*}
\mixed{q}{2} &\hookrightarrow H^{0,-\alpha} \\
\intertext{for $1 \le q \le 2$ and $\alpha > \frac{1}{q} - \half$, and it
follows that}
\Mixed{q}{2} &\hookrightarrow H^{0,-\alpha} \end{align*}
for such $q, \alpha$. Combining this with Proposition
\ref{MixednormHolder}, Lemma \ref{DensityLemma} and
\eqref{EnergyEmbedding}, we get \eqref{SpecialEmbeddingK}.
\item \eqref{SpecialEmbeddingL} follows from \eqref{SpecialEmbeddingK} and
\eqref{SpecialEmbeddingJ}. \end{enumerate}
\section{Main Estimates for \eqref{WMType}}\label{WMProof}
Our aim here is to prove part (a) of the Main Theorem: \eqref{WMType} is
locally well-posed for initial data in $H^{s} \times H^{s-1}$ for $s >
\frac{n}{2}$ and $n \ge 2$.

It suffices to verify the hypotheses of Theorem
\ref{SpecializedWellPosednessTheorem}. Take $\X^{s} = \scrH^{s,\theta}$,
where $\theta = \theta(s,n) > \half$ is to be determined. The nonlinearity
is
$$
\mathcal N^I (u) = - \sum_{J,K = 1}^N \Gamma^I_{JK}(u) Q_0(u^J,u^K),
\qquad 1 \le I \le N. $$
Let us first check that conditions (a)--(e) of Theorem
\ref{SpecializedWellPosednessTheorem} are satisfied. Condition (a) holds
by Proposition \ref{SecondEnergyEmbeddingProposition}; condition (b) is
obviously satisfied; conditions (c) and (d) follow from Theorem
\ref{BasicConditionsTheorem}, in view of Proposition
\ref{BasicSpacesProposition}; finally, condition (e) holds by Proposition
\ref{TimeCutOffProposition}.

It remains to prove \eqref{EpsilonGenericEstimate} and
\eqref{EpsilonGenericDifferenceEstimate}. Let us first prove
\eqref{EpsilonGenericEstimate} in the case where the $\Gamma^{I}_{JK}$'s
are constants, and let us set $\varepsilon = 0$. The proof of the general
case is quite similar; it appears in section \ref{WMDetailedProof}.
\subsection{The Simplified Case}
Since $\Lambda_{+} \Lambda_{-} \scrH^{s,\theta} = H^{s-1,\theta-1}$, what
we want to prove is the following:
\begin{theorem}\label{WMTheorem1}
Suppose $n \ge 2$, $s > \frac{n}{2}$ and $\half < \theta \le s -
\frac{n-1}{2}$. Then
$$
Q_{0}\bigl(\scrH^{s,\theta},\scrH^{s,\theta}\bigr) \hookrightarrow
H^{s-1,\theta-1}.
$$
\end{theorem}
By estimating the Fourier symbol of the null form $Q_{0}$ in absolute
value, it is easy to prove (see Lemma \ref{Q0Estimate} below) that
$$
Q_0 (\phi,\psi) \preceq D_+ D_- (\phi' \psi') + (D_+ D_- \phi') \psi' +
\phi' D_+ D_- \psi' $$
whenever $\phi \preceq \phi'$ and $\psi \preceq \psi'$. Therefore, in view
of Remark \ref{PositivityRemark}, it suffices to prove \begin{align*}
\scrH^{s,\theta} \cdot \scrH^{s,\theta} &\hookrightarrow \scrH^{s,\theta},
\\
H^{s-1,\theta-1} \cdot \scrH^{s,\theta} &\hookrightarrow H^{s-1,\theta-1}.
\end{align*}

To keep the discussion as simple as possible (the complete details appear
in section \ref{WMDetailedProof}), let us for the moment ignore the
difference between $\scrH^{s,\theta}$ and $H^{s,\theta}$. Thus, we want to
prove \begin{align}
\label{WMEmbedding1}
H^{s,\theta} \cdot H^{s,\theta} &\hookrightarrow H^{s,\theta}, \\
\label{WMEmbedding2}
H^{s-1,\theta-1} \cdot H^{s,\theta} &\hookrightarrow H^{s-1,\theta-1}.
\end{align}
These are special cases of the following. %
\begin{theorem}\label{WMProductEstimates} Let $n \ge 2$, $s > \frac{n}{2}$
and $\half < \theta \le s - \frac{n-1}{2}$. Then $$
H^{a,\alpha} \cdot H^{s,\theta} \hookrightarrow H^{a,\alpha} $$
for all $a, \alpha$ satisfying
\begin{align*}
0 &\le \alpha \le \theta,
\\
-s + \alpha &< a \le s.
\end{align*}
(Hence, by duality, for all $- \theta \le \alpha \le 0$ and $-s \le a < s
+ \alpha$.) \end{theorem}
The proof is achieved by interpolating between four different points in
the $(a,\alpha)$-plane. One of these points is $(s,\theta)$, which
corresponds to the estimate \eqref{WMEmbedding1}. We give the proof of the
latter here. The proofs of the remaining estimates are similar, and can be
found in the Appendix.

We may restate \eqref{WMEmbedding1} as follows: %
\begin{theorem}\label{AlgebraTheorem}
$H^{s,\theta}$ is an algebra if $n \ge 2$, $s > \frac{n}{2}$ and $\half <
\theta \le s - \frac{n-1}{2}$. \end{theorem}
For the proof we need the following ``Leibniz rule'', which is an
immediate consequence of the triangle inequality. %
\begin{lemma}\label{EllipticLambdaEstimates} If $\alpha > 0$, then $$
\Lambda^{\alpha}(uv) \precsim (\Lambda^{\alpha} u) v + u \Lambda^{\alpha}
v $$
for all $u$ and $v$ with $\widehat u, \widehat v \ge 0$. Moreover, the
same estimate holds with $\Lambda^{\alpha}$ replaced by either of the
operators $D^\alpha, D_+^\alpha$ or $\Lambda_{+}^{\alpha}$.
\end{lemma}
By Lemma \ref{EllipticLambdaEstimates}, the proof of Theorem
\ref{AlgebraTheorem} reduces to showing
\begin{align*}
H^{0,\theta} \cdot H^{s,\theta} &\hookrightarrow H^{0,\theta}. \\
\intertext{But by Lemma \ref{HyperbolicLambdaEstimates}, the latter
reduces to three estimates:}
H^{0,\theta} \cdot H^{s,0} &\hookrightarrow L^2, \\
L^2 \cdot H^{s,\theta} &\hookrightarrow L^2, \\
R^\theta( H^{0,\theta}, H^{s,\theta} ) &\hookrightarrow L^2. \end{align*}
The first one follows from H\"older's inequality, the energy embedding
\eqref{EnergyEmbedding} and the Sobolev embedding
\eqref{SpaceInftyEmbedding}; the second one holds by H\"older's inequality
and \eqref{InftyEmbedding}; the third one is a special case of Theorem F.
\subsection{The General Case}\label{WMDetailedProof} %
Here we prove the following.
\begin{theorem}\label{WMTheorem}
Let $n \ge 2$, $s > \frac{n}{2}$. Suppose \begin{gather*}
\half < \theta \le \min \Bigl(1, s-\frac{n-1}{2} \Bigr), \\
0 \le \varepsilon \le \min \Bigl(1 - \theta, s-\frac{n-1}{2} - \theta
\Bigr). \end{gather*}
Let $\Gamma : \R^N \to \R$ be smooth. Then there exist continuous
functions $g, h : [0,\infty) \to [0,\infty)$ such that \begin{equation}
\label{WMEstimate1}
\spacetimenorm{\Gamma(u) Q_0(u^J,u^K)}{s-1}{\theta+\varepsilon-1} \le g
\bigl( \spacetimenorm{u}{s}{\theta} \bigr) \Spacetimenorm{u}{s}{\theta}^2
\end{equation}
and
\begin{multline}
\label{WMEstimate2}
\spacetimenorm{\Gamma(u) Q_0(u^J,u^K) - \Gamma(U)
Q_0(U^J,U^K)}{s-1}{\theta+\varepsilon-1} \\
\le h \bigl( \spacetimenorm{u}{s}{\theta} + \spacetimenorm{U}{s}{\theta}
\bigr) \Spacetimenorm{u-U}{s}{\theta}
\end{multline}
for all $\R^N$-valued $u, U \in \scrH^{s,\theta}$ and $1 \le J,K \le N$.
\end{theorem}
As a consequence, we obtain part (a) of the Main Theorem.

We shall need the following.
\begin{lemma}\label{Q0Estimate}
If $0 \le \alpha \le 1$, then
$$
Q_0(\phi,\psi) \precsim D_+^{1-\alpha} D_-^{1-\alpha} (D_+^\alpha \phi'
D_+^\alpha \psi')
+ D_+ D_-^{1-\alpha} \phi' D_+^\alpha \psi' + D_+^\alpha \phi' D_+
D_-^{1-\alpha} \psi' $$
whenever $\phi \preceq \phi'$ and $\psi \preceq \psi'$. \end{lemma}
\begin{proof} The symbol of $Q_0$ is $q_0(\Xi,\Theta) \simeq
\innerprod{\Xi}{\Theta}$. Recall that $\innerprod{\cdot}{\cdot}$ denotes
the Minkowskian inner product on $\R^{1+n}$, while $\abs{\cdot}$ always
denotes the Euclidean norm. Since $$
\innerprod{\Xi+\Theta}{\Xi+\Theta} = \innerprod{\Xi}{\Xi} +
\innerprod{\Theta}{\Theta} + 2 \innerprod{\Xi}{\Theta}, $$
we have
$$
\abs{q_0(\Xi,\Theta)} \lesssim \abs{\innerprod{\Xi}{\Xi}} +
\abs{\innerprod{\Theta}{\Theta}} +
\abs{\innerprod{\Xi+\Theta}{\Xi+\Theta}}. $$
Take this to the power $1-\alpha$, and take the trivial estimate $$
\abs{q_0(\Xi,\Theta)} \lesssim \abs{\Xi} \abs{\Theta} $$
to the power $\alpha$. The product of the left hand sides of the resulting
inequalities is then bounded by the product of the right hand sides, and
keeping in mind that the symbol of $D_+ D_-$ is
$\abs{\innerprod{\Xi}{\Xi}}$, we get the desired estimate. \end{proof}
We first prove Theorem \ref{WMTheorem} in the case of constant $\Gamma$;
the general case is then reduced to this, by virtue of Theorem
\ref{WMProductEstimates} and the following result, which is an analogue of
the Moser inequality \eqref{MoserInequality}. %
\begin{theorem}\label{SchauderLemma}
Assume that $\Gamma \in C^\infty(\R^N)$ and $\Gamma(0) = 0$. If
$n,s,\theta$ are as in Theorem \ref{AlgebraTheorem} and $\theta \le 1$,
there exists a continuous function $g = g_{s,\theta} : [0,\infty) \to
[0,\infty)$ such that
$$
\spacetimenorm{\Gamma(u)}{s}{\theta}
\le g\bigl( \spacetimenorm{u}{n/2+\varepsilon}{\theta} \bigr)
\spacetimenorm{u}{s}{\theta}
$$
for all $\R^{N}$-valued $u \in H^{s,\theta}$, where $\varepsilon = \theta
- 1/2$.
\end{theorem}
This was proved in \cite{Se}.

Let us now prove Theorem \ref{WMTheorem}. Throughout the rest of this
section we assume that $n$, $s$, $\theta$ and $\varepsilon$ satsify the
hypotheses of Theorem \ref{WMTheorem}. %
\paragraph{Step 1.}
We assume $\Gamma = 1$; i.e., we prove
\begin{equation}\label{WMEmbedding}
Q_0 \bigl( \scrH^{s,\theta}, \scrH^{s,\theta} \bigr) \hookrightarrow
H^{s-1,\theta+\varepsilon-1}. \end{equation}
If we apply Lemma \ref{Q0Estimate} with $\alpha = \varepsilon$, and then
apply Lemma \ref{EllipticLambdaEstimates} to the first term on the right
hand side, we get $$
Q_0(\phi,\psi) \precsim D_-^{1-\alpha}
(D_+ \phi' D_+^\alpha \psi')
+ D_+ D_-^{1-\alpha} \phi' D_+^\alpha \psi' + \text{symmetric terms},
$$
where $\phi \preceq \phi'$ and $\psi \preceq \psi'$. Thus
\eqref{WMEmbedding} reduces to two estimates:
\begin{align*}
H^{s-1,\theta} \cdot H^{s-\varepsilon,\theta} &\hookrightarrow
H^{s-1,\theta},
\\
H^{s-1,\theta+\varepsilon-1} \cdot H^{s-\varepsilon,\theta}
&\hookrightarrow H^{s-1,\theta+\varepsilon-1}, \end{align*}
both of which are special cases of Theorem \ref{WMProductEstimates}. %
\paragraph{Step 2.}
We prove the theorem for general $\Gamma$. First, we write $$
\Gamma(u) Q_0(u^J,u^K)
= \bigl\{ \Gamma(u) - \Gamma(0) \bigr\} Q_0(u^J,u^K) + \Gamma(0)
Q_0(u^J,u^K). $$
Since
\begin{equation}\label{SpecialProductEstimate}
H^{s,\theta} \cdot H^{s-1,\theta+\varepsilon-1} \hookrightarrow
H^{s-1,\theta+\varepsilon-1} \end{equation}
by Theorem \ref{WMProductEstimates}, \eqref{WMEstimate1} reduces to
\begin{align}
\tag{\ref{WMEstimate1}a}
\spacetimenorm{\Gamma(u) - \Gamma(0)}{s}{\theta} &\le g \bigl(
\spacetimenorm{u}{s}{\theta} \bigr), \\
\tag{\ref{WMEstimate1}b}
\spacetimenorm{Q_0(u^J,u^K)}{s-1}{\theta+\varepsilon-1} &\lesssim
\Spacetimenorm{u^J}{s}{\theta} \Spacetimenorm{u^K}{s}{\theta}. \end{align}
The former holds by Theorem \ref{SchauderLemma}, the latter by Step 1.

To prove \eqref{WMEstimate2}, write
\begin{multline*}
\Gamma(u) Q_0(u^J,u^K) - \Gamma(U) Q_0(U^J,U^K) \\
= \bigl\{ \Gamma(u) - \Gamma(U) \bigr\} Q_0(u^J,u^K) + \Gamma(U) \bigl\{
Q_0(u^J- U^J,u^K) + Q_0(U^J,u^K - U^K) \bigr\}. \end{multline*}
The second term on the right hand side
is covered by the proof of \eqref{WMEstimate1}, while the first term
reduces, in view of \eqref{SpecialProductEstimate} and
(\ref{WMEstimate1}b), to the estimate
\begin{equation}\label{WMEstimate3}
\spacetimenorm{\Gamma(u) - \Gamma(U)}{s}{\theta} \le h \bigl(
\spacetimenorm{u}{s}{\theta} + \spacetimenorm{U}{s}{\theta} \bigr)
\Spacetimenorm{u-U}{s}{\theta}.
\end{equation}
But
\begin{multline*}
\Gamma(u) - \Gamma(U)
= \int_0^1 d\Gamma \bigl( (1-\lambda)U + \lambda u \bigr) \cdot (u - U) \,
d\lambda
\\
= \int_0^1 \left\{ d\Gamma \bigl( (1-\lambda)U + \lambda u \bigr) -
d\Gamma(0) \right\} \cdot (u - U) \, d\lambda + d\Gamma(0) \cdot (u - U),
\end{multline*}
and since $H^{s,\theta}$ is an algebra (Theorem \ref{AlgebraTheorem}), it
follows that
\begin{multline*}
\spacetimenorm{\Gamma(u) - \Gamma(U)}{s}{\theta} \\
\lesssim
\int_0^1 \spacetimenorm{ d\Gamma \bigl( (1-\lambda)U + \lambda u \bigr) -
d\Gamma(0) }{s}{\theta} \spacetimenorm{u - U}{s}{\theta} \, d\lambda +
\abs{d\Gamma(0)} \spacetimenorm{u - U}{s}{\theta}. \end{multline*}
Thus, \eqref{WMEstimate3} follows after another application of Theorem
\ref{SchauderLemma}.
\section{Main Estimates for \eqref{MKGType}/\eqref{YMType}}\label{YMProof}
Here we prove part (b) of the Main Theorem: \eqref{MKGType}/\eqref{YMType}
are locally well-posed for initial data in $H^{s} \times H^{s-1}(\R^{n})$
for $s > \frac{n-2}{2}$ and $n \ge 4$.

The full details of the proof appear in section \ref{YMDetailedProof}. The
reader who wants to get the gist of the argument, without getting bogged
down in technicalities, is advised to read first the informal discussion
in section \ref{YMInformalProof}.
\subsection{Informal Proof}\label{YMInformalProof} %
In order to prove part (b) of the Main Theorem by iteration in the space
$\scrH^{s,\theta}$, we would need two types of estimates: \begin{align}
\label{YMA}
D^{-1} Q_{ij} \bigl( \scrH^{s,\theta}, \scrH^{s,\theta} \bigr)
&\hookrightarrow H^{s-1,\theta-1}, \\
\label{YMB}
Q_{ij} \bigl( D^{-1} \scrH^{s,\theta}, \scrH^{s,\theta} \bigr)
&\hookrightarrow H^{s-1,\theta-1}, \end{align}
for all $s > \frac{n-2}{2}$, $n \ge 4$ and some $\theta = \theta(s,n) >
\half$.

We shall use the following.
\begin{lemma}\label{QijEstimate}
The estimate
$$
Q_{ij} (\phi,\psi) \preceq
D^\half D_-^\half (D^\half \phi' D^\half \psi') + D^\half( D^\half
D_-^\half \phi' D^\half \psi') + D^\half( D^\half \phi' D^\half D_-^\half
\psi') $$
holds whenever $\phi \preceq \phi'$ and $\psi \preceq \psi'$. \end{lemma}
For the proof, see \cite{Kl-Ma4}.

In view of Lemma \ref{QijEstimate} and Remark \ref{PositivityRemark},
proving \eqref{YMA} reduces to proving, if we replace $D^{-1}$ by
$\Lambda^{-1}$ and ignore the difference between $\scrH^{s,\theta}$ and
$H^{s,\theta}$, \begin{align}
\label{YMA1}
H^{s-\half,\theta} \cdot H^{s-\half,\theta} &\hookrightarrow
H^{s-\frac{3}{2},\theta-\half}, \\
\label{YMA2}
H^{s-\half,\theta-\half} \cdot H^{s-\half,\theta} &\hookrightarrow
H^{s-\frac{3}{2},\theta-1}, \\
\intertext{Similarly, \eqref{YMB} can be reduced to}
\label{YMB1}
H^{s+\half,\theta} \cdot H^{s-\half,\theta} &\hookrightarrow
H^{s-\half,\theta-\half}, \\
\label{YMB2}
H^{s+\half,\theta-\half} \cdot H^{s-\half,\theta} &\hookrightarrow
H^{s-\half,\theta-1}, \\
\label{YMB3}
H^{s+\half,\theta} \cdot H^{s-\half,\theta-\half} &\hookrightarrow
H^{s-\half,\theta-1}. \end{align}
\begin{theorem}\label{NaiveYMTheorem}
The estimates \eqref{YMA1}--\eqref{YMB1} and \eqref{YMB3} hold for all $s
> \frac{n-2}{2}$, $n \ge 4$ and $\half < \theta \le s - \frac{n-3}{2}$.
However, \eqref{YMB2} fails if $s < \frac{n}{2} - \theta$. \end{theorem}
\begin{remarkNoLabel}
The condition $\theta - \half \le s - \frac{n-2}{2}$ is necessary by
scaling. Therefore, in view of the above theorem, \eqref{YMB2} cannot hold
unless $s \ge \frac{n}{2} - \frac{3}{4}$. In fact, we expect that
\eqref{YMA1}--\eqref{YMB3} are all true for $s > \frac{n}{2} -
\frac{3}{4}$, $n \ge 4$ with $\theta = s - \frac{n-3}{2}$. This has been
verified in dimension $n = 3$; see Cuccagna \cite{Cu} and Keel-Tao
\cite{Ke-Ta3}. We do not pursue this question. \end{remarkNoLabel}
\begin{remarkNoLabel}
The failure of our attempt to iterate in $\scrH^{s,\theta}$ when $s <
\frac{n}{2} - \theta$ is not due to any loss of information through the
use of Lemma \ref{QijEstimate}: The proof of the last statement in Theorem
\ref{NaiveYMTheorem} (see the Appendix) shows that \eqref{YMB} also fails
for such $s$. \end{remarkNoLabel}
The first (positive) statement in Theorem \ref{NaiveYMTheorem} is proved
in section \ref{YMDetailedProof}, and the second (negative) statement is
proved in the Appendix. The following heuristic arguments should convince
the reader that the result is reasonable.

Observe that if we consider the idealized case $\theta = \half$, then the
estimates \eqref{YMA1}--\eqref{YMB3} are all of the form (using duality if
necessary)
$$
H^{a,\theta} \cdot H^{b,\theta} \hookrightarrow H^{-c,0}. $$
The latter is morally equivalent to a product estimate for two solutions
of the homogeneous wave equation:
$$
\twonorm{D^{-c}(uv)}{} \lesssim \Sobnorm{f}{a} \Sobnorm{g}{b}, $$
where $\square u = \square v = 0$, $(u,\partial_{t} u) \init = (f,0)$ and
$(v,\partial_{t} v) \init = (g,0)$. By Theorem C, a necessary condition
for this estimate to hold is
\begin{equation}\label{WaveEstimateCondition}
a + b \ge \half.
\end{equation}
Note that if $a + b + c = \frac{n-1}{2}$, then
\eqref{WaveEstimateCondition} is equivalent to $c \le \frac{n-2}{2}$.

Let us reexamine our estimates in the light of condition
\eqref{WaveEstimateCondition}.
Taking $s = \frac{n-2}{2}$ and $\theta = \half$, and using duality where
necessary, \eqref{YMA1}--\eqref{YMB3} reduce to \begin{align}
\tag{\ref{YMA1}$'$}
H^{\frac{n-3}{2},\theta} \cdot H^{\frac{n-3}{2},\theta} &\hookrightarrow
H^{\frac{n-5}{2},0},
\\
\tag{\ref{YMA2}$'$}
H^{\frac{5-n}{2},\theta} \cdot H^{\frac{n-3}{2},\theta} &\hookrightarrow
H^{\frac{3-n}{2},0},
\\
\tag{\ref{YMB1}$'$}
H^{\frac{n-1}{2},\theta} \cdot H^{\frac{n-3}{2},\theta} &\hookrightarrow
H^{\frac{n-3}{2},0},
\\
\tag{\ref{YMB2}$'$}
H^{\frac{3-n}{2},\theta} \cdot H^{\frac{n-3}{2},\theta} &\hookrightarrow
H^{\frac{1-n}{2},0},
\\
\tag{\ref{YMB3}$'$}
H^{\frac{n-1}{2},\theta} \cdot H^{\frac{3-n}{2},\theta} &\hookrightarrow
H^{\frac{3-n}{2},0}.
\end{align}
Condition \eqref{WaveEstimateCondition} is satisfied in all of the above
except (\ref{YMB2}$'$), where $a + b = 0$ (and $c = \frac{n-1}{2}$). The
latter estimate is therefore far from being true (it is half a derivative
off the mark). On the other hand, since $n \ge 4$, it is easily checked
that the other four estimates above are in fact true by Theorem F, if we
take $\theta > \half$. %
Let us now take a closer look at the estimate which fails, namely
\eqref{YMB2}. By Lemma \ref{EllipticLambdaEstimates}, this reduces to
\begin{align}
\notag
H^{1,\theta-\half} \cdot H^{s-\half,\theta} &\hookrightarrow
H^{0,\theta-1}, \\
\label{YMB21}
H^{s+\half,\theta-\half} \cdot H^{0,\theta} &\hookrightarrow
H^{0,\theta-1}. \end{align}
The former is true for $s > \frac{n-2}{2}$, $n \ge 4$ (see section
\ref{YMDetailedProof} for the proof), while the latter fails for $s <
\frac{n}{2} - \theta$ (see the Appendix).

For simplicity, throughout the remainder of section \ref{YMInformalProof}
we will only consider the idealized case where $s = \frac{n-2}{2}$ and
$\theta = \half$. \emph{Thus, when we say that an estimate holds, we mean
up to a logarithmic divergence}. The informal arguments in this section
are easily made rigorous (see section \ref{YMDetailedProof}).

Since we are assuming $s = \frac{n-2}{2}$ and $\theta = \half$, the
problematic estimate \eqref{YMB21} reads
$$
H^{\frac{n-1}{2},0} \cdot H^{0,\theta} \hookrightarrow H^{0,-\theta}. $$
An easy way to fix the problem with this estimate is to replace
$H^{\frac{n-1}{2},0}$ on the left hand side with $$
H^{\frac{n-1}{2},0} \cap \mixed{1}{\infty}. $$
In other words, we claim that
\begin{equation}\label{NewEstimate}
\bigl( H^{\frac{n-1}{2},0} \cap \mixed{1}{\infty} \bigr) \cdot
H^{0,\theta} \hookrightarrow H^{0,-\theta}. \end{equation}
This is a trivial consequence of energy estimates and H\"older's
inequality. Indeed, by \eqref{SpecialEmbeddingA}, \begin{align*}
H^{0,\theta} &\hookrightarrow \mixed{\infty}{2}. \\
\intertext{The dual of this embedding is} \mixed{1}{2} &\hookrightarrow
H^{0,-\theta}. \end{align*}
Since
$$
\mixed{1}{\infty} \cdot \mixed{\infty}{2} \hookrightarrow \mixed{1}{2}
$$
by H\"older's inequality, we obtain \eqref{NewEstimate}.

This suggests taking
\begin{equation}\label{FirstGuessSpace}
\scrH^{s,\theta} \cap \Lambda^{\half} \Lambda_{-}^{-\half} \bigl(
\Mixed{1}{\infty} \bigr)
\end{equation}
as our iteration space. This works for systems of the type
\eqref{MKGType}, but leads to problems for \eqref{YMType} (cf. Remarks
\ref{MKGRemark2} and \ref{YMRemark} below). A better choice turns out to be
\begin{equation}\label{NewSpace}
\X^{s} = \scrH^{s,\theta} \cap \Lambda^{-\gamma} \Lambda_{-}^{-\half}
\bigl( \Mixed{1}{2n} \bigr),
\end{equation}
where $\gamma > 0$ is sufficiently small. Since we are assuming $s =
\frac{n-2}{2}$ and $\theta = \half$, we will take $\gamma = 0$ here.

We are now faced with the task of proving \begin{align*}
D^{-1} Q_{ij}( \X^{s}, \X^{s}) &\hookrightarrow \Lambda_{+} \Lambda_{-}
\X^{s},
\\
Q_{ij}(D^{-1} \X^{s}, \X^{s}) &\hookrightarrow \Lambda_{+} \Lambda_{-}
\X^{s},
\end{align*}
where $\X^{s}$ is given by \eqref{NewSpace}. In fact, we can prove
\begin{align}
\label{YMC}
D^{-1} Q_{ij}( \scrH^{s,\theta}, \scrH^{s,\theta}) &\hookrightarrow
\Lambda_{+} \Lambda_{-} \X^{s},
\\
\label{YMD}
Q_{ij}(D^{-1} \X^{s}, \scrH^{s,\theta}) &\hookrightarrow \Lambda_{+}
\Lambda_{-} \X^{s}.
\end{align}
In view of the definition of $\X^{s}$, \eqref{YMC} is equivalent to
\begin{align}
\label{YMC1}
D^{-1} Q_{ij}( \scrH^{s,\theta}, \scrH^{s,\theta}) &\hookrightarrow
H^{s-1,\theta-1},
\\
\label{YMC2}
D^{-1} Q_{ij}( \scrH^{s,\theta}, \scrH^{s,\theta}) &\hookrightarrow
\Lambda_{+} \Lambda_{-}^{\half} \bigl( \Mixed{1}{2n} \bigr). \end{align}
Similarly, \eqref{YMD} is equivalent to
\begin{align}
\label{YMD1}
Q_{ij}(D^{-1} \X^{s}, \scrH^{s,\theta}) &\hookrightarrow H^{s-1,\theta-1}
\\
\label{YMD2}
Q_{ij}(D^{-1} \X^{s}, \scrH^{s,\theta}) &\hookrightarrow \Lambda_{+}
\Lambda_{-}^{\half} \bigl( \Mixed{1}{2n} \bigr). \end{align}
Observe that \eqref{YMC1} follows from Theorem \ref{NaiveYMTheorem}. Also,
\eqref{YMD1} is the estimate that motivated the introduction of the new
space $\X^{s}$.
It therefore remains to prove \eqref{YMC2} and \eqref{YMD2}. %
\begin{remark}\label{MKGRemark1}
For the system \eqref{MKGType} the estimate \eqref{YMD2} is not needed. In
fact, it is clear from the special structure of \eqref{MKGType} that if
\eqref{YMC} and \eqref{YMD1} are true, then we can iterate in the space $$
\bigl\{(u,v) : u \in \X^{s}, v \in \scrH^{s,\theta} \bigr\}. $$
See also Remark \ref{MKGRemark2} below.
\end{remark}
\subsubsection{Informal Proof of \eqref{YMC2}} %
By Lemma \ref{QijEstimate}, \eqref{YMC2} can be reduced to proving (again
we ignore the difference between $D^{-1}$ and $\Lambda^{-1}$) \begin{align}
\label{YMC21}
H^{s-\half,\theta} \cdot H^{s-\half,\theta} &\hookrightarrow
\Lambda^{\half} \Lambda_{+} \Mixed{1}{2n}, \\
\label{YMC22}
H^{s-\half,\theta-\half} \cdot H^{s-\half,\theta} &\hookrightarrow
\Lambda^{\half} \Lambda_{+} \Lambda_{-}^{\half} \Mixed{1}{2n}. \end{align}
Clearly, $\Lambda^{\half} \Lambda_{+}$ may be replaced by
$\Lambda^{\frac{3}{2}}$ on the right hand side of both estimates.

First, \eqref{YMC21} holds by Theorem E (again up to logarithmic
divergence).

To prove \eqref{YMC22}, we use the following special case of Proposition
\ref{DualKlMaEmbedding} (valid since $n \ge 4$):
\begin{equation}\label{SpecialCaseKlMaEmbedding}
\Mixednorm{u}{1}{2n} \lesssim \bigMixednorm{ \Lambda^{\frac{n-2}{2}}
\Lambda_{-}^{\theta} u}{1}{2}.
\end{equation}
Thus, since we are assuming $s = \frac{n-2}{2}$ and $\theta = \half$,
\eqref{YMC22} reduces to
$$
H^{\frac{n-3}{2},0} \cdot H^{\frac{n-3}{2},\theta} \hookrightarrow
\Lambda^{\frac{5-n}{2}} \Mixed{1}{2}.
$$
For simplicity, we consider only the case $n = 4$; the latter estimate
then reads
\begin{equation}\label{YMC22four}
H^{\half,0} \cdot H^{\half,\theta}
\hookrightarrow
\Lambda^{\half} \Mixed{1}{2}.
\end{equation}
Set $I = \bigMixednorm{\Lambda^{-\half}(uv)}{1}{2}$. We have to prove $I
\lesssim \spacetimenorm{u}{\half}{0}
\spacetimenorm{v}{\half}{\theta}$.
Since the norms involved only depend on the size of the Fourier transform
and are compatible with the relation $\preceq$, we may assume that
$\widehat u, \widehat v \ge 0$.

Now write
\begin{equation}\label{Trick}
\Lambda^{-\half}(uv) = \Lambda^{-\half} \bigl( u \Lambda^{\half}
\Lambda^{-\half} v \bigr),
\end{equation}
and apply the estimate (valid for any $\alpha > 0$)
\begin{equation}\label{VariationEllipticLambdaEstimate}
u \Lambda^{\alpha} v \lesssim \Lambda^{\alpha}(uv) + \Lambda^{\alpha} u
\cdot v.
\end{equation}
(This holds by the triangle inequality.) Thus, $I \lesssim I_{1} + I_{2}$,
where
\begin{align*}
I_{1} &= \bigMixednorm{u \Lambda^{-\half} v}{1}{2}, \\
I_{2} &= \bigMixednorm{\Lambda^{-\half} \bigl( \Lambda^{\half} u \cdot
\Lambda^{-\half} v \bigr)}{1}{2}.
\end{align*}

For $I_{1}$, \eqref{ModifiedMixednormBound} and H\"older's inequality give
$$
I_{1} \le \mixednorm{u}{2}{\frac{8}{3}}
\bigmixednorm{ \Lambda^{-\half} v}{2}{8}. $$
Since $n = 4$, we have
$$
\mixednorm{u}{2}{\frac{8}{3}} \lesssim \spacetimenorm{u}{\half}{0} $$
by Sobolev embedding, and
$$
\bigmixednorm{ \Lambda^{-\half} v}{2}{8} \lesssim
\spacetimenorm{v}{1}{\theta} $$
by the following special case of Theorem D:
\begin{equation}\label{SpecialStrichartz1}
H^{1,\theta} \hookrightarrow \mixed{2}{8} \qquad (n = 4). \end{equation}

For $I_{1}$, \eqref{ModifiedMixednormBound} and Sobolev embedding,
followed by H\"older's inequality, gives $$
I_{2} \lesssim \bigmixednorm{\Lambda^{\half} u \cdot \Lambda^{-\half}
v}{1}{\frac{8}{5}} \le \bigtwonorm{\Lambda^{\half} u}{}
\bigmixednorm{\Lambda^{-\half} v}{2}{8}. $$
And by \eqref{SpecialStrichartz1} again, the right hand side is $\lesssim
\spacetimenorm{u}{\half}{0}
\spacetimenorm{v}{\half}{\theta}$.

This concludes the discussion of \eqref{YMC2}.
\begin{remark}\label{MKGRemark2}
An inspection of the above arguments reveals that up to this point we
could just as well have been working in the space \eqref{FirstGuessSpace}.
For \eqref{YMC} still holds if we let $\X^{s}$ be defined by
\eqref{FirstGuessSpace}, and with essentially the same proof as above
(only a few obvious modifications are needed). It is only when we try to
prove \eqref{YMD} that we run into problems if we choose the space
\eqref{FirstGuessSpace}.
\end{remark}
\subsubsection{Informal Proof of \eqref{YMD2}} %
By Lemma \ref{QijEstimate}, \eqref{YMD2} reduces to \begin{align}
\label{YMD21}
\Lambda^{-\half} \X^{s} \cdot \scrH^{s-\half,\theta} &\hookrightarrow
\Lambda^{\half} \Mixed{1}{2n},
\\
\label{YMD22}
\Lambda^{-\half} \Lambda_{-}^{\half} \X^{s} \cdot \scrH^{s-\half,\theta}
&\hookrightarrow
\Lambda^{\half} \Lambda_{-}^{\half}\Mixed{1}{2n}, \\
\label{YMD23}
\Lambda^{-\half} \X^{s} \cdot \scrH^{s-\half,\theta-\half} &\hookrightarrow
\Lambda^{\half} \Lambda_{-}^{\half} \Mixed{1}{2n}. \end{align}
\begin{remark}\label{YMRemark}
It is the estimate \eqref{YMD22} which necessitates the use of the space
\eqref{NewSpace} rather than \eqref{FirstGuessSpace} (more precisely, the
problem comes up with the estimate \eqref{YMD221} below, which derives
from \eqref{YMD22}). The other two estimates, \eqref{YMD21} and
\eqref{YMD23}, are easily seen to be true also with $\X^{s}$ given by
\eqref{FirstGuessSpace}, and with essentially the same proof as below.
\end{remark}
For \eqref{YMD21} it suffices to prove
$$
H^{s+\half,\theta} \cdot H^{s-\half,\theta} \hookrightarrow
\Lambda^{\half} \Mixed{1}{2n}.
$$
Equivalently,
\begin{equation}\label{AsymmetricYMEstimate}
\bigMixednorm{\Lambda^{-\half} \bigl( \Lambda^{-\half} u \cdot
\Lambda^{\half} v \bigr)}{1}{2n} \lesssim \spacetimenorm{u}{s}{\theta}
\spacetimenorm{v}{s}{\theta}. \end{equation}
By \eqref{VariationEllipticLambdaEstimate}, $$
\bigMixednorm{\Lambda^{-\half} \bigl( \Lambda^{-\half} u \cdot
\Lambda^{\half} v \bigr)}{1}{2n}
\lesssim \bigmixednorm{\Lambda^{-\half} u \cdot v}{1}{2n} +
\bigmixednorm{\Lambda^{-\half} (u v)}{1}{2n}, $$
where we also used \eqref{ModifiedMixednormBound}. By Theorem E,
$$
\bigmixednorm{\Lambda^{-\half} (u v)}{1}{2n} \lesssim
\spacetimenorm{u}{s}{\theta} \spacetimenorm{v}{s}{\theta}. $$
Take $n = 4$ for simplicity. Then by H\"older's inequality and Sobolev
embedding, $$
\bigmixednorm{\Lambda^{-\half} u \cdot v}{1}{8} \le
\bigmixednorm{\Lambda^{-\half} u}{2}{\infty} \mixednorm{v}{2}{8}
\lesssim \mixednorm{u}{2}{8} \mixednorm{v}{2}{8}. $$
Now apply \eqref{SpecialStrichartz1}. This finishes the proof of
\eqref{YMD21}.

By applying \eqref{SpecialCaseKlMaEmbedding}, we reduce \eqref{YMD22} to
$$
\Lambda^{-\half} \Lambda_{-}^{\half} \X^{s} \cdot \scrH^{s-\half,\theta}
\hookrightarrow
\Lambda^{\frac{3-n}{2}} \Mixed{1}{2}.
$$
By Lemma \ref{EllipticLambdaEstimates}, it suffices to prove (recall that
$s = \frac{n-2}{2}$):
\begin{align}
\label{YMD221}
\Lambda^{\frac{n-4}{2}} \Lambda_{-}^{\half} \X^{s} \cdot
H^{\frac{n-3}{2},\theta}
&\hookrightarrow
\Mixed{1}{2},
\\
\label{YMD222}
\Lambda^{-\half} \Lambda_{-}^{\half} \X^{s} \cdot H^{0,\theta}
&\hookrightarrow
\Mixed{1}{2}.
\end{align}
For simplicity, we take $n = 4$ again. Then \eqref{YMD221} becomes
\begin{align*}
\Lambda_{-}^{\half} \X^{s} \cdot H^{\half,\theta} &\hookrightarrow
\Mixed{1}{2},
\\
\intertext{and in view of the definition of $\X^{s}$, it is enough to
prove}
\Mixed{1}{8} \cdot H^{\half,\theta}
&\hookrightarrow
\Mixed{1}{2}.
\end{align*}
Equivalently,
$$
\Mixednorm{uv}{1}{2} \lesssim \Mixednorm{u}{1}{8}
\spacetimenorm{v}{\half}{\theta}.
$$
In view of Lemma \ref{DensityLemma}, we may assume that $v \in \Schwartz$
and $\widehat v \ge 0$. Therefore, by Proposition \ref{MixednormHolder},
$$
\Mixednorm{uv}{1}{2} \le \Mixednorm{u}{1}{8}
\mixednorm{v}{\infty}{\frac{8}{3}},
$$
and by \eqref{SpecialEmbeddingH} (since $n = 4$), $$
\mixednorm{v}{\infty}{\frac{8}{3}} \lesssim
\spacetimenorm{v}{\half}{\theta}. $$
This proves \eqref{YMD221}. As for \eqref{YMD222}, it is enough to prove
$$
\Lambda^{-\half} \Mixed{1}{8} \cdot H^{0,\theta} \hookrightarrow
\Mixed{1}{2}.
$$
Reasoning as above, we have
$$
\Mixednorm{uv}{1}{2} \le \Mixednorm{u}{1}{\infty} \mixednorm{v}{\infty}{2}.
$$
By Sobolev embedding, or more accurately by \eqref{SpecialEmbeddingJ}, $$
\Mixednorm{u}{1}{\infty} \lesssim \bigMixednorm{\Lambda^{\half} u}{1}{8}.
$$
By the energy embedding \eqref{EnergyEmbedding}, $$
\mixednorm{v}{\infty}{2} \lesssim \spacetimenorm{v}{0}{\theta}. $$
This finishes the proof of \eqref{YMD22}.

As in the proof of \eqref{YMD22}, by applying
\eqref{SpecialCaseKlMaEmbedding} followed by Lemma
\ref{EllipticLambdaEstimates}, \eqref{YMD23} reduces to proving (recall
that $s = \frac{n-2}{2}$ and $\theta = \half$) \begin{align}
\label{YMD231}
\Lambda^{\frac{n-4}{2}} \X^{s} \cdot
H^{\frac{n-3}{2},0}
&\hookrightarrow
\Mixed{1}{2},
\\
\label{YMD232}
\Lambda^{-\half} \X^{s} \cdot
L^{2}
&\hookrightarrow
\Mixed{1}{2}.
\end{align}
Again we take $n = 4$. For \eqref{YMD231} it is then enough to prove $$
H^{1,\theta} \cdot H^{\half,0} \hookrightarrow \Mixed{1}{2}. $$
In view of \eqref{ModifiedMixednormBound}, we may replace $\Mixed{1}{2}$
on the right by $\mixed{1}{2}$, and by H\"older's inequality,
$$
\mixednorm{uv}{1}{2} \le \mixednorm{u}{2}{8} \mixednorm{v}{2}{\frac{8}{3}}
$$
To the first factor on the right we apply \eqref{SpecialStrichartz1}, to
the second factor we apply Sobolev embedding.

For \eqref{YMD232} it is enough to prove $$
H^{\frac{3}{2},\theta} \cdot L^{2} \hookrightarrow \Mixed{1}{2}. $$
As above, we simply note that by H\"older, $$
\mixednorm{uv}{1}{2} \le \mixednorm{u}{2}{\infty} \twonorm{v}{}.
$$
Apply \eqref{SpecialEmbeddingD} (with $n = 4$) to the first factor on the
right. This finishes the proof of \eqref{YMD23}. %
\subsection{Proof of Main Theorem, part (b)}\label{YMDetailedProof} %
We shall prove the following.
\begin{theorem}\label{YMTheorem}
Let $n \ge 4$, $s > \frac{n-2}{2}$. Assume that $\theta$ and $\varepsilon$
satisfy
\begin{gather*}
\half < \theta \le \min \left( \frac{3}{4}, \half + s-\frac{n-2}{2}
\right), \\
0 \le \varepsilon \le \frac{1}{8} \min \left( \frac{3}{2} - 2\theta, s -
\frac{n-2}{2} + 1 - 2\theta \right), \end{gather*}
and let $\gamma = \theta - \half + 3\varepsilon$. Let $\X^s$ be the Banach
space given by the norm
$$
\norm{u} = \Spacetimenorm{u}{s}{\theta} + \bigMixednorm{\Lambda^{\gamma}
\Lambda_-^{\half} u}{1}{2n}.
$$
Then
\begin{align}
\label{YME}
D^{-1}Q_{ij} \bigl( \scrH^{s,\theta}, \scrH^{s,\theta} \bigr)
&\hookrightarrow \Lambda_+ \Lambda_-^{1-\varepsilon} \X^s, \\
\label{YMF}
Q_{ij} \bigl( D^{-1} \X^s, \scrH^{s,\theta} \bigr) &\hookrightarrow
\Lambda_+ \Lambda_-^{1-\varepsilon} \X^s. \end{align}
\end{theorem}
This implies part (b) of the Main Theorem, in view of Theorem
\ref{SpecializedWellPosednessTheorem}, since conditions (a)--(e) of the
latter are satisfied by the space $\X^{s}$ (condition (a) holds by
Proposition \ref{SecondEnergyEmbeddingProposition}; condition (b) is
obviously satisfied; conditions (c) and (d) follow from Theorem
\ref{BasicConditionsTheorem}, in view of Proposition
\ref{BasicSpacesProposition}; finally, condition (e) holds by Proposition
\ref{TimeCutOffProposition}).

By the definition of $\X^{s}$, \eqref{YME} is equivalent to two estimates:
\begin{align*}
D^{-1}Q_{ij} \bigl( \scrH^{s,\theta}, \scrH^{s,\theta} \bigr)
&\hookrightarrow H^{s-1,\theta+\varepsilon-1}, \\
D^{-1}Q_{ij} \bigl( \scrH^{s,\theta}, \scrH^{s,\theta} \bigr)
&\hookrightarrow \Lambda^{-\gamma}
\Lambda_+ \Lambda_-^{\half-\varepsilon} \Mixed{1}{2n}, \\
\intertext{and \eqref{YMF} is equivalent to} Q_{ij} \bigl( D^{-1} \X^s,
\scrH^{s,\theta} \bigr) &\hookrightarrow H^{s-1,\theta+\varepsilon-1}, \\
Q_{ij} \bigl( D^{-1} \X^s, \scrH^{s,\theta} \bigr) &\hookrightarrow
\Lambda^{-\gamma}
\Lambda_+ \Lambda_-^{\half-\varepsilon} \Mixed{1}{2n}. \end{align*}
We split these four estimates into what we call high and low frequency
cases. The high frequency estimates are the ones obtained by replacing
$D^{-1}$ by $\Lambda^{-1}$:
\begin{align}
\label{YME1}
\Lambda^{-1}Q_{ij} \bigl( \scrH^{s,\theta}, \scrH^{s,\theta} \bigr)
&\hookrightarrow H^{s-1,\theta+\varepsilon-1}, \\
\label{YME2}
\Lambda^{-1}Q_{ij} \bigl( \scrH^{s,\theta}, \scrH^{s,\theta} \bigr)
&\hookrightarrow \Lambda^{-\gamma}
\Lambda_+ \Lambda_-^{\half-\varepsilon} \Mixed{1}{2n}, \\
\label{YMF1}
Q_{ij} \bigl( \Lambda^{-1} \X^s, \scrH^{s,\theta} \bigr) &\hookrightarrow
H^{s-1,\theta+\varepsilon-1}, \\
\label{YMF2}
Q_{ij} \bigl( \Lambda^{-1} \X^s, \scrH^{s,\theta} \bigr) &\hookrightarrow
\Lambda^{-\gamma}
\Lambda_+ \Lambda_-^{\half-\varepsilon} \Mixed{1}{2n}. \end{align}
In the low frequency estimates, $D^{-1}$ is replaced by $\Lambda^{-M}
D^{-1}$, where $M > 0$ can be chosen arbitrarily large. In view of the
trivial estimates\footnote{These follow from the fact that the symbol of
$Q_{ij}$ is bounded in absolute value by $\abs{\xi \wedge \eta}$, where
$\xi \wedge \eta$ is the exterior product of vectors in $\R^n$. We have
$\abs{\xi \wedge \eta} \le \abs{\xi} \abs{\eta}$, and since $\xi \wedge
\eta = \xi \wedge (\xi + \eta) = (\xi + \eta) \wedge \eta$, we also have
$\abs{\xi \wedge \eta} \le \abs{\xi} \abs{\xi + \eta}, \abs{\xi + \eta}
\abs{\eta}$. By combining these we get the desired estimates.}
\begin{align*}
Q_{ij}(\phi,\psi) &\precsim D ( D^\half \phi' D^\half \psi'), \\
Q_{ij}(\phi,\psi) &\precsim D^\half ( D \phi' D^\half \psi'), \end{align*}
where $\phi \preceq \phi'$ and $\psi \preceq \psi'$, the low frequency
estimates reduce to \begin{align}
\label{YME1low}
H^{s-\half,\theta} \cdot H^{s-\half,\theta} &\hookrightarrow H^{-M,0},
\\
\label{YME2low}
H^{s-\half,\theta} \cdot H^{s-\half,\theta} &\hookrightarrow \Lambda^{-M}
\Mixed{1}{2n}, \\
\label{YMF1low}
H^{M,\theta} \cdot H^{s-\half,\theta}
&\hookrightarrow H^{s-\half,0},
\\
\label{YMF2low}
H^{M,\theta} \cdot H^{s-\half,\theta}
&\hookrightarrow \Lambda^{1-\gamma} \Mixed{1}{2n}, \end{align}
where $M > 0$ can be taken arbitrarily large.

The estimates \eqref{YME2low} and \eqref{YMF2low} hold by Theorem E, while
\eqref{YME1low} and \eqref{YMF1low} are special cases of the following
theorem, which is essentially a corollary of Theorem F (see the Appendix
for the proof). %
\begin{theorem}\label{YMProductTheorem}
Let $n \ge 4$ and $\theta > \half$. Then $$
H^{a,\theta} \cdot H^{b,\theta} \hookrightarrow H^{-c,0}. $$
for all $a,b,c$ satisfying
\begin{align*}
a,b &\ge -c,
\\
a + b &\ge \half,
\\
a + b + c &\ge \frac{n-1}{2}.
\end{align*}
\end{theorem}
Let us now turn to the proofs of \eqref{YME1}--\eqref{YMF2}. By Lemma
\ref{QijEstimate}, \eqref{YME1} reduces to \begin{align}
\label{YME11}
\scrH^{s-\half,\theta} \cdot \scrH^{s-\half,\theta} &\hookrightarrow
H^{s-\frac{3}{2},\theta+\varepsilon-\half}, \\
\label{YME12}
\scrH^{s-\half,\theta-\half} \cdot \scrH^{s-\half,\theta} &\hookrightarrow
H^{s-\frac{3}{2},\theta+\varepsilon-1}; \end{align}
\eqref{YME2} reduces to
\begin{align}
\label{YME21}
\scrH^{s-\half,\theta} \cdot \scrH^{s-\half,\theta} &\hookrightarrow
\Lambda^{\half-\gamma}
\Lambda_+ \Lambda_-^{-\varepsilon} \Mixed{1}{2n}, \\
\label{YME22}
\scrH^{s-\half,\theta-\half} \cdot \scrH^{s-\half,\theta} &\hookrightarrow
\Lambda^{\half-\gamma}
\Lambda_+ \Lambda_-^{\half-\varepsilon} \Mixed{1}{2n}; \end{align}
\eqref{YMF1} reduces to
\begin{align}
\label{YMF11}
\scrH^{s+\half,\theta} \cdot \scrH^{s-\half,\theta} &\hookrightarrow
H^{s-\half,\theta+\varepsilon-\half}, \\
\label{YMF12}
\Lambda^{-\half} \Lambda_{-}^{\half} \X^{s} \cdot \scrH^{s-\half,\theta}
&\hookrightarrow H^{s-\half,\theta+\varepsilon-1}, \\
\label{YMF13}
\scrH^{s+\half,\theta} \cdot \scrH^{s-\half,\theta} &\hookrightarrow
H^{s-\half,\theta+\varepsilon-1}; \end{align}
and \eqref{YMF2} reduces to
\begin{align}
\label{YMF21}
\scrH^{s+\half,\theta} \cdot \scrH^{s-\half,\theta} &\hookrightarrow
\Lambda^{-\half-\gamma} \Lambda_+ \Lambda_-^{-\varepsilon} \Mixed{1}{2n},
\\
\label{YMF22}
\Lambda^{-\half} \Lambda_{-}^{\half} \X^s \cdot \scrH^{s-\half,\theta}
&\hookrightarrow \Lambda^{-\half-\gamma} \Lambda_+
\Lambda_-^{\half-\varepsilon} \Mixed{1}{2n}, \\
\label{YMF23}
\scrH^{s+\half,\theta} \cdot \scrH^{s-\half,\theta-\half} &\hookrightarrow
\Lambda^{-\half-\gamma} \Lambda_+ \Lambda_-^{\half-\varepsilon}
\Mixed{1}{2n}. \end{align}
For the proofs, we need a few technical lemmas. %
\begin{lemma}\label{TrivialLambdaMinusEstimate} Let $\alpha > 0$. Then
$$
\Lambda_{-}^{\alpha}(uv) \precsim \Lambda_{+}^{\alpha} u
\Lambda_{+}^{\alpha} v
$$
for all $u$ and $v$ with $\widehat u, \widehat v \ge 0$. \end{lemma}
The trivial proof is omitted.
\begin{lemma}\label{YetAnotherLambdaMinusEstimate} Let $\alpha,\beta \ge
0$. Then
$$
\Lambda_{-}^{\alpha}(uv) \precsim (\Lambda_{+}^{\alpha} u)v +
(\Lambda^{-\beta} u) \Lambda_{-}^{\alpha+\beta} v $$
for all $u$ and $v$ with $\widehat u, \widehat v \ge 0$. \end{lemma}
\begin{proof}
Since $R^{\alpha} (u,v) \precsim (D^{\alpha} u) v$ by the triangle
inequality, Lemma \ref{HyperbolicLambdaEstimates} implies
$$
\Lambda_{-}^{\alpha}(uv) \precsim (\Lambda_{-}^{\alpha} u)v + u
\Lambda_{-}^{\alpha} v + (D^{\alpha} u) v. $$
To finish the proof, combine this with
$$
u \Lambda_{-}^{\alpha} v \precsim (\Lambda^{\alpha} u) v +
(\Lambda^{-\beta} u) \Lambda_{-}^{\alpha + \beta} v. $$
The latter is proved by considering two cases: $\hypwt{\lambda}{\eta} \le
\abs{\xi}$ and $\hypwt{\lambda}{\eta} > \abs{\xi}$, where $(\tau,\xi)$ and
$(\lambda,\eta)$ are the frequencies of $u$ and $v$ respectively.
\end{proof}
\begin{lemma}\label{NegativePowerLambdaMinusEstimate} Let $\alpha,\beta
\ge 0$. Then
{ \renewcommand{\theenumi}{\alph{enumi}}
\renewcommand{\labelenumi}{(\theenumi)}
\begin{enumerate}
\item
$\Lambda_{-}^{-\beta}(uv) \precsim \Lambda_{-}^{-\alpha -
\beta}(\Lambda_{+}^{\alpha} u \Lambda_{+}^{\alpha} v)$, %
\item
$\Lambda_{-}^{-\beta}(uv) \precsim
\Lambda_{-}^{-\alpha-\beta}(u \Lambda^{\alpha} v) + u \Lambda^{-\beta} v$,
\end{enumerate}
}
\noindent
for all $u$ and $v$ with $\widehat u, \widehat v \ge 0$. \end{lemma}
\begin{proof}
Part (a) follows from Lemma \ref{TrivialLambdaMinusEstimate}. Part (b) is
proved by considering two cases: $\hypwt{\tau+\lambda}{\xi+\eta} \le
\abs{\eta}$ and $\hypwt{\tau+\lambda}{\xi+\eta}
>\abs{\eta}$, where $(\tau,\xi)$ and $(\lambda,\eta)$ are the frequencies
of
$u$ and $v$ respectively.
\end{proof}
\paragraph{Proof of \eqref{YME11}.}
Set $\delta = \theta + \varepsilon - 1/2$. In view of the hypotheses of
Theorem \ref{YMTheorem},
$$
\delta \le \half \min \left( \half, s - \frac{n-2}{2} \right). $$
By Lemma \ref{TrivialLambdaMinusEstimate}, it suffices to prove $$
H^{s-\half-\delta,\theta} \cdot H^{s-\half-\delta,\theta} \hookrightarrow
H^{s-\frac{3}{2},0}.
$$
This estimate holds by Theorem \ref{YMProductTheorem}. %
\paragraph{Proof of \eqref{YME12}.}
Set $\zeta = \theta + 2\varepsilon - 1/2$. In view of the hypotheses of
Theorem \ref{YMTheorem},
$$
\zeta \le \half \min \left( \half, s - \frac{n-2}{2} \right). $$
By Lemma \ref{NegativePowerLambdaMinusEstimate}(a), it suffices to prove
\begin{align*}
H^{s-\half-\zeta,\theta-\half} \cdot H^{s-\half-\zeta,\theta}
&\hookrightarrow H^{s-\frac{3}{2},-\half-\varepsilon}. \\
\intertext{Since $\theta > \half$, this estimate is weaker than}
H^{s-\half-\zeta,0} \cdot H^{s-\half-\zeta,\theta} &\hookrightarrow
H^{s-\frac{3}{2},-\half-\varepsilon}, \\
\intertext{which by duality is equivalent to}
H^{\frac{3}{2}-s,\half+\varepsilon} \cdot H^{s-\half-\zeta,\theta}
&\hookrightarrow H^{-s+\half+\zeta,0}.
\end{align*}
The latter holds by Theorem \ref{YMProductTheorem}. %
\paragraph{Proof of \eqref{YMF11}.}
Again we let $\delta = \theta + \varepsilon - 1/2$. By Lemma
\ref{EllipticLambdaEstimates}, it suffices to prove: \begin{align}
\label{YMF111}
\scrH^{1,\theta} \cdot \scrH^{s-\half,\theta} &\hookrightarrow
H^{0,\theta+\varepsilon-\half}, \\
\label{YMF112}
\scrH^{s + \half,\theta} \cdot \scrH^{0,\theta} &\hookrightarrow
H^{0,\theta+\varepsilon-\half}. \\
\intertext{By Lemma \ref{TrivialLambdaMinusEstimate}, \eqref{YMF111}
reduces to}
\notag
H^{1-\delta,\theta} \cdot H^{s-\half-\delta,\theta} &\hookrightarrow L^2,
\\
\intertext{which holds by Theorem \ref{YMProductTheorem}. By Lemma
\ref{YetAnotherLambdaMinusEstimate}, \eqref{YMF112} reduces to}
\notag
H^{s + \half - \delta,\theta} \cdot H^{0,\theta} &\hookrightarrow L^2,
\\
\notag
H^{s + \half + \theta - \delta,\theta} \cdot L^2 &\hookrightarrow L^2.
\end{align}
The former holds by Theorem \ref{YMProductTheorem}, the latter by the
embedding \eqref{InftyEmbedding}. %
\paragraph{Proof of \eqref{YMF13}.}
Set $\zeta = \theta + 2\varepsilon - 1/2$. By Lemma
\ref{EllipticLambdaEstimates}, it suffices to prove: \begin{align}
\label{YMF131}
\scrH^{1,\theta} \cdot \scrH^{s-\half,\theta-\half} &\hookrightarrow
H^{0,\theta+\varepsilon-1}, \\
\label{YMF132}
\scrH^{s + \half,\theta} \cdot \scrH^{0,\theta-\half} &\hookrightarrow
H^{0,\theta+\varepsilon-1}. \\
\intertext{By Lemma \ref{NegativePowerLambdaMinusEstimate}(a),
\eqref{YMF131} reduces to}
\notag
H^{1-\zeta,\theta} \cdot H^{s-\half-\zeta,0} &\hookrightarrow
H^{0,-\half-\varepsilon} \\
\intertext{which by duality is equivalent to}
\notag
H^{1-\zeta,\theta} \cdot H^{0,\half+\varepsilon} &\hookrightarrow
H^{-s+\half+\zeta,0}.
\\
\intertext{This estimate holds by Theorem \ref{YMProductTheorem}. By Lemma
\ref{NegativePowerLambdaMinusEstimate}(b), \eqref{YMF132} reduces to}
\notag
H^{s + \half - \zeta,\theta} \cdot L^2
&\hookrightarrow H^{0,-\half-\varepsilon}, \\
\notag
H^{s + \frac{3}{2}-\theta-\varepsilon,\theta} \cdot L^2 &\hookrightarrow
L^2.
\end{align}
The former holds by Theorem \ref{YMProductTheorem}, the latter by the
embedding \eqref{InftyEmbedding}. %
\paragraph{Proof of \eqref{YME21}.}
This reduces to
\begin{align*}
H^{s-\half,\theta} \cdot H^{s-\half,\theta} &\hookrightarrow
D^{\frac{3}{2} - \gamma-\varepsilon} \Mixed{1}{2n}, \\
\intertext{and in view of \eqref{ModifiedMixednormBound} it suffices to
prove}
H^{s-\half,\theta} \cdot H^{s-\half,\theta} &\hookrightarrow
D^{\frac{3}{2} - \gamma-\varepsilon} L_t^1(L_x^{2n}). \end{align*}
The last estimate holds by Theorem E (since $n \ge 4$). %
\paragraph{Proof of \eqref{YME22}.}
Since $n \ge 4$, we may apply Proposition \ref{DualKlMaEmbedding}. Thus,
it suffices to prove
$$
\scrH^{s-\half,\theta-\half} \cdot \scrH^{s-\half,\theta} \hookrightarrow
\Lambda^{\frac{3}{2}-\frac{n}{2}-\gamma} \Lambda_+
\Lambda_-^{-2\varepsilon} \Mixed{1}{2}.
$$
Replace $\Lambda_+ \Lambda_-^{-2\varepsilon}$ on the right hand side by
$\Lambda^{1-2\varepsilon}$ and apply Lemma \ref{EllipticLambdaEstimates},
thereby reducing to (since $s > \frac{n-2}{2} + \gamma + 2\varepsilon$)
\begin{equation}\label{YME221}
H^{\frac{n-3}{2},0} \cdot
H^{\frac{n-3}{2},\theta}
\hookrightarrow \Lambda^{\frac{5}{2}-\frac{n}{2}} \Mixed{1}{2}.
\end{equation}
We consider the cases $n = 4$ and $n \ge 5$ separately.

If $n = 4$, \eqref{YME221} is just \eqref{YMC22four}, and in view of
\eqref{Trick} and \eqref{VariationEllipticLambdaEstimate} (with $\alpha =
\half$), it suffices to prove: \begin{align}
\label{YME2211}
H^{\half,0} \cdot H^{1,\theta}
&\hookrightarrow L_t^1(L_x^{2}),
\\
\label{YME2212}
L^{2} \cdot H^{1,\theta}
&\hookrightarrow \Lambda^{\half} L_t^1(L_x^{2}). \end{align}
Here we also used \eqref{ModifiedMixednormBound}. By H\"older's
inequality, \eqref{YME2211} reduces to \begin{align*}
H^{\half,0} &\hookrightarrow L_t^2(L_x^{\frac{8}{3}}), \\
H^{1,\theta} &\hookrightarrow L_t^2(L_x^{8}). \end{align*}
These are just \eqref{SpecialEmbeddingF} and \eqref{SpecialEmbeddingB} in
dimension $n = 4$. By Sobolev embedding we reduce \eqref{YME2212} to $$
L^{2} \cdot H^{1,\theta}
\hookrightarrow L_t^1(L_x^{\frac{8}{5}}). $$
But this holds by \eqref{SpecialEmbeddingB}.

Now assume $n \ge 5$. By Lemma \ref{EllipticLambdaEstimates},
\eqref{YMF221} reduces to
\begin{align}
\label{YMF2213}
H^{1,0} \cdot H^{\frac{n-3}{2},\theta}
&\hookrightarrow L_t^1(L_x^{2}),
\\
\notag
H^{\frac{n-3}{2},0} \cdot H^{1,\theta}
&\hookrightarrow L_t^1(L_x^{2}).
\end{align}
Using H\"older's inequality, these follow from \eqref{SpecialEmbeddingA},
\eqref{SpecialEmbeddingG}, \eqref{SpecialEmbeddingF} and
\eqref{SpecialEmbeddingB}.
\paragraph{Proof of \eqref{YMF12}.}
By Lemma \ref{EllipticLambdaEstimates}, it suffices to prove \begin{align}
\label{YMF121}
\Lambda^{s-1} \Lambda_{-}^{\half} \X^{s} \cdot \scrH^{s-\half,\theta}
&\hookrightarrow H^{0,\theta+\varepsilon-1}, \\
\label{YMF122}
\Lambda^{-\half} \Lambda_{-}^{\half} \X^{s} \cdot \scrH^{0,\theta}
&\hookrightarrow H^{0,\theta+\varepsilon-1}. \end{align}

Set $\zeta = \theta + 2\varepsilon - 1/2$. By Lemma
\ref{NegativePowerLambdaMinusEstimate}(a), \eqref{YMF121} reduces to
\begin{align*}
H^{1-\zeta,0} \cdot H^{s-\half-\zeta,\theta} &\hookrightarrow
H^{0,-\half-\varepsilon}, \\
\intertext{which by duality is equivalent to}
H^{0,\half+\varepsilon} \cdot H^{s-\half-\zeta,\theta} &\hookrightarrow
H^{-1+\zeta,0}.
\end{align*}
The latter holds by Theorem \ref{YMProductTheorem}.

By Lemma \ref{NegativePowerLambdaMinusEstimate}(b), \eqref{YMF122} reduces
to two estimates:
\begin{align}
\label{YMF1221}
\Lambda^{\theta + \varepsilon-\frac{3}{2}} \Lambda_{-}^{\half} \X^{s}
\cdot H^{0,\theta}
&\hookrightarrow L^{2},
\\
\label{YMF1222}
\Lambda^{\zeta-\half} \Lambda_{-}^{\half} \X^{s} \cdot H^{0,\theta}
&\hookrightarrow H^{0,-\half-\varepsilon}. \\
\intertext{(Recall that $\zeta = \theta + 2\varepsilon - 1/2$.) In view of
\eqref{SpaceInftyEmbedding} and \eqref{EnergyEmbedding},}
\notag
H^{s+\frac{3}{2}-\theta-\varepsilon,0} \cdot H^{0,\theta} &\hookrightarrow
L^{2},
\end{align}
which implies \eqref{YMF1221}. Since $\zeta + \varepsilon = \gamma$,
\eqref{YMF1222} follows from \eqref{SpecialEmbeddingL}. %
\paragraph{Proof of \eqref{YMF21}.}
This reduces to
$$
H^{s+\half,\theta} \cdot H^{s-\half,\theta} \hookrightarrow
\Lambda^{\half-\gamma-\varepsilon} \Mixed{1}{2n}. $$
In view of \eqref{Trick} and \eqref{VariationEllipticLambdaEstimate} (with
$\alpha = \half - \gamma - \varepsilon$), this reduces to \begin{align*}
H^{s+\half,\theta} \cdot
H^{s-\gamma-\varepsilon,\theta}
&\hookrightarrow
\mixed{1}{2n},
\\
H^{s+\gamma+\varepsilon,\theta} \cdot
H^{s-\gamma-\varepsilon,\theta}
&\hookrightarrow
\Lambda^{\half-\gamma-\varepsilon} \mixed{1}{2n}. \end{align*}
The former holds by \eqref{SpecialEmbeddingD} and
\eqref{SpecialEmbeddingE}, the latter holds by Theorem E. %
\paragraph{Proof of \eqref{YMF22}.}
By Lemma \ref{NegativePowerLambdaMinusEstimate}(b), this reduces to
\begin{align}
\label{YMF221}
\Lambda^{-\half} \Lambda_{-}^{\half} \X^s \cdot
H^{s-\half-2\varepsilon,\theta}
&\hookrightarrow \Lambda^{-\half-\gamma} \Lambda_+
\Lambda_-^{\half+\varepsilon} \Mixed{1}{2n}, \\
\label{YMF222}
\Lambda^{-\half} \Lambda_{-}^{\half} \X^s \cdot H^{s-\varepsilon,\theta}
&\hookrightarrow \Lambda^{-\half-\gamma} \Lambda_+ \Mixed{1}{2n}. \\
\intertext{By Proposition \ref{DualKlMaEmbedding}, \eqref{YMF221} reduces
to}
\notag
\Lambda^{-\half} \Lambda_{-}^{\half} \X^s \cdot
H^{s-\half-2\varepsilon,\theta}
&\hookrightarrow \Lambda^{\frac{3-n}{2}-\gamma} \Mixed{1}{2}.
\\
\intertext{By Lemma \ref{EllipticLambdaEstimates}, the latter reduces to
two estimates:}
\label{YMF2211}
\Lambda^{\frac{n-4}{2} + \gamma} \Lambda_{-}^{\half} \X^s \cdot
H^{\frac{n-3}{2},\theta}
&\hookrightarrow \Mixed{1}{2}.
\\
\label{YMF2212}
\Lambda^{-\half} \Lambda_{-}^{\half} \X^s \cdot H^{0,\theta}
&\hookrightarrow \Mixed{1}{2}.
\\
\intertext{If $n = 4$, then \eqref{YMF2211} follows from}
\notag
\Mixed{1}{8} \cdot H^{\half,\theta}
&\hookrightarrow \Mixed{1}{2},
\end{align}
which holds by Proposition \ref{MixednormHolder}, Lemma \ref{DensityLemma}
and \eqref{SpecialEmbeddingH}. If $n \ge 5$, then \eqref{YMF2211} reduces
to \eqref{YMF2213}.

In view of \eqref{EnergyEmbedding}, \eqref{YMF2212} reduces to $$
\Lambda^{-\half-\gamma} \Mixed{1}{2n} \cdot \mixed{\infty}{2}
\hookrightarrow \Mixed{1}{2}, $$
which holds by Proposition \ref{MixednormHolder}, Lemma \ref{DensityLemma}
and \eqref{SpecialEmbeddingJ} (since $\gamma > 0$).

By Lemma \ref{EllipticLambdaEstimates}, \eqref{YMF222} follows from
\begin{align*}
H^{\frac{n-1}{2},0} \cdot H^{\frac{n-2}{2},\theta} &\hookrightarrow
D^{\half} \mixed{1}{2n}. \\
\intertext{By Sobolev embedding, $\mixed{1}{n} \hookrightarrow D^{\half}
\mixed{1}{2n}$, so it suffices to have}
H^{\frac{n-1}{2},0} \cdot H^{\frac{n-2}{2},\theta} &\hookrightarrow
\mixed{1}{n}.
\end{align*}
This follows from \eqref{SpecialEmbeddingC} and \eqref{SpecialEmbeddingE}.
\paragraph{Proof of \eqref{YMF23}.}
By Proposition \ref{DualKlMaEmbedding}, this reduces to \begin{align*}
\Lambda^{-\half} \X^s \cdot \scrH^{s-\half,\theta-\half} &\hookrightarrow
\Lambda^{\frac{3-n}{2}-\gamma-2\varepsilon} \Mixed{1}{2}.
\\
\intertext{By Lemma \ref{EllipticLambdaEstimates}, it suffices to prove}
H^{1,\theta} \cdot H^{\frac{n-3}{2},0}
&\hookrightarrow \mixed{1}{2}.
\\
H^{\frac{n-1}{2} + \varepsilon,\theta} \cdot L^{2} &\hookrightarrow
\mixed{1}{2}.
\end{align*}
The first of these holds by \eqref{SpecialEmbeddingB} and
\eqref{SpecialEmbeddingF}. The second follows from
\eqref{SpecialEmbeddingD}.
\section{Main Estimates for \eqref{CFWMType}}\label{WMMProof}
Here we prove the following.
\begin{theorem}\label{CFWMTheorem}
Let $n \ge 3$, $s > \frac{n-2}{2}$. Assume that $\theta, \varepsilon$ and
$q$ satisfy
\begin{gather*}
\half < \theta < \min \left( 1, \half + \half \left[ s-\frac{n-2}{2}
\right] \right),
\\
0 \le \varepsilon < \min \left( 1 - \theta, \theta - \half, \half \left[
s-\frac{n-2}{2} + 1 - 2\theta \right] \right), \\
\half \le \frac{1}{q} < \frac{3}{2} - \theta - \varepsilon. \end{gather*}
Let $\X^s$ be the Banach space given by the norm $$
\norm{u} = \Spacetimenorm{u}{s}{\theta} + \Mixednorm{\Lambda^{-1}
\Lambda_- u}{q}{\infty},
$$
then
\begin{equation}\label{CFWMA}
\widetilde Q \bigl( \X^s, \X^s \bigr) \hookrightarrow \Lambda_+
\Lambda_-^{1-\varepsilon} \X^s,
\end{equation}
where $\widetilde Q$ is the null form appearing in \eqref{CFWMType}.
\end{theorem}
This implies part (c) of the Main Theorem, in view of Theorem
\ref{SpecializedWellPosednessTheorem}, since conditions (a)--(e) of the
latter are satisfied by the space $\X^{s}$ (condition (a) holds by
Proposition \ref{SecondEnergyEmbeddingProposition}; condition (b) is
obviously satisfied; conditions (c) and (d) follow from Theorem
\ref{BasicConditionsTheorem}, in view of Proposition
\ref{BasicSpacesProposition}; finally, condition (e) holds by Proposition
\ref{TimeCutOffProposition}).

By the definition of $\X^{s}$, \eqref{CFWMA} is equivalent to two
estimates: \begin{align}
\label{CFWM1}
\Lambda^{-1} \widetilde Q(\X^s,\X^s) &\hookrightarrow
H^{s,\theta+\varepsilon-1}, \\
\label{CFWM2}
\Lambda^{-1} \widetilde Q(\X^s,\X^s) &\hookrightarrow \Lambda_+
\Lambda_-^{-\varepsilon} \Mixed{q}{\infty}.
\end{align}
The latter can be proved without using the null structure of $\widetilde
Q$. In fact, we shall rely on the following crude estimate: $$
\Lambda^{-1} \widetilde Q(\phi,\psi)
\precsim D^{-1} \Lambda_+ \phi' D^{-1} \Lambda_+ \psi' \quad
\text{whenever} \quad \text{$\phi \preceq \phi'$, $\psi \preceq \psi'$}. $$
(The trivial proof of this is omitted.) Thus, \eqref{CFWM2} reduces to
\begin{align*}
D^{-1} H^{s-1,\theta} \cdot D^{-1} H^{s-1,\theta} &\hookrightarrow
\Lambda^{1-\varepsilon} \Mixed{q}{\infty}. \\
\intertext{By Sobolev embedding, this reduces to}
D^{-1} H^{s-1,\theta} \cdot D^{-1} H^{s-1,\theta} &\hookrightarrow
\Lambda^{1-\frac{n}{r}-2\varepsilon} \Mixed{q}{r} \end{align*}
for any $2 \le r < \infty$. The latter holds by Theorem E, if we take $r$
so large (to ensure that $(2q,2r)$ is wave admissible) that $$
\frac{1}{q} \le \frac{n-1}{2} - \frac{n-1}{2r}. $$
(We can do this since $q > 1$ and $n \ge 3$.)

To prove \eqref{CFWM1} we need to take into account the null structure. In
fact, proving \eqref{CFWM1} can be reduced to proving four estimates:
\begin{align}
\label{CFWM3}
D^{-1} H^{M,\theta} \cdot H^{0,\theta} &\hookrightarrow
H^{0,\theta+\varepsilon-1},
\\
\label{CFWM4}
\scrH^{1,\theta-1} \cdot H^{s,\theta} &\hookrightarrow
H^{0,\theta+\varepsilon-1},
\\
\label{CFWM5}
\Mixed{q}{\infty} \cdot H^{0,\theta} &\hookrightarrow
H^{0,\theta+\varepsilon-1},
\\
\label{CFWM6}
R( H^{s+1,\theta} , H^{0,\theta} ) &\hookrightarrow
H^{0,\theta+\varepsilon-1}.
\end{align}
Here $M > 0$ can be taken arbitrarily large (\eqref{CFWM3} is a low
frequency estimate which comes up because we want to replace $D^{-1}$ by
$\Lambda^{-1}$ in certain places). The important estimates are
\eqref{CFWM4}--\eqref{CFWM6}.

We will need the following theorem (essentially a corollary of Theorem F;
see the Appendix for the proof). %
\begin{theorem}\label{CFWMProductTheorem} Let $n \ge 3$ and $\theta >
\half$. Then $$
H^{a,\theta} \cdot H^{b,\theta} \hookrightarrow H^{-c,\theta} $$
holds for all $a,b,c$ satisfying
\begin{align*}
a,b,c &\ge 0,
\\
c &< \frac{n-1}{2},
\\
a + b + c &\ge \frac{n-1}{2} + \theta.
\end{align*}
\end{theorem}
The basic estimate for the null form $\widetilde Q$ is as follows. %
\begin{lemma}\label{QtildeEstimate}
The estimate
$$
D^{-1} \widetilde Q(\phi,\psi) \precsim
D^{-1} D_- \phi' \cdot \psi' + R( D^{-1} \phi', \psi') + \text{symmetric
terms}
$$
holds whenever $\phi \preceq \phi'$ and $\psi \preceq \psi'$. \end{lemma}
The proof can be found in \cite[Lemmas 2.3 and 2.4]{Kl-Ma3}. We shall also
make use of the following:
\begin{lemma}\label{CFWMNullFormEstimate} Let $\gamma, M \ge 0$. Then
\begin{multline*}
\Lambda^{\gamma-1} \widetilde Q (\phi,\psi) \precsim \Lambda^{-M} D^{-1}
\Lambda_+ \phi' \cdot \Lambda^\gamma \psi' + \Lambda^{\gamma-1} \Lambda_-
\phi' \cdot \psi' \\
+ \Lambda^{-1} \Lambda_- \phi' \cdot \Lambda^\gamma \psi' + R(
\Lambda^{-1} \phi', \Lambda^\gamma \psi') + \text{symmetric terms}
\end{multline*}
whenever $\phi \preceq \phi'$ and $\psi \preceq \psi'$. \end{lemma}
\begin{proof}
By Lemma \ref{QtildeEstimate},
$$
\Lambda^{-1} \widetilde Q (\phi,\psi) \precsim D^{-1} \Lambda_- \phi'
\cdot \psi' + R( D^{-1} \phi', \psi') + \text{symmetric terms}.
$$
Since $D^{-1} u \lesssim \Lambda^{-1} u + D^{-1} \Lambda^{-M} u$ whenever
$\widehat u \ge 0$, we conclude that \begin{multline*}
\Lambda^{-1} \widetilde Q (\phi,\psi) \precsim \Lambda^{-M} D^{-1}
\Lambda_+ \phi' \cdot \psi' + R(\Lambda^{-M} D^{-1} \phi', \psi')
\\
+ \Lambda^{-1} \Lambda_- \phi' \cdot \psi' + R( \Lambda^{-1} \phi', \psi')
+ \text{symmetric terms}. \end{multline*}
Since $R(u,v) \precsim (Du) v$ whenever $\widehat u, \widehat v \ge 0$,
the second term on the right hand side is subsumed in the first. It is
easy to see that
$$
\Lambda^\gamma R( \Lambda^{-1} u, v) \precsim R( \Lambda^\gamma u,
\Lambda^{-1} v ) + R( \Lambda^{-1} u, \Lambda^\gamma v )
$$
provided $\widehat u, \widehat v \ge 0$. Combining this with Lemma
\ref{EllipticLambdaEstimates} yields the desired estimate. \end{proof}
We shall also need the following estimate for the operator $R$. %
\begin{lemma}\label{CFWMEstimateForR}
Let $\alpha \in [0,1]$, $\delta \ge 0$. Then $$
R (\phi,\psi) \precsim
\Lambda_-^{1-\alpha} R^{\alpha}(\phi',\psi') + \Lambda_- \phi' \cdot \psi'
+ \Lambda_{-}^{-\delta} ( \Lambda^{\delta} \phi' \cdot \Lambda_{-} \psi')
$$
whenever $\phi \preceq \phi'$ and $\psi \preceq \psi'$. \end{lemma}
\begin{proof}
It is readily verified that the symbol $r$ of $R$ satisfies $$
r(\tau,\xi;\lambda,\eta) \le A + B + C,
$$
where $A = \hypwt{\tau+\lambda}{\xi+\eta}$, $B = \hypwt{\tau}{\xi}$ and $C
= \hypwt{\lambda}{\eta}$. We consider three cases, corresponding to $A$,
$B$ or $C$ being the maximum of the three.

If $A$ is the maximum, then $r \lesssim A^{1-\alpha} r^{\alpha}$.

If $B$ is the maximum, then $r \lesssim B$.

If $C$ is the maximum, we consider two subcases: (i) $\abs{\xi} \ge A$;
and (ii) $\abs{\xi} < A$. In case (i), $r \lesssim A^{-\delta} A^{\delta}
C \lesssim A^{-\delta} \abs{\xi}^{\delta} C$. In case (ii), $r \lesssim
\abs{\xi}^{1-\alpha} r^{\alpha} \lesssim A^{1-\alpha} r^{\alpha}$.
\end{proof}
Applying Lemma \ref{CFWMNullFormEstimate} to \eqref{CFWM1} (with $\gamma =
s$), we see that \eqref{CFWM1} reduces to \eqref{CFWM3}--\eqref{CFWM6}. %
\paragraph{Proof of \eqref{CFWM3}.}
This is weaker than
$$
D^{-1} H^{M,0} \cdot H^{0,\theta} \hookrightarrow L^{2}.
$$
By H\"older's inequality and \eqref{EnergyEmbedding}, the latter reduces to
$$
D^{-1} H^{M,0} \hookrightarrow \mixed{2}{\infty}, $$
which holds by Sobolev embedding for $M > \frac{n-2}{2}$, since $n \ge 3$.
\paragraph{Proof of \eqref{CFWM4}.}
By Lemma \ref{NegativePowerLambdaMinusEstimate}(b), this reduces to
\begin{align}
\label{CFWM41}
\scrH^{2-\theta-\varepsilon,\theta-1} \cdot H^{s,\theta} &\hookrightarrow
L^2,
\\
\label{CFWM42}
H^{1-\delta,\theta-1} \cdot H^{s,\theta} &\hookrightarrow
H^{0,-\half-\varepsilon},
\\
\intertext{where $\delta = \theta + \varepsilon - \half$. Observe that
\eqref{CFWM41} is weaker than}
\notag
H^{1-\theta-\varepsilon,\theta} \cdot H^{s,\theta} &\hookrightarrow L^2,
\\
\intertext{which holds by Theorem F. \eqref{CFWM42} is equivalent to}
\notag
H^{0,\half+\varepsilon} \cdot H^{s,\theta} &\hookrightarrow
H^{\delta-1,1-\theta},
\end{align}
and the latter holds by Theorem \ref{CFWMProductTheorem}. %
\paragraph{Proof of \eqref{CFWM5}.}
This holds by \eqref{SpecialEmbeddingK}. %
\paragraph{Proof of \eqref{CFWM6}.}
By Lemma \ref{CFWMEstimateForR} (with $\alpha = \theta+\varepsilon$ and
$\delta = \theta + \varepsilon - \half$), this reduces to three estimates,
one of which is \eqref{CFWM5}; the other two are \begin{align*}
R^{\theta + \varepsilon}( H^{s+1,\theta} , H^{0,\theta} ) &\hookrightarrow
L^{2},
\\
H^{s+1-\delta,\theta} \cdot H^{0,\theta-1} &\hookrightarrow
H^{0,-\half-\varepsilon}. \\
\intertext{The former holds by Theorem F, and the latter is equivalent to}
H^{s+1-\delta,\theta} \cdot H^{0,\half+\varepsilon} &\hookrightarrow
H^{0,1-\theta},
\end{align*}
which holds by Theorem \ref{CFWMProductTheorem}.
\section{Further Results, Open Problems and Historical
Remarks}\label{FurtherResults}
(1)\,\,\,
The results discussed in the Main Theorem
(see section \ref{BasicEquations}) confirm part (i)
of the General WP Conjecture (see section \ref{Motivation}) for the 
equations  we consider, with two notable exceptions\footnote{Strictly
 speaking the results we refer to, in connection
to gauge theories,  concern only the model problems \eqref{MKGType} and
\eqref{YMType}. There is  however little doubt that the results
can be extended to the full (MKG) and \eqref{YM} systems. See 
\cite{Kl-Ma0.2},\cite{Kl-Ma0.3} as well as  \cite{Cu} and \cite{Ke-Ta3}
for a full treatment of the equations.}. 
The first concerns (MKG) and \eqref{YM}, as well as the simplified
model problems \eqref{MKGType} and \eqref{YMType}, in dimension $n=3$.
In both cases the critical WP exponent is $s_c = \half$, and one can prove
local well-posedness for $s > s_c + \frac{1}{4}$ using the $H^{s,\theta}$ spaces
with $\theta = s$; see \cite{Cu} and \cite{Ke-Ta3}. It is easy to
see that while the first iterate, in all the above-mentioned cases, is
well-posed for $s > \half$, it fails to belong to the corresponding $H^{s,\theta}$
space for $\half < s\le\frac{3}{4}$, $\theta > \half$. One can
show that any strategy based on norms which depend only on the size of the
Fourier transform, such as those used in this survey, is bound to fail to
prove well-posedness for $s$ in the range $\half < s \le \frac{3}{4}$.
Is it possible that well-posedness fails in that range?

The second case in which our present techniques do
not allow us to go all the way to the critical exponent
is the \eqref{CFWMType} equation in dimension $n=2$. Even though the system
\eqref{WaveMapsModelA} is a equivalent, through a simple transformation,
to the standard wave maps equation \eqref{WMType}, for which we can prove 
well-posedness for all $s > s_c(\mathrm{WM}) = 1$ ($n = 2$), we cannot treat it 
for $s$ close to the corresponding critical exponent $s_{c} 
(\mathrm{WMM}) = 0$ ($n = 2$). In fact, the best result proved so far 
for \eqref{CFWMType} in dimension $n = 2$, is that it is
well-posed for $s > \frac{1}{4}$; see \cite{Se}. In comparison with the 
higher-dimensional ($n \ge 3$) case, discussed in section 
\ref{WMMProof}, we remark that the estimate \eqref{CFWM1} is true for 
all $s > 0$ in dimension $n = 2$, whereas the ``$\mixed{1}{\infty}$'' 
estimate \eqref{CFWM2} fails for $0 < s < \frac{1}{8}$, by a 
counterexample given in \cite{Se}.
For a more in-depth discussion of the \eqref{CFWMType} equation, see 
\cite{Se, Se3}.

Another interesting class of semilinear wave equations,
which was not discussed in this survey, is provided by
Klein-Gordon-Dirac and Maxwell-Dirac. 
Improved  well posedness results, based on bilinear estimates, were proved  
in \cite{Bour1} and \cite{Bour2}, see also \cite{Kl-Ma0.1a}; the question of
optimal  well-posedness remains however wide open.

\medskip
\noindent
(2)\,\,\, Close to nothing is known concerning part (ii) of the General WP
Conjecture, except for semilinear scalar wave equations of the type
$\square\phi+V'(\phi)=0$; see \cite{Ka}, \cite{Sh-St2}. An important
advance was made recently by D. Tataru \cite{Ta2, Ta3} who was able
to prove, in the case of \eqref{WMType} equations, global well-posedness
for data in
the Besov space $B^{2, 1}_{\frac{n}{2}}$, in any dimension $n \ge 2$. It
would be interesting to extend Tataru's result to the other cases covered
by the Main Theorem. We expect that all classical field theories are
globally well-posed (in the strong sense of our Main Theorem) for small
data in $B^{2, 1}_{s_c}$ with $s_c$ the
critical WP exponent. The fundamental problem of well-posedness in
$H^{s_c}$ is far more difficult (see the relevant discussion on well-posedness
and its connection to the issue of global
regularity in \cite{Kl}). Ultimately the issue of optimal well-posedness
must be tied to that of global regularity for all finite energy data, in
the case of critical nonlinearities, or that of spontaneous formation of
singularities\footnote{See \cite{St1} for an up to date survey concerning
weak solutions and formation of singularities, as well as known results in
the case of equivariant or spherically symmetric wave maps.}
in the case of supercritical equations.

\medskip
\noindent
(3)\,\,\, Recently M. Keel and T. Tao, see \cite{Ke-Ta3}, were able to prove
global existence for the full (MKG) system for arbitrarily large $H^s$ initial
data with
$\frac{3}{4}<s<1$. Local well-posedness in the same range was dealt with
in \cite{Cu}.
Global $H^1$ well-posedness for the harder case of the (YM) equations,
corresponding to the energy norm, was treated in \cite{Kl-Ma0.3}.

\medskip
\noindent
(4)\,\,\, The issue of optimal well-posedness for quasilinear wave
equations has only very recently started to be investigated. We refer the
interested reader to the works of Chemin-Bahouri
\cite{ChBa1, ChBa2}, Tataru \cite{Ta4, Ta5} and Klainerman \cite{Kl2}. See
also the relevant discussion in \cite{Kl}. It is not difficult to predict
that this  very important area of activity  will play a
 predominant  role in the  future.

\medskip
\noindent
(5)\,\,\, The first improved\footnote{ By comparison to what can be derived
by Strichartz-type estimates. } space-time regularity results for null
quadratic forms appear in \cite{Kl-Ma0.1}. Those estimates were used, by
virtue of Duhamel's principle, to set up an iteration procedure with
respect to a space $\X^s_T$ (see the discussion
in section \ref{TheIterationSpace}) defined by the space-time norm
\begin{equation}\label{neweq1}
  \norm{u}_{L_t^\infty([0,T] , H^s)} + \norm{\partial_t
  u}_{L_t^\infty([0,T], H^{s-1})}
  + T^{\frac{1}{2}}\norm{D^{s-1}\square u}_{L^2([0,T]\times\R^n)}
\end{equation}
and derive improved well-posedness results for a general class of
nonlinear wave equations verifying the null condition (see \cite{Kl0}),
including \eqref{WMType}. The same type of estimates and a similar version of
the $\X^s_T$ iteration space were used
in \cite{Kl-Ma0.2, Kl-Ma0.3} to derive global well-posedness
results in the energy norm for the full (MKG) and (YM) systems. Observe
that the norm \eqref{neweq1} is essentialy the same as that of the spaces
$\scrH^{s,\theta}$ for $\theta=1$. The  case of the Yang-Mills-Higgs
equations, with critical power for the scalar Higgs component,
was treated  using a clever localization  of  the same norm  in \cite{Ke}.
Variations of the same techniques were  also used, see \cite{Bour1,Bour2},
to derive nontrivial results for the Maxwell-Dirac and Klein-Gordon-Dirac equations. 

It is noteworthy that in \cite{Kl-Ma0.1} there appear also sharper bilinear estimates
corresponding to the homogeneous $H^{s,\theta }$, $\theta =\frac{1}{2}$
spaces (see section \ref{WaveSobolevSpaces}).
These better estimates could not, however, be used in
an iterative procedure; $\theta=\frac{1}{2}$ leads to an obvious
logarithmic divergence. The use of the
$H^{s,\theta}$ spaces, for $\theta > \frac{1}{2}$, was initiated in
\cite{Kl-Ma1} under the influence of the works of Bourgain
\cite{B} and Kenig-Ponce-Vega \cite{KPV} for dispersive equations. The
new idea, provided by these works, was to introduce a time cut-off
function which allows one to replace $\square$ by
$\Lambda_{+}\Lambda_{-}$. The inhomogeneous $H^{s,\theta}$ for $s>s_c$,
$\theta>\frac{1}{2}$, avoids the above logarithmic divergences and allowed
one to prove a well-posedness result for $s>s_c$, in the case
of the \eqref{WMType} system, for $n\ge 3$. The case $n=2$ was treated
later in
\cite{Kl-Se} with the help of the new bilinear estimates proved in
\cite{Kl-Ma2}. The \eqref{WMType} system is the only one for which the
$\scrH^{s,\theta}$ spaces alone suffice to prove optimal WP results for
$s>s_c$. More precisely, for (MKG), (YM) $n \ge 4$ and
\eqref{CFWMType} $n \ge 3$, some of the product properties of the
$\scrH^{s,\theta}$
spaces, necessary to carry the step by step iteration, fail by a lot. The
starting point in
\cite{Kl-Ma3},\cite{Kl-Ma4} and \cite{Kl-Ma5} was the observation that,
despite this failure, one can check nevertheless that the second iterates
belong
to $\scrH^{s,\theta}$. Moreover the second iterates satisfy additional
trilinear properties which, when taken into account, allow one to prove
inductively that all iterates belong both to $\scrH^{s,\theta}$ and
satisfy the same trilinear conditions. In the wake of the bilinear
estimates of Theorem B, proved in \cite{Kl-Ta}, it
became clear
that the additional trilinear conditions can be more conveniently
rephrased in
terms of the $\Mixed{q}{r}$ spaces discussed in this paper.
See also \cite{Ma} for a review of this circle of ideas.

\medskip
\noindent
(6)\,\,\, The spaces $H^{s,\theta}$ are by no means new in PDE. Before
\cite{Kl-Ma0.1} and the systematic use of such spaces by Bourgain \cite{B}
in the study of optimal well-posedness for
periodic initial conditions for KdV and nonlinear Schr\"odinger
equations (see also \cite{KPV}), such spaces were used in microlocal
analysis, in particular
in the the study of propagation of singularities for nonlinear wave
equations; see \cite{Be}. The novel idea, in both \cite{Kl-Ma0.1} and
\cite{B}, was to estimate directly, in a space-time $L^2$ norm,
the principal quadratic part of the nonlinear term\footnote{In the case of
the KdV equation $u_t+u_{xxx}+uu_x=0$, treated in \cite{B}, this was
$uu_x$.
In \cite{Kl-Ma0.1} one relies on space-time $L^2$ estimates for the null
quadratic forms $Q_0$ and $Q_{\alpha\beta}$.}. These new types of
estimates\footnote{Previous attempts to prove optimal well-posedness
results relied on the idea of treating the nonlinear part of the equation
as a source term and using the best available estimates, such as
Strichartz, for the corresponding linear inhomogeneous equation (see e.g.
\cite{KPV0}, \cite{PS}). In some situations, such as nonlinear wave
equations of the type $\square\phi=\pm\phi^p$,
this procedure is in fact optimal, see \cite{LS}.}, which we now refer to
as bilinear, provide additional regularity information in connection with the issue of
optimal well-posedness. The $L^2$ set-up of the $ H^{s,\theta}$ spaces is
the most convenient\footnote{Indeed, in view of the Plancherel identity
it suffices to estimate bilinear weighted convolutions. This idea
has its origin in  the  proof of the restriction theorem  for the
 special  of the  $L^4$-norm;
see for example \cite{CaS}.   } way to
take into account possible cancellations between the symbol of the special
quadratic part of the nonlinear equation and the symbol of the corresponding
linear operator.
Aside from simplicity there is in fact no reason to stop at $L^2$; as we have
seen above, additional information can be provided by combining the
$H^{s,\theta}$ norm with a suitable $\Mixed{q}{r}$ norm. Further progress in this
respect may be expected from the $\mixed{q}{r}$ bilinear estimates
conjectured in \cite{Fo-Kl} and proved partially in \cite{W} and
\cite{Tao1}. Clever modifications of the $ H^{s,\theta}$ spaces appear
also in \cite{Ta2} and \cite{Ta3}. 
\appendix
\section{Appendix}
\subsection{Counterexamples}
Here we prove the negative statement in Theorem \ref{NaiveYMTheorem}. The
argument below is a slight modification of the counterexample used in
\cite{Kl-Ma2}. We construct, for all sufficiently large $L > 0$, functions
$u_{L}$ and $v_{L}$ such that for any $s$ and $\theta$,
\begin{gather}
\label{RHSnorms}
\Spacetimenorm{u_{L}}{s}{\theta} \sim L^{s + \theta + \frac{n}{2}}, \qquad
\Spacetimenorm{v_{L}}{s}{\theta} \sim L^{2s + \frac{n+1}{2}}, \\
\label{FourierSize}
\Fourier \left\{ Q_{ij}(D^{-1}u_{L},v_{L}) \right\} \sim L^{2}
\Fourier(u_{L} v_{L})
\sim \Fourier \left\{ \Lambda^{\half}( \Lambda^{-\half}
\Lambda_{-}^{\half} u_{L} \cdot \Lambda^{\half} v_{L}) \right\}.
\end{gather}
Moreover, for all $(\tau,\xi)$ in a certain set $C$ with measure $\sim
L^{n+1}$, \begin{equation}\label{FourierSizeB}
\Fourier\left\{ \Lambda^{s-1} \Lambda_{-}^{\theta-1} (u_{L} v_{L})
\right\} (\tau,\xi) \sim L^{2(s-1) + n}. \end{equation}
It follows from \eqref{FourierSize}, \eqref{FourierSizeB} and Plancherel's
theorem that
\begin{multline*}
\spacetimenorm{Q_{ij}(D^{-1}u_{L},v_{L})}{s-1}{\theta-1} \sim
\spacetimenorm{\Lambda^{\half}( \Lambda^{-\half} \Lambda_{-}^{\half} u_{L}
\cdot \Lambda^{\half} v_{L})}{s-1}{\theta-1} \\
\gtrsim L^{2(s-1) + n + 2} \sqrt{\abs{C}} \sim L^{2s + n + \frac{n+1}{2}}.
\end{multline*}
But by \eqref{RHSnorms},
$$
\Spacetimenorm{u_{L}}{s}{\theta} \Spacetimenorm{v_{L}}{s}{\theta} \sim
L^{3s + \theta + n + \half}.
$$
We conclude that the estimates
\begin{align*}
\spacetimenorm{Q_{ij}(D^{-1}u,v)}{s-1}{\theta-1} &\lesssim
\Spacetimenorm{u}{s}{\theta} \Spacetimenorm{v}{s}{\theta}, \\
\spacetimenorm{\Lambda^{\half}( \Lambda^{-\half} \Lambda_{-}^{\half} u
\cdot \Lambda^{\half} v)}{s-1}{\theta-1} &\lesssim
\Spacetimenorm{u}{s}{\theta} \Spacetimenorm{v}{s}{\theta} \end{align*}
must fail if $s < \frac{n}{2} - \theta$.

We may of course take $i = 1 < j$. Let $A$ be the set of $(\lambda,\eta)
\in \R^{1+n}$ such that
$$
\abs{\lambda-\eta_{1}} \le 1, \qquad \frac{L}{2} \le \eta_{1} \le L,
\qquad \frac{L}{2} \le \abs{\eta'} \le L, $$
where we write $\eta = (\eta_{1},\dots,\eta_{n})$ and $\eta' =
(\eta_{2},\dots,\eta_{n})$. With the same notation, let $B$ be the set of
$(\tau,\xi)$ such that
$$
\bigabs{\tau - \abs{\xi}} \le 8, \qquad \frac{L^{2}}{2} \le \xi_{1} \le
4L^{2},
\qquad \abs{\xi'} \le 2L,
$$
and let $C$ be the set determined by
$$
\bigabs{\tau - \abs{\xi}} \le 1, \qquad L^{2} \le \xi_{1} \le 2L^{2},
\qquad \abs{\xi'} \le L.
$$
Let $\widehat{u_{L}}$ and $\widehat{v_{L}}$ be the characteristic
functions of $A$ and $B$ respectively.

Clearly, \eqref{RHSnorms} is satisfied. Also,
\begin{equation}\label{ConvolutionSupports}
(\lambda,\eta) \in A, (\tau,\xi) \in C \implies (\tau-\lambda,\xi-\eta)
\in B,
\end{equation}
since
\begin{align*}
\bigabs{\tau-\lambda-\abs{\xi-\eta}}
&\le \bigabs{\tau-\abs{\xi}} + \abs{\lambda-\eta_{1}} + \abs{\xi}-\xi_{1}
+ \abs{\xi-\eta}-(\xi_{1}-\eta_{1}) \\
&\le 2 + \frac{\abs{\xi'}^{2}}{\abs{\xi}+\xi_{1}} +
\frac{\abs{\xi'-\eta'}^{2}}{\abs{\xi-\eta}+\xi_{1}-\eta_{1}} \le 2 + 1 + 5
= 8.
\end{align*}
Observe that \eqref{ConvolutionSupports} implies $$
\Fourier(u_{L} v_{L}) (\tau,\xi) = \abs{A} \sim L^{n} \quad \text{for all}
\quad (\tau,\xi) \in C, $$
and \eqref{FourierSizeB} follows.

To prove \eqref{FourierSize}, write
\begin{multline*}
\Fourier \left\{ Q_{1j}(D^{-1}u_{L},v_{L}) \right\}(\tau,\xi) \\
= \int_{\R^{1+n}} \left( \frac{\eta_{j}}{\abs{\eta}} (\xi_{1}-\eta_{1})
- \frac{\eta_{1}}{\abs{\eta}} (\xi_{j} - \eta_{j}) \right)
\widehat{u_{L}}(\lambda,\eta) \widehat{v_{L}}(\tau-\lambda,\xi-\eta) \,
d\lambda \, d\eta.
\end{multline*}
Obviously,
$$
(\lambda,\eta) \in A, (\tau,\xi) \in B \implies
\frac{\eta_{j}}{\abs{\eta}} \xi_{1} - \frac{\eta_{1}}{\abs{\eta}} \xi_{j}
\sim L^{2},
$$
whence
$$
\Fourier \left\{ Q_{1j}(D^{-1}u_{L},v_{L}) \right\} \sim L^{2}
\Fourier(u_{L} v_{L}).
$$
It is also easy to see that
$$
\Fourier \left\{ \Lambda^{\half}( \Lambda^{-\half} \Lambda_{-}^{\half}
u_{L} \cdot \Lambda^{\half} v_{L}) \right\} \sim L^{2} \Fourier(u_{L}
v_{L}),
$$
so we have established \eqref{FourierSize}. This concludes the proof. %
\subsection{Proof of Theorem \ref{WMProductEstimates}} %
We first prove the following.
\begin{proposition}\label{TrivialProductEstimates} Let $n \ge 1$ and
$a,b,c,\alpha,\beta,\gamma \ge 0$. Then $$
H^{a,\alpha} \cdot H^{b,\beta} \hookrightarrow H^{-c,-\gamma}, $$
provided $a+b+c > \frac{n}{2}$ and $\alpha + \beta + \gamma > \half$.
\end{proposition}
\begin{proof}
Assume that $s = a+b+c > \frac{n}{2}$ and $\theta = \alpha + \beta +
\gamma > \half$ ($a,b,c,\alpha,\beta,\gamma \ge 0$).

By H\"older's inequality,
\begin{align}
\notag
L^2 \cdot L^\infty &\hookrightarrow L^2, \\
\notag
L_t^\infty(L_x^2) \cdot L_t^2(L_x^\infty) &\hookrightarrow L^2, \\
\intertext{so by \eqref{InftyEmbedding},\eqref{EnergyEmbedding} and
\eqref{SpaceInftyEmbedding} , we have}
\label{ThmHa1}
L^2 \cdot H^{s,\theta} &\hookrightarrow L^2, \\
\label{ThmHa2}
H^{0,\theta} \cdot H^{s,0} &\hookrightarrow L^2. \end{align}
Once we have these estimates, the others follow by interpolation: %
\paragraph{Step 1.}
Assume $\gamma = 0$ (so $\theta = \alpha + \beta > \half$). Interpolation
between \eqref{ThmHa1} and \eqref{ThmHa2} gives $$
H^{0,\alpha} \cdot H^{s,\beta} \hookrightarrow L^2. $$
\paragraph{Step 2.}
By Step 1, we have
\begin{align*}
H^{0,\alpha+\gamma} \cdot H^{s,\beta} &\hookrightarrow L^2, \\
L^{2} \cdot H^{s,\beta} &\hookrightarrow H^{0,-\alpha-\gamma}. \\
\intertext{Interpolation between these yields}
H^{0,\alpha} \cdot H^{s,\beta} &\hookrightarrow H^{0,-\gamma}. \end{align*}
\paragraph{Step 3.}
Assume $c = 0$ (so $s = a + b > \frac{n}{2}$). By Step 2,
\begin{align*}
H^{0,\alpha} \cdot H^{s,\beta} &\hookrightarrow H^{0,-\gamma}, \\
H^{s,\alpha} \cdot H^{0,\beta} &\hookrightarrow H^{0,-\gamma}, \\
\intertext{and by interpolation,}
H^{a,\alpha} \cdot H^{b,\beta} &\hookrightarrow H^{0,-\gamma}. \end{align*}
\paragraph{Step 4.}
By Step 3,
\begin{align*}
H^{a+c,\alpha} \cdot H^{b,\beta} &\hookrightarrow H^{0,-\gamma}, \\
H^{0,\alpha} \cdot H^{b,\beta} &\hookrightarrow H^{-a-c,-\gamma}, \\
\intertext{so interpolation gives}
H^{a,\alpha} \cdot H^{b,\beta} &\hookrightarrow H^{-c,-\gamma}.
\end{align*}
This concludes the proof.
\end{proof}
We now turn to the proof of Theorem \ref{WMProductEstimates}, which we
restate here for convenience.
\begin{plaintheorem}
Let $n \ge 2$, $s > \frac{n}{2}$ and $\half < \theta \le s -
\frac{n-1}{2}$. Then $$
H^{a,\alpha} \cdot H^{s,\theta} \hookrightarrow H^{a,\alpha} $$
for all $a, \alpha$ satisfying
\begin{align*}
0 &\le \alpha \le \theta,
\\
-s + \alpha &< a \le s.
\end{align*}
(Hence, by duality, for all $- \theta \le \alpha \le 0$ and $-s \le a < s
+ \alpha$.) \end{plaintheorem}
By interpolation, it suffices to prove:
\begin{align}
\tag{A}
H^{s,\theta} \cdot H^{s,\theta} &\hookrightarrow H^{s,\theta}, \\
\tag{B}
H^{s,0} \cdot H^{s,\theta} &\hookrightarrow H^{s,0}, \\
\tag{C}
H^{-s,0} \cdot H^{s,\theta} &\hookrightarrow H^{-s,0}, \\
\tag{D}
H^{a,\theta} \cdot H^{s,\theta} &\hookrightarrow H^{a,\theta}, \qquad -s +
\theta < a < 0.
\end{align}
The estimates A,B,C and D correspond to the vertices of a trapezoid in the
$(a,\alpha)$-plane.
\paragraph{Estimate A.}
By Lemma \ref{EllipticLambdaEstimates} it suffices to prove \begin{align*}
H^{0,\theta} \cdot H^{s,\theta} &\hookrightarrow H^{0,\theta}. \\
\intertext{By Lemma \ref{HyperbolicLambdaEstimates}, the last estimate
reduces to three estimates:}
H^{0,\theta} \cdot H^{s,0} &\hookrightarrow L^2, \\
L^2 \cdot H^{s,\theta} &\hookrightarrow L^2, \\
R^\theta( H^{0,\theta}, H^{s,\theta} ) &\hookrightarrow L^2. \end{align*}
The first two are covered by Proposition \ref{TrivialProductEstimates};
the third follows from Theorem F. %
\paragraph{Estimates B and C.}
These are equivalent by duality, so it suffices to prove B, which by lemma
\ref{EllipticLambdaEstimates} reduces to: \begin{align*}
H^{s,0} \cdot H^{0,\theta} &\hookrightarrow L^{2}, \\
L^{2} \cdot H^{s,\theta} &\hookrightarrow L^{2}. \end{align*}
Both of these are covered by Proposition \ref{TrivialProductEstimates}. %
\paragraph{Estimate D.}
Since D is equivalent to
\begin{align}
\notag
H^{-a,-\theta} \cdot H^{s,\theta} &\hookrightarrow H^{-a,-\theta} \\
\intertext{by duality, and since $a < 0$, by Lemma
\ref{EllipticLambdaEstimates} it suffices to prove}
\notag
H^{-a,-\theta} \cdot H^{s+a,\theta} &\hookrightarrow H^{0,-\theta}, \\
\notag
H^{0,-\theta} \cdot H^{s,\theta} &\hookrightarrow H^{0,-\theta}. \\
\intertext{By duality, the last two estimates are equivalent to}
\notag
H^{0,\theta} \cdot H^{s+a,\theta} &\hookrightarrow H^{a,\theta}, \\
\notag
H^{0,\theta} \cdot H^{s,\theta} &\hookrightarrow H^{0,\theta}. \\
\intertext{The last estimate was proved above (estimate for A), and the
second to last reduces, by Lemma \ref{HyperbolicLambdaEstimates}, to three
estimates:}
\label{DD}
H^{0,\theta} \cdot H^{s+a,0} &\hookrightarrow H^{a,0}, \\
\label{DDD}
L^2 \cdot H^{s+a,\theta} &\hookrightarrow H^{a,0}, \\
\label{DDDD}
R^\theta( H^{0,\theta}, H^{s+a,\theta} ) &\hookrightarrow H^{a,0}. \\
\intertext{By interpolation between the estimates}
\notag
H^{0,\theta} \cdot H^{s,0} &\hookrightarrow L^2, \\
\notag
H^{0,\theta} \cdot L^2 &\hookrightarrow H^{-s,0}, \end{align}
which are dual to each other and hold by Proposition
\ref{TrivialProductEstimates}, we get \eqref{DD}. Proposition
\ref{TrivialProductEstimates} also covers \eqref{DDD} (via duality).
Finally, for \eqref{DDDD} we consider two cases: \begin{enumerate}
\item If $s = \frac{n-1}{2} + \theta$, then $a > - \frac{n-1}{2}$, and
\eqref{DDDD} holds by Theorem F.
\item If $s > \frac{n-1}{2} + \theta$, then $-s + \theta < -
\frac{n-1}{2}$, so we may assume that $- s + \theta < a < - \frac{n-1}{2}$
(then the estimate for $- \frac{n-1}{2} \le a < 0$ follows by
interpolation with estimate A). Choose $\varepsilon > 0$ so small that
$\theta + \varepsilon < s + a$ and $\varepsilon \le \frac{n-1}{2}$. Then
by Theorem F,
$$
R^\theta( H^{0,\theta}, H^{\theta + \varepsilon,\theta} ) \hookrightarrow
H^{-\frac{n-1}{2} + \varepsilon,0},
$$
which implies \eqref{DDDD}.
\end{enumerate}
This concludes the proof of Theorem \ref{WMProductEstimates}. %
\subsection{Proof of Theorem \ref{YMProductTheorem}} %
Let $n \ge 4$, $\theta > \half$. Assume
\begin{align*}
a,b &\ge -c,
\\
a + b &\ge \half,
\\
a + b + c &\ge \frac{n-1}{2}.
\end{align*}
We want to prove
\begin{equation}\label{ThmHA}
H^{a,\theta} \cdot H^{b,\theta} \hookrightarrow H^{-c,0}. \end{equation}
\paragraph{Step 1.}
Assume $c \le 0$. Then by Lemma \ref{EllipticLambdaEstimates},
\eqref{ThmHA} reduces to
\begin{align}
\notag
H^{a+c,\theta} \cdot H^{b,\theta} &\hookrightarrow L^2. \\
\intertext{This can be reduced to the extreme case}
\label{ThmHB}
H^{0,\theta} \cdot H^{\frac{n-1}{2},\theta} &\hookrightarrow L^2,
\end{align}
which holds by Theorem F.
\paragraph{Step 2.}
Assume $- c < 0 \le a,b$. If $a + b \ge \frac{n-1}{2}$, then \eqref{ThmHA}
follows from
\begin{align*}
H^{a,\theta} \cdot H^{b,\theta} &\hookrightarrow L^2, \\
\intertext{which again reduces to \eqref{ThmHB}. If $a + b <
\frac{n-1}{2}$, set $\gamma = a + b - \frac{n-1}{2}$. Then \eqref{ThmHA}
follows from}
H^{a,\theta} \cdot H^{b,\theta} &\hookrightarrow H^{\gamma,0}, \end{align*}
and the latter holds by Theorem F.
\paragraph{Step 3.}
Assume $- c \le a < 0$. By Lemma \ref{EllipticLambdaEstimates},
\eqref{ThmHA} reduces to two estimates:
\begin{align*}
H^{0,\theta} \cdot H^{a+b,\theta} &\hookrightarrow H^{-c,0}, \\
H^{0,\theta} \cdot H^{b,\theta} &\hookrightarrow H^{-a-c,0}. \end{align*}
These estimates hold by Steps 1 and 2.
\subsection{Proof of Theorem \ref{CFWMProductTheorem}} %
Let $n \ge 3$, $\theta > \half$. Assume
\begin{align*}
a,b,c &\ge 0,
\\
c &< \frac{n-1}{2},
\\
a + b + c &\ge \frac{n-1}{2} + \theta.
\end{align*}
We want to prove
\begin{align*}
H^{a,\theta} \cdot H^{b,\theta} &\hookrightarrow H^{-c,\theta}. \\
\intertext{By Lemma \ref{HyperbolicLambdaEstimates}, this reduces to}
H^{a,0} \cdot H^{b,\theta} &\hookrightarrow H^{-c,0}, \\
H^{a,\theta} \cdot H^{b,0} &\hookrightarrow H^{-c,0}, \\
R^{\theta}( H^{a,\theta}, H^{b,\theta} ) &\hookrightarrow H^{-c,0}.
\end{align*}
The first two hold by Proposition \ref{TrivialProductEstimates}, the last
one by Theorem F (take $\gamma = -c$, $\gamma_{-} = \theta$ and choose $0
\le s_{1} \le a$, $0 \le s_{2} \le b$ such that $c + s_{1} + s_{2} =
\frac{n-1}{2} + \theta$). %
\subsection{Analysis of the First Iterate} %
Here we work out in more detail the examples considered in section
\ref{Motivation}.
\paragraph{Step 1.}
If $u$ solves $\square u = F$ with vanishing initial data at time $t = 0$,
then
\begin{align*}
\bigabs{\widehat{u(t)}(\xi)} &\le \frac{C_{t}}{\abs{\xi}} \int_\R
\frac{\bigabs{\widehat{F}(\tau,\xi)}} {1 + \hypwt{\tau}{\xi}} \, d\tau,
\\
\bigabs{\widehat{u(t)}(\xi)} &\le t^2
\int_\R \bigabs{\widehat{F}(\tau,\xi)} \, d\tau, \end{align*}
for all $t > 0$.
The first estimate is an immediate consequence of the formula $$
\widehat{u(t)}(\xi)
= \int_{\R} \frac{\widehat{F}(\tau,\xi)}{4 \pi \abs{\xi}} \left(
\frac{e^{it\tau}-e^{it\abs{\xi}}}{\tau - \abs{\xi}} -
\frac{e^{it\tau}-e^{-it\abs{\xi}}}{\tau + \abs{\xi}} \right) \, d\tau,
$$
which is easily derived from Duhamel's principle (see, e.g., \cite[Section
3.6.3]{Se}). As for the second estimate, Duhamel's principle implies $$
\bigabs{\widehat{u(t)}(\xi)} \le t \int_0^t \bigabs{ \widehat{F(t')}(\xi)
} \, dt', $$
and clearly,
$\bigabs{ \widehat{F(t')}(\xi) } \le \int_\R
\bigabs{\widehat{F}(\tau,\xi)} \, d\tau$. %
\paragraph{Step 2.}
Let $B$ be a bilinear operator of the form $$
\Fourier\left(B(u,v)\right)(\tau,\xi) = \int_{\R^{1+n}}
b(\tau-\lambda,\xi-\eta;\lambda,\eta)
\widehat{u}(\tau-\lambda,\xi-\eta) \widehat{v}(\lambda,\eta) \, d\lambda
\, d\eta.
$$
Assume that
$$
\square v = \square w = 0, \qquad (v,\partial_{t} v) \init = (v_{0},0),
\qquad (w,\partial_{t} w) \init = (w_{0},0).
$$
As in section \ref{Notation}, we decompose $v$ and $w$ into half-waves.
Thus,
$$
B(v,w) \simeq B(v_{+},w_{+}) + B(v_{+},w_{-}) + B(v_{-},w_{+}) +
B(v_{-},w_{-}),
$$
where $v_{\pm} = e^{\pm itD} v_{0}$ and $w_{\pm} = e^{\pm itD} w_{0}$. It
suffices to consider the first two terms on the right hand side. Since
$\widehat{v_{+}}(\tau,\xi) \simeq \delta(\tau-\abs{\xi})
\widehat{v_{0}}(\xi)$ and
$\widehat{w_{\pm}}(\tau,\xi) \simeq \delta(\tau\mp\abs{\xi})
\widehat{w_{0}}(\xi)$, we have
\begin{align*}
&\Fourier B(v_{+},w_{\pm}) (\tau,\xi)
\\
&= \int b(\tau-\lambda,\xi-\eta;\lambda,\eta) \delta(\tau - \lambda -
\abs{\xi-\eta})
\delta(\lambda \mp \abs{\eta})
\widehat{v_{0}}(\xi-\eta)
\widehat{w_{0}}(\eta) \, d\lambda \, d\eta \\
&= \int k_{\pm}(\xi-\eta,\eta)
\delta(\tau \mp \abs{\eta} - \abs{\xi-\eta}) \widehat{v_{0}}(\xi-\eta)
\widehat{w_{0}}(\eta) \, d\lambda \, d\eta. \end{align*}
where $k_\pm(\xi,\eta) = b(\abs{\xi},\xi;\pm \abs{\eta},\eta)$. %
\paragraph{Step 3.}
Let $u_{\pm}$ be the solution of $\square u_{\pm} = B(v_{+},w_{\pm})$ with
vanishing initial data. Set $f(\xi) = \elwt{\xi}^{s}
\abs{\widehat{v_{0}}(\xi)}$ and $g(\xi) = \elwt{\xi}^{s}
\abs{\widehat{w_{0}}(\xi)}$. By Steps 1 and 2,
$$
\elwt{\xi}^{s} \bigabs{\widehat{u_{\pm}(t)}(\xi)} \le C_{t} \int
K(\xi-\eta,\eta) f(\xi-\eta) g(\eta) \, d\eta, $$
where
$$
K(\xi,\eta) = \frac{ \elwt{\xi+\eta}^{s} \abs{k_{\pm}}(\xi,\eta)}{
\elwt{\xi}^{s} \elwt{\eta}^{s} } \min \left( 1 , \frac{1}{\abs{\xi+\eta}
\bigl( 1 + \Delta_{\pm}(\xi,\eta) \bigr)} \right)
$$
and $\Delta_{\pm}$ is defined by \eqref{DeltaPlusMinus}. %
\paragraph{Step 4.}
By Step 3, proving $\Sobnorm{u_{\pm}(t)}{s} \le C_{t} \Sobnorm{v_{0}}{s}
\Sobnorm{w_{0}}{s}$ reduces to proving
\begin{equation}\label{AppendixIntegralEstimate}
\int_{\R^{n} \times \R^{n}} K(\xi,\eta) f(\xi) g(\eta) h(\xi + \eta) \,
d\xi \, d\eta
\lesssim \twonorm{f}{} \twonorm{g}{} \twonorm{h}{} \end{equation}
for all $f,g, h \in L^{2}(\R^{n})$, where $K$ is as in Step 3. Write $K =
K_1 + K_2 + K_3$, where $K_1$, $K_2$ and $K_3$ are supported in the
mutually disjoint regions
\begin{align*}
\Omega_1 &= \{ (\xi,\eta) : \abs{\xi+\eta} < 1 \}, \\
\Omega_2 &= \{ (\xi,\eta) : \abs{\xi+\eta} \ge 1, \abs{\eta} < \abs{\xi}
\}, \\
\Omega_3 &= \{ (\xi,\eta) : \abs{\xi+\eta} \ge 1, \abs{\eta} \ge \abs{\xi}
\} \end{align*}
respectively.

Obviously, $K_{1} \lesssim \abs{k_{\pm}} \elwt{\xi}^{-s}
\elwt{\eta}^{-s}$. Assuming that
\begin{equation}\label{AppA}
\abs{k_{\pm}} \lesssim \elwt{\xi}^{s} \elwt{\eta}^{s}, \end{equation}
it follows that $K_{1}$ is bounded, whence $$
\sup_{\eta} \int K_{1}^{2}(\xi-\eta,\eta) \, d\xi < \infty. $$
Therefore \eqref{AppendixIntegralEstimate} holds by Lemma
\ref{IntegralOperatorLemma} below (after a linear change of variables).
Observe that \eqref{AppA} is satisfied in the examples we consider (we
always have $s \ge 1$ and $\abs{k_{\pm}} \lesssim \abs{\xi}\abs{\eta}$).

Next, for $K_{2}$ we have
$$
K_{2}(\xi,\eta) \lesssim \frac{\abs{k_{\pm}}(\xi,\eta)}{ \elwt{\xi}
\elwt{\eta}^{s} \bigl( 1 + \Delta_{\pm}(\xi,\eta) \bigr)} $$
Let us now consider this expression for the operators $B(v,w)$ appearing
in our examples (the estimates for $K_{3}$ are the same, since the
operators are symmetric).
\begin{enumerate}
\item In Example \ref{FirstIterateExample1}, $B(v,w) = \partial_{t} v
\cdot \partial_{t} w$, so $b(\tau,\xi; \lambda,\eta) \simeq \tau \lambda$.
Therefore,
$$
\abs{k_{\pm}} \le \abs{\xi} \abs{\eta}
$$
and
$$
K_{2}(\xi,\eta) \lesssim \frac{1}{ \elwt{\eta}^{s-1} \bigl( 1 +
\Delta_{\pm}(\xi,\eta) \bigr)}. $$
\item
In Example \ref{FirstIterateExample2}, $B(v,w) = Q_{0}(v,w)$, so
$b(\tau,\xi; \lambda,\eta) \simeq \tau \lambda - \xi \cdot \eta$.
Therefore,
$$
k_{\pm}(\xi,\eta) \simeq
\pm \abs{\xi} \abs{\eta}
- \xi \cdot \eta
= \pm \abs{\xi} \abs{\eta} \left( 1 \mp
\cos \angle(\xi,\eta) \right).
$$
And in view of Lemma \ref{DeltaEstimate} below, this implies $$
\abs{k_{\pm}} \le \max(\abs{\xi},\abs{\eta}) \Delta_{\pm}(\xi,\eta). $$
Hence,
$$
K_{2}(\xi,\eta) \lesssim \elwt{\eta}^{-s}. $$
\item
In Example \ref{FirstIterateExample3}, $B(v,w) = Q_{ij}(v,w)$, so
$b(\tau,\xi; \lambda,\eta) \simeq \xi_{i} \eta_{j} - \xi_{j} \eta_{i}$.
Therefore,
$$
\abs{k_{\pm}(\xi,\eta)}
\lesssim \abs{\xi \wedge \eta}
= \abs{\xi} \abs{\eta} \sqrt{1-\cos^{2} \theta} \lesssim \abs{\xi}
\abs{\eta} \sqrt{1 \mp \cos \theta}.
$$
In view of Lemma \ref{DeltaEstimate}, this implies
$$
K_{2}(\xi,\eta) \lesssim \frac{1}
{\elwt{\eta}^{s-\half} \bigl( 1 + \Delta_{\pm}(\xi,\eta) \bigr)^{\half}}.
$$
\end{enumerate}
\subsubsection{Proof of Proposition \ref{IntegralEstimateProposition}} %
Instead of proving Proposition \ref{IntegralEstimateProposition} as
stated, we prove a homogeneous version. The proof is easily modified to
give the inhomogeneous statement in Proposition
\ref{IntegralEstimateProposition}. We first prove two lemmas. %
\begin{lemma}\label{DeltaEstimate}
Let $\Delta_{\pm}$ be defined by \eqref{DeltaPlusMinus}. Then $$
\min(\abs{\xi},\abs{\eta}) (1 \pm \cos \theta) \le 2
\Delta_{\mp}(\xi,\eta),
$$
where $\theta$ is the angle between $\xi$ and $\eta$. \end{lemma}
\begin{proof}
We have
$$
\Delta_{+}(\xi,\eta)
= \frac{\Delta_{+}(\xi,\eta)(\abs{\xi}+\abs{\eta} + \abs{\xi + \eta})}
{\abs{\xi}+\abs{\eta} + \abs{\xi + \eta}} \ge \frac{\abs{\xi}\abs{\eta}(1
- \cos \theta)}{\abs{\xi} + \abs{\eta}}, $$
and a similar computation gives the proof for $\Delta_{-}$. \end{proof}
\begin{lemma}\label{IntegralOperatorLemma} If $K$ is a measurable function
on
$\R^{n} \times \R^{n}$ such that at least one of the numbers
$$
\sup_{\xi} \int K^{2}(\xi,\eta) \, d\eta, \quad \sup_{\eta} \int
K^{2}(\xi,\eta) \, d\xi $$
is finite, then
$$
\int_{\R^{n} \times \R^{n}}
K(\xi,\eta) f(\xi) g(\eta) h(\xi + \eta) \, d\xi \, d\eta \le
C\twonorm{f}{} \twonorm{g}{} \twonorm{h}{} $$
for all $f,g,h \in L^{2}$.
\end{lemma}
To prove this, simply apply the Cauchy-Schwarz inequality twice. %
\begin{plainproposition}
If $a,b,c \ge 0$ and $\Delta(\xi,\eta)$ is either of the expressions
defined in \eqref{DeltaPlusMinus}, then
$$
\int_{\R^{n} \times \R^{n}} \frac{f(\xi) g(\eta) h(\xi+\eta)} {
\abs{\xi}^{a} \abs{\eta}^{b}
\Delta^{c}(\xi,\eta)} \, d\xi d\eta
\le C \twonorm{f}{} \twonorm{g}{} \twonorm{h}{} $$
for all $f,g,h \in L^{2}(\R^{n})$, provided $$
a + b + c = \frac{n}{2}, \quad a,b < \frac{n}{2} - c, \quad c <
\frac{n-1}{4}.
$$
\end{plainproposition}
\begin{proof}
Set $K(\xi,\eta) = \abs{\xi}^{-a} \abs{\eta}^{-b} \Delta^{-c}(\xi,\eta)$
and write $K = K_{1} + K_{2}$, where $K_{1}$ is supported in the region
$\abs{\eta} \le \abs{\xi}$ and $K_{2}$ is supported in $\abs{\eta} >
\abs{\xi}$. By lemma \ref{DeltaEstimate},
$$
K_{1}(\xi,\eta) \le \frac{2}{\abs{\xi}^{a} \abs{\eta}^{b+c} (1 \pm \cos
\theta)^{c}},
$$
where $\theta$ is the angle between $\xi$ and $\eta$. Thus, for all $\xi$,
integration in polar coordinates $(r,\omega) =
(\abs{\eta},\eta/\abs{\eta})$ yields
\begin{align*}
\int K_{1}^{2}(\xi,\eta) \, d\eta
&\le \frac{4}{\abs{\xi}^{2a}}
\int_{0}^{\abs{\xi}} r^{n-1-2(b+c)} \, dr \int_{S^{n-1}}
\frac{d\sigma(\omega)}{(1\pm\cos\theta)^{2c}} \\
&= \frac{4}{n-2(b+c)} \int_{S^{n-1}}
\frac{d\sigma(\omega)}{(1\pm\cos\theta)^{2c}}, \end{align*}
and the last integral is finite iff $2c < (n-1)/2$. By symmetry, this
implies that $\sup_{\eta} \int K_{2}^{2}(\xi,\eta) \, d\xi$ is also
finite, so we may apply lemma \ref{IntegralOperatorLemma}. \end{proof}
\end{document}